\documentclass[12pt]{article}

\usepackage[normalem]{ulem}

\makeatletter
\newcommand*\nobreakhyphen{\hbox{-}\nobreak\hskip\z@skip}
\makeatother

\DeclareFontFamily{U}{matha}{}
\DeclareFontShape{U}{matha}{m}{n}{
	<-5.5>    matha5
	<5.5-6.5> matha6
	<6.5-7.5> matha7
	<7.5-8.5> matha8
	<8.5-9.5> matha9
	<9.5-11>  matha10
	<11->     matha12
}{}
\DeclareSymbolFont{matha}{U}{matha}{m}{n}
\DeclareFontSubstitution{U}{matha}{m}{n}
\DeclareFontFamily{U}{mathx}{\hyphenchar\font45}
\DeclareFontShape{U}{mathx}{m}{n}{<-> mathx10}{}
\DeclareSymbolFont{mathx}{U}{mathx}{m}{n}
\DeclareFontSubstitution{U}{mathx}{m}{n}

\DeclareMathAccent{\widecheck}{0}{mathx}{"71}

\usepackage{tikz}
\usepackage{graphicx}
\usepackage{color}

\usepackage{amssymb,amsmath}
\usepackage{amsfonts}
\usepackage{comment}
\usepackage{graphicx}
\usepackage{wasysym}
\usepackage{enumerate}

\usepackage[margin=1.5cm]{geometry}

\newtheorem{thm}{Theorem}

\newtheorem{lem}{Lemma}
\newtheorem{exa}{Example}

\newcommand{\bpi}{\mbox{\boldmath$\pi$}}

\newcommand{\CC}{\mbox{$\mathbb C$}}

\newcommand{\bmu}{\mbox{\boldmath$\mu$}}
\newcommand{\bnu}{\mbox{\boldmath$\nu$}}

\newtheorem{Remark}{Remark}

\title{Matrix-Analytic Methods for the analysis of \\Stochastic Fluid-Fluid Models}

\author{
Nigel G. Bean\thanks{Australian Research Centre of Excellence for Mathematical and Statistical Frontiers. School of Mathematical Sciences, University of Adelaide, SA 5005, Australia (nigel.bean@adelaide.edu.au).}
\and
Ma{\l}gorzata M. O'Reilly\thanks{Australian Research Centre of Excellence for Mathematical and Statistical Frontiers. Discipline of Mathematics, University of Tasmania, TAS 7005, Australia (malgorzata.oreilly@utas.edu.au).}
\thanks{Ma{\l}gorzata M. O'Reilly would like to thank the Australian Research Council for funding this research through Linkage Project LP140100152.}
\and
Zbigniew Palmowski\thanks{Faculty of Pure and Applied Mathematics, Wroc{\l}aw University of Science and Technology, ul. Wybrze\.{z}e Wyspia\'{n}skiego 27,
	50-370 Wroc{\l}aw, Poland (zbigniew.palmowski@pwr.edu.pl).}
\thanks{This work is partially supported by the National Science Centre under the grant 2018/29/B/ST1/00756 (2019-2022).}
}

\date{\normalsize \today}

\begin{document}
\maketitle
\begin{abstract}
Stochastic fluid-fluid models (SFFMs) offer powerful modeling ability for a wide range of real-life systems of significance. The existing theoretical framework for this class of models is in terms of operator-analytic methods. For the first time, we establish matrix-analytic methods for the efficient  analysis of SFFMs. We illustrate the theory with numerical examples.
\end{abstract}


\section{Introduction}\label{S_intro}

Many real-life systems have an element of uncertainty which is modelled by probabilistic methods.
Inspired by engineering problems, economical theories, telecommunication networks, health care, insurance and manufacturing/management systems, environmental problems or biological applications, in this paper
we consider stochastic fluid-fluid models (SFFMs) introduced by Bean and O'Reilly in~\cite{BeaOreFF}. The main goal of this paper is to establish matrix-analytic methods for the efficient numerical analysis of this class of models.

The SFFMs is an extension of the stochastic fluid models (SFMs), a class of models proposed by Anick, Mitra and Sondhi in~\cite{anms82} and analysed in a series of papers on the stationary and transient behaviour~\cite{AJR,AR,AR2,2020AHT,asmu95,BOT4,BOT,BOT3,rama96,rama99,1994Rog,SOB17}, algorithms for the computation of key measures~\cite{AR3,BOT2,BOT5,2017GUM}, with great application potential in areas that include manufacturing~\cite{anms82}, risk processes in insurance~\cite{BBD2}, hydro power generation systems~\cite{Bea-Ore}, distribution of resource in peer-to-peer file sharing applications~\cite{2007GT}, maintenance in continuously deteriorating systems~\cite{2017SHOB}, SIR epidemics~\cite{2020S}, and queues with abandonment in passenger-taxi service systems and organ transplantation systems~\cite{2020WQ}.

\begin{figure}
	\begin{center}
		\begin{tikzpicture}[>=stealth,redarr/.style={->}]
		\filldraw[draw=black,fill=lightgray] (4,4) rectangle (6,6);
		\draw[draw=black,fill=white] (4,6) rectangle (6,7);
		\draw (3.2,4.3) node[anchor=north, below=-0.17cm] {\footnotesize{\color{black} $X$, $c_i$}};
		
		\draw (6.8,6.2) node[anchor=north, below=-0.17cm] {\scriptsize{\color{black} $x$}};
		\draw (6.8,4) node[anchor=north, below=-0.17cm] {\scriptsize{\color{black} $0$}};
		\draw [ ->] (7,4) -- (12,4);
		\draw [ ->] (7,4) -- (7,7);
		\draw (6.6,7) node[anchor=north, below=-0.17cm] {\scriptsize{\color{black} $X(t)$}};
		\draw (11.8,4.2) node[anchor=north, below=-0.17cm] {\scriptsize{\color{black} $t$}};
		
		\draw [blue, thick, ->] (7,6.2) -- (7.5,6) -- (8,6.5) -- (8.5,5) -- (9,5.5) -- (9.5,4) -- (10,4) -- (10.5,4.5);
		
		\node at (7.5,6) [red,circle,fill,inner sep=1.3pt]{};
		\node at (8,6.5) [red,circle,fill,inner sep=1.3pt]{};
		\node at (8.5,5) [red,circle,fill,inner sep=1.3pt]{};		
		\node at (9,5.5) [red,circle,fill,inner sep=1.3pt]{};
		\node at (10,4) [red,circle,fill,inner sep=1.3pt]{};
		
		\node at (10,6.5) [red,circle,fill,inner sep=1.3pt]{\color{white}{\bf 1}};
		\node at (11,6.5) [red,circle,fill,inner sep=1.3pt]{\color{white}{\bf 2}};
		\node at (12,6.5) [red,circle,fill,inner sep=1.8pt]{\color{white}{\bf n}};
		
		\draw [black,->] (10,6.7) to [out=90,in=90] (11,6.7);
		\draw [black,->] (11,6.7) to [out=90,in=90] (12,6.7);
		\draw [black,->] (10,6.7) to [out=90,in=90] (12,6.7);
		
		\draw [black,->] (11,6.3) to [out=270,in=270] (10,6.3);	
		\draw [black,->] (12,6.3) to [out=270,in=270] (11,6.3);
		\draw [black,->] (12,6.3) to [out=270,in=270] (10,6.3);
		
		\draw (11.5,6.5) node[anchor=north, below=-0.17cm] {\large{\color{black} $\ldots$}};
		
		\draw (9.4,7.1) node[anchor=north, below=-0.17cm] {{\color{black} \footnotesize $\varphi(t)$}};

		\node at (11,5) [black,circle,fill,inner sep=1.3pt]{};
		\draw [black, ->] (11,5) -- (11.5,5);
		\draw [black, ->] (11,5) -- (11.5,5.5);
		\draw [black, ->] (11,5) -- (11.5,4.5);
		\end{tikzpicture}
		\caption{\footnotesize Evolution of an SFM $\{(\varphi(t),X(t)):t\geq 0\}$ modelling the level in buffer $X$. The process is driven by a continuous-time Markov chain $\{\varphi(t):t\geq 0\}$ with state space $\mathcal{S}=\{1,\ldots,n\}$ and some rates $c_i$, $i\in\mathcal{S}$ so that $dX(t)/dt=c_{\varphi(t)}\times I(X(t)>0)+ \max \{c_{\varphi(t)},0\}\times I(X(t)=0)$, where $I(\cdot)$ is the indicator function. Possible directions of movement are indicated by the black arrows. }		
		\label{fig:1dSFM}
	\end{center}
\end{figure}
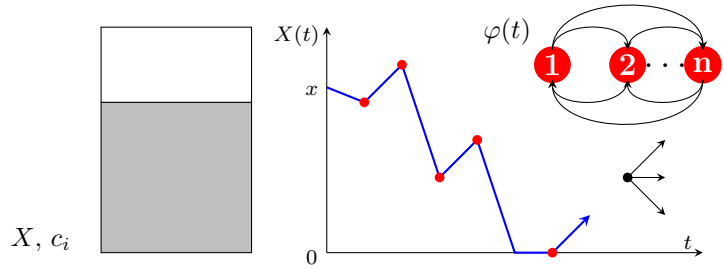

Another useful illustration of an application of the SFMs is a water reservoir~\cite{BeaOreFF}, whose level $X(t)$ at time~$t$ is changing according to the fluid rate $c_i$ driven by the phase $\varphi(t)=i$ of an underlying environment (modelling water usage as well as supply), and so the level in the reservoir may increase, decrease, or remain constant, see Figure~\ref{fig:1dSFM}.

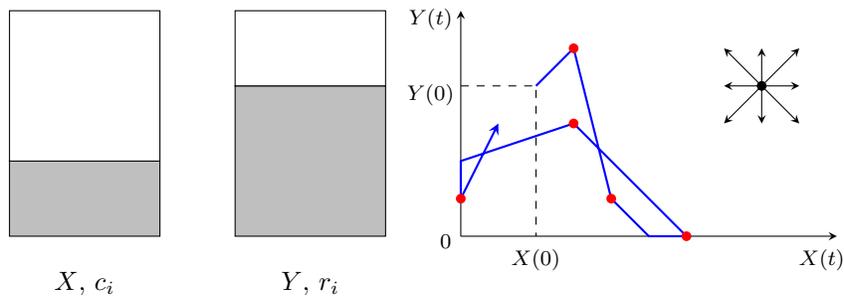
\begin{figure}
	\begin{center}
		\begin{tikzpicture}[>=stealth,redarr/.style={->}]

		\filldraw[draw=black,fill=lightgray] (1,4) rectangle (3,5);
		\draw[draw=black,fill=white] (1,5) rectangle (3,7);
		\draw (2,3.5) node[anchor=north, below=-0.17cm] {\footnotesize{\color{black} $X$, $c_i$}};
		
		\filldraw[draw=black,fill=lightgray] (4,4) rectangle (6,6);
		\draw[draw=black,fill=white] (4,6) rectangle (6,7);
		\draw (5,3.5) node[anchor=north, below=-0.17cm] {\footnotesize{\color{black} $Y$, $r_i$}};		
		\draw (6.6,6) node[anchor=north, below=-0.17cm] {\scriptsize{\color{black} $Y(0)$}};
		\draw (8,3.8) node[anchor=north, below=-0.17cm] {\scriptsize{\color{black} $X(0)$}};
		\draw (6.8,4) node[anchor=north, below=-0.17cm] {\scriptsize{\color{black} $0$}};
		\draw [ ->] (7,4) -- (12,4);
		\draw [ ->] (7,4) -- (7,7);
		\draw (6.6,7) node[anchor=north, below=-0.17cm] {\scriptsize{\color{black} $Y(t)$}};
		\draw (11.8,3.8) node[anchor=north, below=-0.17cm] {\scriptsize{\color{black} $X(t)$}};
		
		\draw [blue, thick, ->] (8,6) -- (8.5,6.5) -- (9,4.5) -- (9.5,4) -- (10,4)
		-- (8.5,5.5) -- (7,5) -- (7,4.5) -- (7.5,5.5);
		
		\draw [dashed] (7,6) -- (8,6) -- (8,4);
		\node at (8.5,6.5) [red,circle,fill,inner sep=1.3pt]{};
		\node at (9,4.5) [red,circle,fill,inner sep=1.3pt]{};
		\node at (10,4) [red,circle,fill,inner sep=1.3pt]{};
		\node at (8.5,5.5) [red,circle,fill,inner sep=1.3pt]{};
		\node at (7,4.5) [red,circle,fill,inner sep=1.3pt]{};

		\node at (11,6) [black,circle,fill,inner sep=1.3pt]{};
		\draw [black, ->] (11,6) -- (11.5,6);
		\draw [black, ->] (11,6) -- (10.5,6);
		\draw [black, ->] (11,6) -- (11,6.5);
		\draw [black, ->] (11,6) -- (11,5.5);
		\draw [black, ->] (11,6) -- (11.5,6.5);
		\draw [black, ->] (11,6) -- (10.5,5.5);
		\draw [black, ->] (11,6) -- (11.5,5.5);
		\draw [black, ->] (11,6) -- (10.5,6.5);

		\end{tikzpicture}
	\end{center}
	\caption{\footnotesize Evolution of an SFFM $\{(\varphi(t),X(t),Y(t)):t\geq 0\}$ modelling the level in buffers $X$ and $Y$. The process is driven by a continuous-time Markov chain $\{\varphi(t):t\geq 0\}$ with state space $\mathcal{S}=\{1,\ldots,n\}$ and some rates $c_i$, $r_i$. Various directions of movement on the quadrant may be possible depending on the model parameters, as indicated by the arrows.		
	}
	\label{fig:SFFMquadrant}
\end{figure}

Stochastic fluid-fluid model (SFFM) introduced in~\cite{BeaOreFF} is built on the concept of the Markov Chains. Its state space is three-dimensional and consists of the discrete {\em phase} variable $\varphi(t)\in{\mathcal S}$, where ${\mathcal S}$ is some finite set, and two continuous {\em level} variables $X(t)$
and $Y(t)$. Phase is used to model the state of some real-life environment/system, while levels are used to model some continuous performance measures of the system, such as the content of the buffers $X$ and $Y$ at time $t$.

The SFFMs further extend the application potential of the SFMs to systems such as tandem networks~\cite{2010WR}, process of growth and bleaching of
coral reefs~\cite{BeaOreFF}, resource-sharing in computer processing~\cite{2013LNP}, and in general, any system that can be modelled using an SFM $\{(\varphi(t),X(t)):t\geq 0\}$, with an additional level variable $Y(t)$ introduced to model some continuous performance measure of that system~\cite{BeaOreFF}.

The evolution of the level $Y(t)$ may depend on the value of the level $X(t)$ at time $t$, see an example of a tandem fluid queue analysed by O'Reilly and Scheinhardt in~\cite{2017_OS}. To model such scenarios, Bean and O'Reilly \cite{BeaOreFF} considered an SFFM in which the second fluid, $Y(t)$, is driven by rates $r_i(x)$, $i\in\mathcal{S}$, $x\geq 0$, that depend on $\varphi(t)=i$ and may also depend on $X(t)=x$. In the model considered here we assume constant rates $r_i(x)=r_i$ for all $x$. We illustrate this in Figure~\ref{fig:SFFMquadrant}.

More formally, in this work we consider an SFFM $\{(\varphi(t),X(t),Y(t)):t\geq 0\}$ with a three-dimensional state space which consists of the phase variable $\varphi(t)\in\mathcal{S}$ and two continuous level variables $X(t)\geq 0$, $Y(t)\geq 0$, each bounded from below by $0$, and real-valued, {\em nonzero} rates $c_i$, $r_i$, for all $i\in\mathcal{S}$, such that the following assumptions are met.
\begin{itemize}
\item The phase process $\{(\varphi(t)):t\geq 0\}$ is an irreducible, continuous-time Markov Chain (CTMC) with finite state space $\mathcal{S}$ and generator ${\bf T}=[T_{ij}]_{i,j\in\mathcal{S}}$.
\item When $\varphi(t)=i$ and $X(t)>0$, the rate of change of $X(t)$ at time $t$ is given by $dX(t)/dt=c_i$. When $\varphi(t)=i$ and $X(t)=0$, the rate of change of $X(t)$ at time $t$ is given by $dX(t)/dt=\max\{0,c_i\}$.
\item When $\varphi(t)=i$ and $Y(t)>0$, the rate of change of $Y(t)$ at time $t$ is given by $dY(t)/dt=r_i$. When $\varphi(t)=i$ and $Y(t)=0$, the rate of change of $X(t)$ at time $t$ is given by $dY(t)/dt=\max\{0,r_i\}$.
\end{itemize}
Let $\mathcal{S}^+=\{j\in\mathcal{S}:r_j>0\}$, $\mathcal{S}^-=\{j\in\mathcal{S}:r_j < 0\}$ and similarly, $\mathcal{S}_+=\{j\in\mathcal{S}:c_j>0\}$, $\mathcal{S}_-=\{j\in\mathcal{S}:c_j < 0\}$. We have $\mathcal{S}^0=\{j\in\mathcal{S}:r_j=0\}=\varnothing$, $\mathcal{S}_0=\{j\in\mathcal{S}:c_j=0\}=\varnothing$. That is, $\mathcal{S}=\mathcal{S}^+\cup\mathcal{S}^-
	=\mathcal{S}_+\cup\mathcal{S}_-$.

We refer to levels $X$ and $Y$ as buffers where fluid enters according to intensities $c_i$ and $r_i$, respectively. In order to derive the results, we will also deal with an unbounded level process $\widetilde Y(\cdot)\in(-\infty,+\infty)$ which evolves with intensity $r_i$ when $\varphi(t)=i$, that is,
$d\widetilde{Y}(t)/dt=r_{\varphi(t)}$.

Bean and O'Reilly considered a related two-dimensional stochastic fluid model in~\cite{Bea-Ore2}, in which one of the level variables was assumed to be not bounded from below, and the analysis was performed using efficient matrix-analytic methods~\cite{Lat-Ram}. However, the methods developed for the model in~\cite{Bea-Ore2} may not be directly applied for the analysis of the SFFMs, since here, both fluids are bounded from below, and so instead, alternative techniques are required. Further, the theoretical results for the SFFMs in~\cite{BeaOreFF} were developed using operator-analytic methods, which are not directly applicable to numerical computations using matrices.

In this paper, we address this gap and develop matrix-analytic methods for the analysis of a class of the SFFMs with nonzero rates $r_i$ and $c_i$ under the assumptions outlined below, and derive expressions for the distributions at level-dependent stopping times corresponding to the process $Y(\cdot)$.

{\bf Assumption 1.} In Theorems~\ref{exp_distr}-\ref{th:main} we assume that the initial position $X(0)$ of the level process $X(\cdot)$ follows a distribution which consists of a density at $X(0)=x>0$ which has an exponential form, and a possible point mass at~$X(0)=0$. Our choice of the exponential form for the density is due to the fact that it is a classic distribution, present in many applications, which is convenient to handle and which plays an important role in the statistical theory of reliability and lifetime analysis~\cite{2021YTT}. Furthermore, we note that by letting the intensity of the exponential distribution approach infinity, we obtain a fluid process $X(\cdot)$ which starts from level $0$ with probability~$1$.

{\bf Assumption 2.} In Theorem~\ref{lem_ex} we assume that the absolute values of rates $r_i$ and $c_i$ are proportional to one another according to $|r_i|=\gamma |c_i|$ for some constant $\gamma>0$, and that generator ${\bf T}$  is such that
\begin{eqnarray*}
	|{\bf R}|^{-1}{\bf T}&=&	
	\left[
	\begin{array}{cc}
		-(b+\beta){\bf I}&(|{\bf R}_+|^{-1}{\bf T}_{+-})\\
		(|{\bf R}_{-}|^{-1}{\bf T}_{- +})&-b{\bf I}
	\end{array}
	\right]
\end{eqnarray*}
for some $b,\beta > 0$, where ${\bf T}_{+-}=[T_{ij}]_{i\in\mathcal{S}_+,j\in\mathcal{S}_-}$, ${\bf T}_{-+}=[T_{ij}]_{i\in\mathcal{S}_-,j\in\mathcal{S}_+}$, ${\bf R}=diag(r_j)_{j\in\mathcal{S}}$, ${\bf R}_{+}=diag(r_j)_{j\in\mathcal{S}_+}$, ${\bf R}_{-}=diag(r_j)_{j\in\mathcal{S}_-}$, and so the sign of $c_{\varphi(t)}$ changes at the time of jump in $\varphi(\cdot)$.

Under these assumptions, we give new representation for the Laplace-Stieltjes transform (LST) of the time taken for the total amount of fluid that has flowed into or out of the buffer $\widetilde Y(\cdot)$ to reach some level $y$. This result allows us to find the key quantity $\Psi$ recording the probabilities that the fluid level ${Y}(\cdot)$ first returns to level~$0$ and does so in some phase $j$ and with $X(\cdot)$ contained in some chosen set. Using these results we perform some numerical analysis which we illustrate through several examples.

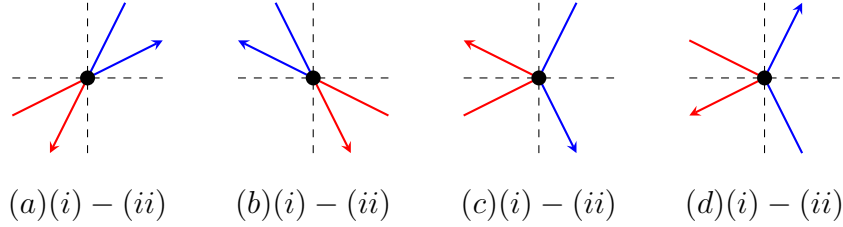
\begin{figure}
	\begin{center}
		\begin{tikzpicture}[>=stealth,redarr/.style={->}]

		\draw [black, dashed] (0,6) -- (2,6);
		\draw [black, dashed] (1,5) -- (1,7);
		\draw [red,thick, ->] (0,5.5) -- (1,6) -- (0.5,5);
		\draw [blue,thick, ->] (1.5,7) -- (1,6) -- (2,6.5);
		\node at (1,6) [black,circle,fill,inner sep=2pt]{};
		\draw (1,4.5) node[black,anchor=north, below=-0.17cm] {$(a)(i)-(ii)$};
		
		\draw [black, dashed] (3,6) -- (5,6);
		\draw [black, dashed] (4,5) -- (4,7);
		\draw [red,thick, ->] (5,5.5) -- (4,6) -- (4.5,5);
		\draw [blue,thick, ->] (3.5,7) -- (4,6) -- (3,6.5);
		\node at (4,6) [black,circle,fill,inner sep=2pt]{};
		\draw (4,4.5) node[black,anchor=north, below=-0.17cm] {$(b)(i)-(ii)$};

		\draw [black, dashed] (6,6) -- (8,6);
		\draw [black, dashed] (7,5) -- (7,7);
		\draw [red,thick, ->] (6,5.5) -- (7,6) -- (6,6.5);
		\draw [blue,thick, ->] (7.5,7) -- (7,6) -- (7.5,5);
		\node at (7,6) [black,circle,fill,inner sep=2pt]{};
		\draw (7,4.5) node[black,anchor=north, below=-0.17cm] {$(c)(i)-(ii)$};

		\draw [black, dashed] (9,6) -- (11,6);
		\draw [black, dashed] (10,5) -- (10,7);
		\draw [red,thick, ->] (9,6.5) -- (10,6) -- (9,5.5);
		\draw [blue,thick, ->]  (10.5,5) -- (10,6) -- (10.5,7);
		\node at (10,6) [black,circle,fill,inner sep=2pt]{};
		\draw (10,4.5) node[black,anchor=north, below=-0.17cm] {$(d)(i)-(ii)$};

		\end{tikzpicture}
	\end{center}
	\caption{\footnotesize Possible changes in the direction in an SFFM $\{(\varphi(t),X(t),Y(t)):t\geq 0\}$ satisfying Assumption 2 when $i\to j$: (a) $i,j\in\mathcal{S}_e$ with (i) $c_i,r_i>0$ and $c_j,r_j<0$ in red and (ii) $c_i,r_i<0$ and $c_j,r_j>0$ in blue; (b) $i,j\notin\mathcal{S}_e$ with (i) $c_i<0,r_i>0$ and $c_j>0,r_j<0$ in red and (ii) $c_i>0,r_i<0$ and $c_j<0,r_j>0$ in blue; (c) $i\in\mathcal{S}_e$, $j\notin\mathcal{S}_e$ with (i) $c_i,r_i>0$ and $c_j<0,r_j<0$ in red and (ii) $c_i,r_i<0$ and $c_j>0,r_j<0$ in blue; (d) $i\notin\mathcal{S}_e$, $j\in\mathcal{S}_e$ with (i) $c_i>0,r_i<0$ and $c_j,r_j<0$ in red and (ii) $c_i<0,r_i>0$ and $c_j,r_j>0$ in blue. }
	\label{fig:SFFMquadrantA1A2}
\end{figure}

Although our key aim is to develop theoretical results, the application potential of this class of SFFMs is also worth emphasizing. To motivate Assumption 2, we consider the following application example. Suppose that a tandem consisting of two telecommunication buffers $X$ and $Y$ driven by an underlying environment $\varphi(\cdot)$ is such that the direction of flow of the fluid is the same in both buffers, with $r_i=\gamma c_i$, or the opposite to one another, with $r_i=-\gamma c_i$, for some $\gamma>0$, depending on whether an underlying phase $\varphi(t)=i$ at time~$t$ is within some specified subset $\mathcal{S}_e=\{i,j\in\mathcal{S}:c_i\times r_i>0\}$ of $\mathcal{S}$ or not, respectively.
\begin{itemize}
	\item As example, suppose that whenenver $\varphi(t)=i\in\mathcal{S}_e$, the data, represented as fluid, enters or leaves such tandem at same rate $d_i$ shared proportionally between the two buffers according to $c_i=d_i/(1+\gamma)$, $r_i= d_i\gamma/(1+\gamma)=\gamma c_i$, with $r_i+c_i=d_i$. Therefore, the direction of flow of the fluid is the same in both buffers, and buffer $Y$ gets a larger, equal, or smaller proportion of the activity depending on whether we set $\gamma>1$, $\gamma=1$, or $\gamma<1$, respectively.
	\item Under assumption $\gamma=1$, suppose that whenever $\varphi(t)=i\notin\mathcal{S}_e$, the fluid level in buffer $X$ changes at rate $c_i$, and in buffer $Y$ in the opposite direction at rate $r_i=-c_i$, with the rate of change of $X+Y$ being $d_i=r_i+c_i=0$ as long as the buffer with the negative rate of change is nonempty. We may interpret this as the fluid moving from one buffer (the one with the negative rate) to another whenever possible.
	\item Under assumption $\gamma\not= 1$, suppose that whenever $\varphi(t)=i\notin\mathcal{S}_e$, the fluid level in buffer $X$ changes at rate $c_i=d_i/(1-\gamma)$ for some $d_i\not= 0$, and in buffer buffer $Y$ in the opposite direction at rate $r_i=-\gamma d_i/(1-\gamma)=-\gamma c_i$. Therefore, the data enters one of the two buffers and leaves the other, and the rate of change of $X+Y$ is $d_i=r_i+c_i$ as long as the buffer with the negative rate of change is nonempty.
	\item Also, assume that the following changes in the direction of movement are possible. Whenever the process transitions from some phase $i$ to $j$ with either $i,j\in\mathcal{S}_e$ or $i,j\notin\mathcal{S}_e$ then the fluid level starts moving in the opposite directions to before in both buffers (the signs of both $r_{\varphi(t)}$ and $c_{\varphi(t)}$ change). However, if the process transitions from some phase $i$ to $j$ with either $i\in\mathcal{S}_e,j\notin\mathcal{S}_e$ or $i\notin\mathcal{S}_e,j\in\mathcal{S}_e$ then the fluid level starts moving in the opposite directions to before only in buffer $X$ (only the sign of $c_{\varphi(t)}$ changes). Potential changes in the directions of the movement are illustrated in Figure~\ref{fig:SFFMquadrantA1A2}.
\end{itemize}
We note that such tandem can be equivalently described as a movement of an individual on a quadrant $\{(x,y):x\geq 0,y\geq 0\}$ with a lot of flexibility in the possible directions of movement, and can be analysed using the class of SFFMs studied here.

The rest of the paper is organized as follows. In the next section we remind key results and notations concerning SFFM queues. Then in Section \ref{sec:mainresult} we state the main result identifying the key generator $D$ of the SFFM. Finally, in Section~\ref{sec:returnzero} we give the formula for the key quantity $\Psi$ corresponding to the first return time to zero in the unbounded level process $\widetilde Y(\cdot)$ as well as in the process ${Y}(\cdot)$. In both these sections we present numerical examples to illustrate the theory as well. We end our paper with conclusions in Section~\ref{sec:conc}.

\section{Preliminaries}\label{Prelim}


In this section we gather some useful results established for the SFMs in Samuelson, O'Reilly and Bean~\cite{SOB17}, and for the SFFMs in Bean and O'Reilly~\cite{BeaOreFF}, with minor notational changes to suit the analysis here. We will build on these results in Sections~\ref{sec:mainresult}--\ref{sec:returnzero}. The key idea in our approach is to
\begin{itemize}
\item consider suitable level-dependent stopping times in paths in the SFM $\{(\varphi(t),Y(t)):t\geq 0\}$; and
\item evaluate expressions for the distribution of the level $X(\cdot)$ at those stopping times.
\end{itemize}

We note that the choice of the level $Y(\cdot)$ versus $X(\cdot)$ is arbitrary. By symmetry, our arguments can be applied to the stopping times in the sample paths in the SFM $\{(\varphi(t),X(t)):t\geq 0\}$ and the distribution of the level $Y(\cdot)$ at those stopping times.

Also, we note that $\mathcal{S}^{\pm}\not=\mathcal{S}_{\pm}$ in general, since the rates $r_i$ and $c_i$ may have different signs, as we illustrate in the numerical examples with a range of behaviours in Sections~\ref{sec:numeex1} and~\ref{sec:numeex2}.

\begin{Remark}
The methodology developed here for the class of SFFMs with $\mathcal{S}^0=\mathcal{S}_0=\varnothing$ can be also applied when $\mathcal{S}^0=\mathcal{S}_0\not=\varnothing$. Suppose an SFFM $\{(\bar\varphi(t),\bar X(t),\bar Y(t)):t\geq 0\}$ driven by a CTMC $\{(\bar\varphi(t)):t\geq 0\}$ with state space $\mathcal{S}^+ \cup \mathcal{S}^-\cup\mathcal{S}^0$, generator $\bar{\bf T}=[\bar{T}_{ij}]_{i,j\in\mathcal{S}^+ \cup \mathcal{S}^-\cup\mathcal{S}^0}$ and real-valued rates $c_i$, $r_i$, $i\in\mathcal{S}^+ \cup \mathcal{S}^-\cup\mathcal{S}^0$, is such that $\mathcal{S}^0=\mathcal{S}_0\not=\varnothing$.

Consider a related SFFM $\{(\varphi(t), X(t), Y(t)):t\geq 0\}$ driven by CTMC $\{(\varphi(t)):t\geq 0\}$ with state space $\mathcal{S}^+ \cup \mathcal{S}^-$, generator ${\bf T}=[T_{ij}]_{i,j\in\mathcal{S}^+ \cup \mathcal{S}^-}$ and nozero rates $c_i$, $r_i$, $i\in\mathcal{S}^+ \cup \mathcal{S}^-$, such that
	\begin{eqnarray*}
		{\bf T}&=&
		\left[
		\begin{array}{cc}
			{\bf T}^{++}&
			{\bf T}^{+-}\\
			{\bf T}^{-+}&
			{\bf T}^{--}
		\end{array}	
		\right]	
		=
		\left[
		\begin{array}{cc}
			\bar{\bf T}^{++}+\bar{\bf T}^{+0}(-\bar{\bf T}^{00})^{-1}\bar{\bf T}^{0+}&
			\bar{\bf T}^{+-}+\bar{\bf T}^{+0}(-\bar{\bf T}^{00})^{-1}\bar{\bf T}^{0-}\\
			\bar{\bf T}^{-+}+\bar{\bf T}^{-0}(-\bar{\bf T}^{00})^{-1}\bar{\bf T}^{0+}&
			\bar{\bf T}^{--}+\bar{\bf T}^{-0}(-\bar{\bf T}^{00})^{-1}\bar{\bf T}^{0-}
		\end{array}	
		\right],
	\end{eqnarray*}
where for $\ell,m\in\{+,-,0\}$, the block matrix $\bar{\bf T}^{\ell m}=[\bar{T}_{ij}^{\ell m}]_{i\in\mathcal{S}^{\ell},j\in\mathcal{S}^{m}}$ of the generator $\bar{\bf T}$ records the transition rates from phases in $\mathcal{S}^{\ell}$ to phases in $\mathcal{S}^m$, so that for $\ell,m\in\{+,-\}$, the block matrix
$${\bf T}^{\ell m}
=
\bar{\bf T}^{\ell m}+\bar{\bf T}^{\ell 0}(-\bar{\bf T}^{00})^{-1}\bar{\bf T}^{0m}=\bar{\bf T}^{\ell m}+\bar{\bf T}^{\ell 0}\int_{t=0}^{\infty}e^{\bar{\bf T}^{00}t}dt\bar{\bf T}^{0m}$$
of the generator ${\bf T}$ records the transition rates from phases in $\mathcal{S}^{\ell}$ to phases in $\mathcal{S}^m$ via a possible visit to the set $\mathcal{S}^0$, as described by Asmussen in~\cite{asmu95}.

Next, if $|c_i|=|r_i|=1$ for all $i\in\mathcal{S}^+ \cup \mathcal{S}^-$, then by~\cite{asmu95}, the SFFM $\{(\bar\varphi(t),\bar X(t),\bar Y(t)):t\geq 0\}$ is statistically equivalent to the SFFM $\{(\varphi(t), X(t), Y(t)):t\geq 0\}$, that is, the distibution at time $t$ is the same in both processes for all $t\geq 0$.

Furthermore, by Bean, O'Reilly and Taylor~\cite{BOT}, the key matrix $\bar{\bf Q}=|{\bf R}|^{-1}{\bf T}$, ${\bf R}=diag(r_j)_{j\in\mathcal{S}_+\cup\mathcal{S}_-}$, referred to as the fluid generator of the SFM $\{(\bar\varphi(t),\bar Y(t)):t\geq 0\}$, is equal to the fluid generator ${\bf Q}=|{\bf R}|^{-1}{\bf T}$ of the SFM $\{(\varphi(t),Y(t)):t\geq 0\}$, and so by~\cite{BOT} the distributions of quantities at level-dependent stopping times, which directly follow from the expressions involving the fluid generator $|{\bf R}|^{-1}{\bf T}$, are the same in both processes, given the same initial distribution.

Therefore, since there is no change in $X$ or $Y$ when $i\in\mathcal{S}^0=\mathcal{S}_0$, it follows that all quantities derived here for the SFFM $\{(\varphi(t), X(t), Y(t)):t\geq 0\}$ with $\mathcal{S}^0=\mathcal{S}_0=\varnothing$, can be used for the SFFM $\{(\bar\varphi(t), \bar X(t), \bar Y(t)):t\geq 0\}$ with $\mathcal{S}^0=\mathcal{S}_0\not=\varnothing$, since they all record some distributions at the level-dependent stopping times.

\end{Remark}

\subsection{Level-dependent stopping times in the SFMs}\label{sec:stop_SFM}

First, we summarise some relevant results for the SFMs from~\cite{SOB17}. The aim of this section is to present physical interpretations of some key quantities from the theory of the SFMs and the intuitions which we will later apply in the proofs in Sections~\ref{sec:mainresult}-\ref{sec:returnzero}.

Consider the SFM $\{(\varphi(t),Y(t)):t\geq 0\}$ with rates $r_i$, defined in the Introduction, and let $\{(\varphi(t),\widehat Y(t)):t\geq 0\}$ be the corresponding SFM with rates $|r_i|$ so that the fluid level $\widehat Y(\cdot)$ may only increase due to $d \widehat Y(t)/dt = |r_{\varphi(t)}| > 0$. We refer to $\widehat Y(\cdot)$ as the in-out fluid of the process $Y(\cdot)$. We interpret $\widehat Y(t)$ as the total amount of fluid that has flowed into or out of the buffer $\widetilde Y(\cdot)$ at time $t$, where $\widetilde Y(.)$ is the process obtained by removing a lower boundary in $Y(\cdot)$ defined in the Introduction.

Following~\cite{BOT}, denote
\begin{equation*}\label{eq:omegay}
\omega(y)=\inf\left\{t>0: \int_{u=0}^t |r_{\varphi(u)}| du =y\right\},
\end{equation*}
where $\widehat Y(t)=\int_{u=0}^t |r_{\varphi(u)}|du$ is interpreted as the total in-out fluid $\widehat Y(\cdot)$ of the process $Y(\cdot)$ at time~$t$, and $\omega(y)$ is the first time at which the in-out fluid $\widehat Y(\cdot)$ hits level $z+y$ given start in $\widehat Y(0)=z$, see Figure~\ref{fig:omega(y)}.

\begin{figure}
	\begin{center}
		\begin{tikzpicture}[>=stealth,redarr/.style={->}]

		\draw [ ->] (0,0) -- (0,7);
		\draw [ ->] (0,0) -- (10,0);
		\draw (-0.6,7) node[anchor=north, below=-0.17cm] {\footnotesize{\color{black} $\widehat Y(t)$}};
		\draw (10,0.5) node[anchor=north, below=-0.17cm] {\footnotesize{\color{black} $t$}};
		\draw (-0.3,1) node[anchor=north, below=-0.17cm] {\footnotesize{\color{black} $z$}};
		\draw (-0.8,6) node[anchor=north, below=-0.17cm] {\footnotesize{\color{black} $z+y$}};
		\draw (7,-0.2) node[anchor=north, below=-0.17cm] {\footnotesize{\color{black} $\omega(y)$}};
		
		\draw [dashed] (0,1) -- (10,1);
		\draw [dashed] (0,6) -- (10,6);
		\draw [dashed] (7,6) -- (7,0);
		
		\draw [ ->] (0,1) -- (2,1.5)-- (2.5,4) -- (5,4.5) -- (7,6);
		
		\node at (2,1.5) [red,circle,fill,inner sep=1.3pt]{};
		\node at (2.5,4) [red,circle,fill,inner sep=1.3pt]{};
		\node at (5,4.5) [red,circle,fill,inner sep=1.3pt]{};

		\end{tikzpicture}
	\end{center}
	\caption{\footnotesize In-out fluid $\widehat Y(t)$ of the SFM $\{(\varphi(t),Y(t)):t\geq 0\}$ with rates $r_i$ is an SFM with rates $|r_i|$. }
	\label{fig:omega(y)}
\end{figure}
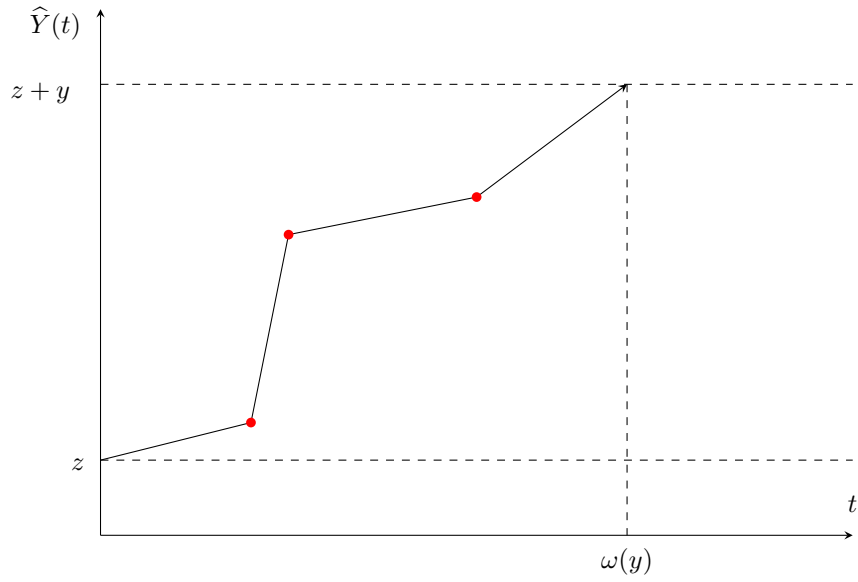

\begin{figure}
	\begin{center}
		\begin{tikzpicture}[>=stealth,redarr/.style={->}]

		\draw [ ->] (0,0) -- (0,7);
		\draw [ ->] (0,0) -- (10,0);
		\draw (-0.8,7) node[anchor=north, below=-0.17cm] {\footnotesize{\color{black} $ Y(t)$}};
		\draw (10,0.5) node[anchor=north, below=-0.17cm] {\footnotesize{\color{black} $t$}};
		\draw (7,-0.2) node[anchor=north, below=-0.17cm] {\footnotesize{\color{black} $\theta_z$}};
		\draw (0,-0.2) node[anchor=north, below=-0.17cm] {\footnotesize{ $0$}};
		
		\draw [dashed] (7,6) -- (7,0);

		
		\draw [blue, ->] (0,1) -- (1,2) -- (2,1.5)-- (2.5,5) -- (3.5,4) -- (6,3) -- (7,1);
		\node at (1,2) [red,circle,fill,inner sep=1.3pt]{};
		\node at (2,1.5) [red,circle,fill,inner sep=1.3pt]{};
		\node at (2.5,5) [red,circle,fill,inner sep=1.3pt]{};
		\node at (3.5,4) [red,circle,fill,inner sep=1.3pt]{};
		\node at (6,3) [red,circle,fill,inner sep=1.3pt]{};
		\draw (-0.3,1) node[anchor=north, below=-0.17cm] {\footnotesize{\color{black} $z$}};
		\draw [dashed] (0,1) -- (10,1);

		\end{tikzpicture}
	\end{center}
	\caption{\footnotesize The SFM $\{(\varphi(t),Y(t)):t\geq 0\}$ starts from level $z$ in some phase $i\in\mathcal{S}^+$ and first returns to level $z$ at time $\theta_z$ and does so in some phase $j\in\mathcal{S}^-$. }
	\label{fig:theta}
\end{figure}
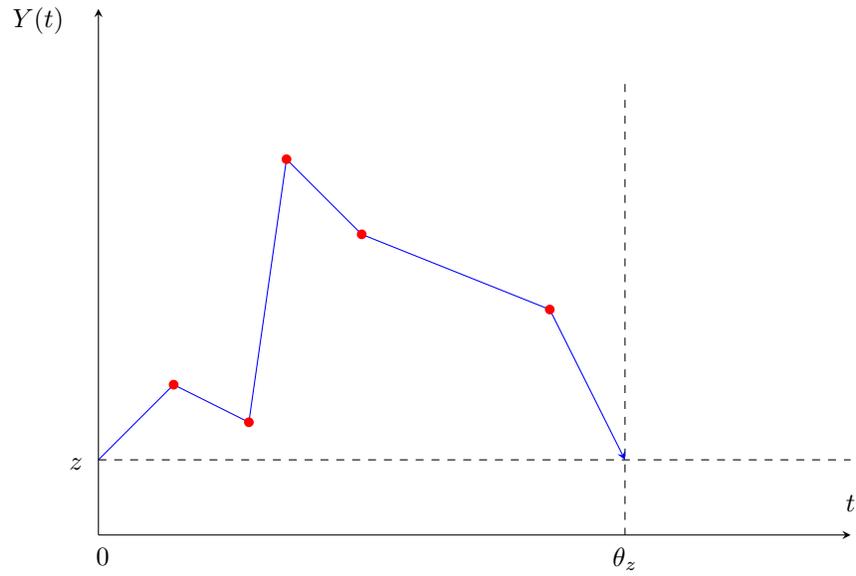

Next, to describe the distribution of the busy period~\cite{1994Abusy,asmu95,2001BST,BOT,rama99}, denote
\begin{equation}\label{eq:thetaz}
\theta_z=\inf\left\{t>0: Y(t)=z\right\},
\end{equation}
interpreted as the first time at which the fluid level hits level $z$. Define matrix ${\bf\Psi}=[\Psi_{ij}]_{i\in\mathcal{S}^+,j\in\mathcal{S}^{-}}$ recording the probabilities $\Psi_{ij}=\mathbb{P}(\varphi(\theta_z)=j \ | \ \varphi(0)=i,Y(0)=z)$ of first return to level $z$ and doing so in phase $j\in\mathcal{S}^{-}$ assuming start from level $z$ in phase $i\in\mathcal{S}^{+}$ in the process $\{(\varphi(t),Y(t)):t\geq 0\}$. We have dropped the index $z$ from the notation since the paths contributing to $\Psi_{ij}$ do not depend on~$z$, see Figure~\ref{fig:theta}.

Also let ${\bf \Xi}=[\Xi_{ij}]_{i\in\mathcal{S}^-,j\in\mathcal{S}^{+}}$, $\Xi_{ij}=\mathbb{P}(\varphi(\theta_z)=j \ | \ \varphi(0)=i,\widetilde Y(0)=z)$, be a matrix with a symmetrical meaning for the process $\widetilde Y(\cdot)$ without a lower boundary, defined in the Introduction. That is, $\Xi_{ij}$ is the probability of first return to level $z$ and doing so in phase $j\in\mathcal{S}^{+}$ assuming start from level $z$ in phase $i\in\mathcal{S}^{-}$ in the process $\{(\varphi(t),\widetilde Y(t)):t\geq 0\}$.

Further, as introduced in~\cite{SOB17}, denote
\begin{equation*}
h_+(t)=\int_{u=0}^{t} \max \{0,r_{\varphi(u)}\} du,
\end{equation*}
interpreted as the total upward shift in $Y(\cdot)$ at time $t$. The term ``upward shift'' reflects the fact that the quantity $h_+(\cdot)$ accumulates only at times $u$ such that $r_{\varphi(u)}>0$ and remains constant otherwise.

Then, $h_+(\omega(y))$ is the total upward shift in $Y(\cdot)$ at the first time at which the in-out fluid $\widehat Y(\cdot)$ of the process $Y(\cdot)$ hits level $z+y$ given start in $\widehat Y(0)=z$.

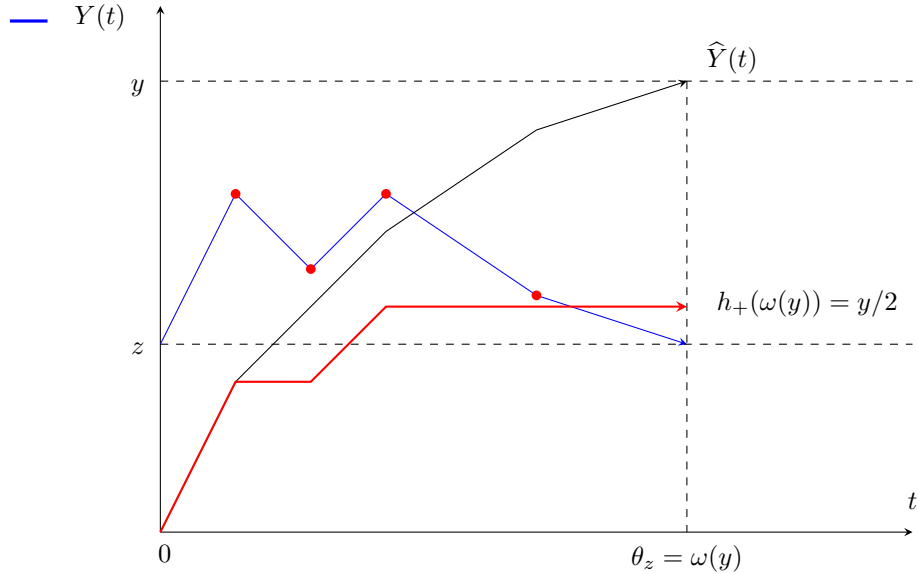
\begin{figure}
	\begin{center}
		\begin{tikzpicture}[>=stealth,redarr/.style={->}]

		\draw [ ->] (0,0) -- (0,7);
		\draw [ ->] (0,0) -- (10,0);
		
		\draw (-0.8,7) node[anchor=north, below=-0.17cm] {\footnotesize{\color{black} $ Y(t)$}};
		\draw [blue,very thick] (-2,6.8) -- (-1.5,6.8);
		
		\draw (7.6,6.5) node[anchor=north, below=-0.17cm] {{\color{black} \footnotesize$ \widehat Y(t)$}};

		\draw (10,0.5) node[anchor=north, below=-0.17cm] {\footnotesize{\color{black} $t$}};
		\draw (7,-0.2) node[anchor=north, below=-0.17cm] {\footnotesize{$\theta_z=\omega(y)$}};
		\draw (0,-0.2) node[anchor=north, below=-0.17cm] {\footnotesize{ $0$}};
		
		\draw [dashed] (7,6) -- (7,0);
	
		\draw [blue, ->] (0,2.5) -- (1,4.5) -- (2,3.5)-- (3,4.5) -- (5,3.15) -- (7,2.5);
		\node at (1,4.5) [red,circle,fill,inner sep=1.3pt]{};
		\node at (2,3.5) [red,circle,fill,inner sep=1.3pt]{};
		\node at (3,4.5) [red,circle,fill,inner sep=1.3pt]{};
		\node at (5,3.15) [red,circle,fill,inner sep=1.3pt]{};

\draw [black, ->] (0,0) -- (1,2) -- (2,3)-- (3,4) -- (5,5.35) -- (7,6);

		\draw (-0.3,2.5) node[anchor=north, below=-0.17cm]
		{\footnotesize{\color{black} $z$}};
		\draw [dashed] (0,2.5) -- (10,2.5);
		
		\draw (-0.3,6) node[anchor=north, below=-0.17cm] {\footnotesize{\color{black} $y$}};
		\draw [dashed] (0,6) -- (10,6);
		
		\draw [red, thick, ->] (0,0) -- (1,2) -- (2,2) -- (3,3) --(7,3);
		
		\draw (8.6,3.2) node[anchor=north, below=-0.17cm] {{\color{black} \footnotesize$ h_+(\omega(y))=y/2$}};

		\end{tikzpicture}
	\end{center}
	\caption{\footnotesize The SFM $\{(\varphi(t),Y(t)):t\geq 0\}$ starts from level $z$ in some phase $i\in\mathcal{S}^+$ and first returns to level $z$  and does so in some phase $j\in\mathcal{S}^-$.  }
	\label{fig:hplus}
\end{figure}

The quantity $h_+(\omega(y))$ is of interest in the analysis due to the useful property that in the paths contributing to $\Psi_{ij}$, at time $\theta_z$ the total upward shift is equal to the half of the total in-out fluid, that is,
\begin{equation*}
\theta_z =  \omega(y) \iff h_+(\omega(y))=y/2 ,
\end{equation*}
see Figure~\ref{fig:hplus}. This observation was used in~\cite{SOB17} to derive the following results.

For $0\leq x\leq y$, let $\widetilde {\bf f}_y(x)=[\widetilde f_y(x)_{ij}]_{i,j\in  \mathcal{S}_+ \cup \mathcal{S}_-}$ be the matrix such that, for all $i,j \in  \mathcal{S}^+ \cup \mathcal{S}^-$,
\begin{eqnarray*}
	\widetilde f_y(x)_{ij} &=& \frac{d}{dx}\mathbb{P}
	(h_+(\omega(y))\leq x,   \varphi(\omega(y))=j\
	|\ Y(0)=0, \varphi(0)=i)
\end{eqnarray*}
is the probability density that the total upward shift in $Y(\cdot)$ at time $\omega(y)$ is $h_+(\omega(y))=x$ and the phase is $\varphi(\omega(y))=j$, given that the process starts in phase $i$ at time zero.


Recall that ${\bf R}=diag(r_j)_{j\in\mathcal{S}^+ \cup \mathcal{S}^-}$ and let $\widecheck{\bf R}=diag(  I(r_j>0)  )_{j\in \mathcal{S}^+ \cup \mathcal{S}^-}$. Then, as shown in~\cite[Corollary 3.1]{SOB17}, for $s>0$, $y>0$,
\begin{eqnarray*}
	[e^{|{\bf R}|^{-1}({\bf T}-s \widecheck{\bf R})y}]_{ij}
	&=&\mathbb{E}(e^{-sh_+(\omega(y))}  I(\varphi(\omega(y))=j)\;|\;\varphi(0)=i)
	=\int_{x=0}^{\infty}e^{-sx}\widetilde f_y(x)_{ij}dx,
\end{eqnarray*}
and so $\widetilde f_y(x)_{ij}$ is the inverse of the LST $[e^{|{\bf R}|^{-1}({\bf T}-s \widecheck{\bf R})y}]_{ij}$ of the distribution of the total upward shift in $Y(\cdot)$ accumulated by the time the total in-out fluid of the process $Y(\cdot)$ first reaches level $y$ and does so in phase $j$, given the process starts at time zero in phase $i$.

This result was achieved in~\cite{SOB17} by interpreting the upward shift as a reward earned at rates $r_i$ whenever the fluid level in $Y(\cdot)$ is increasing, and then noticing that the Laplace-Stieltjes transform of the distribution of this reward accumulates at the rate $|r_i|^{-1}(T_{ii}-s r_i)$, with respect to the in-out fluid level, whenever the process $\varphi(\cdot)$ remains in some $i\in\mathcal{S}^+$. This results intuitively makes sense since by~\cite{BOT2,BOT,BOT5}, $e^{(|{\bf R}|^{-1}{\bf T})y}$ is the conditional probability matrix of the distribution of $\varphi(\cdot)$ at time $\omega(y)$,  with $[e^{|{\bf R}|^{-1}({\bf T})y}]_{ij}
	=\mathbb{P}( \varphi(\omega(y))=j\;|\;\varphi(0)=i)$.

Note that $\widetilde f_y(y/2)_{ij}$ is the probability density that the total upward shift in $Y(\cdot)$ at time $\omega(y)$ is $h_+(\omega(y))=y/2$ and the phase is $\varphi(\omega(y))=j$, given that the process starts in phase $i$ at time zero. We then have $\widetilde Y(\omega(y))=0$, which leads to the following quantity of interest.

Let ${\bf M}=[M_{ij}]$ be a matrix defined in~\cite{SOB17} by
\begin{equation*}\label{Mdef}
M_{ij}=\int_{y=0}^{\infty}\widetilde f_y(y/2)_{ij}dy,
\end{equation*}
and partitioned according to $\left(\mathcal{S}^+\cup\mathcal{S}^-\right) \times \left(\mathcal{S}^+\cup\mathcal{S}^-\right)$ as
\begin{equation*}
{\bf M}=
\left[
\begin{array}{cc}
{\bf M}^{++}&{\bf M}^{+-}\\
{\bf M}^{-+}&{\bf M}^{--}
\end{array}
\right].
\end{equation*}

Then, as shown in~\cite{SOB17}, if the SFM $\{(\varphi(t),\widetilde Y(t)):t\geq 0\}$ with $\widetilde Y(\cdot)\in(-\infty,+\infty)$, generator ${\bf T}$ and rates $r_i$ is transient, then $M_{ij}$ is the expected number of visits to level $0$ and doing so in phase $j$ given start from level $0$ in phase $i$, in such process over the infinite time horizon $[0,\infty)$, and
\begin{eqnarray}
	{\bf M}&=&{\bf\Phi}
	+{\bf M}{\bf\Phi}	
	=
	{\bf\Phi}
	\left({\bf I}-{\bf\Phi}
	\right)^{-1}
	\label{getM},
	\nonumber\\
	{\bf\Phi}&=&
	\left(
	{\bf I} + {\bf M}
	\right)^{-1}{\bf M}
	,
	\label{Meq}
\end{eqnarray}
where
\begin{eqnarray*}
	{\bf\Phi}&=&\left[
	\begin{array}{cc}
		{\bf 0}& {\bf\Psi}\\
		{\bf \Xi}& {\bf 0}
	\end{array}
	\right].
\end{eqnarray*}


\subsection{Distribution of $X(\cdot)$ at the stopping time $\omega(y)$ in $\{(\varphi(t),Y(t)):t\geq 0\}$}

We now adapt the results in~\cite[Lemmas 3 and 4]{BeaOreFF} to our assumption that the rates~$r_i(x)=r_i$ are constant and that $\mathcal{S}^0=\varnothing$, $\mathcal{S}_0=\varnothing$. The results in~\cite[Lemmas 3 and 4]{BeaOreFF} involve operator-analytic methods which are required for the theoretical analysis of the general model, but no algorithmic methods for a practical analysis are given there. Below, we adapt these results to the model considered here and then develop matrix-analytic methods which lead to computational methods for a numerical analysis.

Here and throughout, we write $\mathcal{A}_v$ to denote a set $\mathcal{A}_v=[0,v]$ for some $v>0$. We note that the quantities analysed here may be evaluated for any sets $\mathcal{A}\in\{(u,v),(u,v],[u,v),[u,v]\}$ in terms of sets $\mathcal{A}_v$ and $\mathcal{A}_u$, since for $v>u$, we have $\mu \Gamma ([u,v])=\mu \Gamma ([u,v))=\mu\Gamma ([0,v])-\mu\Gamma ([0,u])$, where $\mu$ is a measure and $\Gamma$ is an operator, for any quantities of type $\mu\Gamma$ considered here.

Let $\mu=[\mu_i]_{i\in\mathcal{S}}$ be some measure with density $\nu(x)=[\nu_i(x)]_{i\in\mathcal{S}}$ for all $x>0$ and atom $p(0)=[p_i(0)]_{i\in\mathcal{S}}$ such that
 \begin{eqnarray*}
\mu_i(\mathcal{A})=\int_{x\in \mathcal{A}}\nu_i(x)dx + p_i(0)I(0\in\mathcal{A})
 \end{eqnarray*}
 for any $i\in\mathcal{S}$ and any set $\mathcal{A}\in\{(u,v),(u,v],[u,v),[u,v]\}$, $v\geq u\geq 0$.

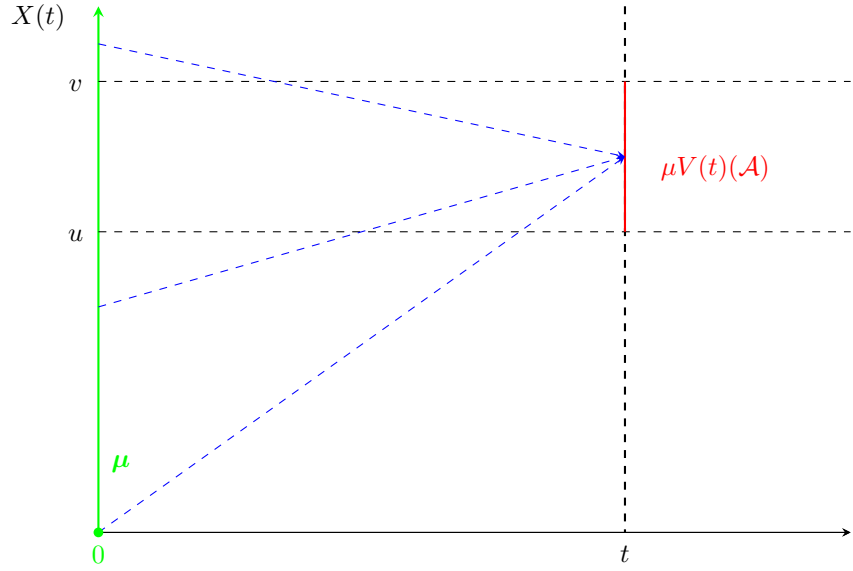
\begin{figure}
	\begin{center}
		\begin{tikzpicture}[>=stealth,redarr/.style={->}]

		\draw [thick, green, ->] (0,0) -- (0,7);
		\draw [ ->] (0,0) -- (10,0);
		\draw (-0.8,7) node[anchor=north, below=-0.17cm] {\footnotesize{\color{black} $X(t)$}};
		\draw (-0.3,4) node[anchor=north, below=-0.17cm] {\footnotesize{\color{black} $u$}};
		\draw (-0.3,6) node[anchor=north, below=-0.17cm] {\footnotesize{\color{black} $v$}};
		\draw (7,-0.2) node[anchor=north, below=-0.17cm] {\footnotesize{\color{black} $t$}};
		\draw (0,-0.2) node[anchor=north, below=-0.17cm] {\footnotesize{\color{green} $0$}};
		
		\draw (0.3,1) node[anchor=north, below=-0.17cm] {\footnotesize{\color{green} $\bmu$}};
		\draw (8.2,5) node[anchor=north, below=-0.17cm] {\footnotesize{\color{red} $\mu V(t) (\mathcal{A})$}};

		\draw [thick,dashed] (7,7) -- (7,0);
		\draw [dashed] (0,4) -- (10,4);
		\draw [dashed] (0,6) -- (10,6);
		
		\draw [red, thick] (7,6) -- (7,4);
		
		\draw [ ->] (0,0) -- (10,0);
		
		\draw [ blue, dashed,->] (0,0) -- (7,5);
		\draw [ blue, dashed,->] (0,3) -- (7,5);
		\draw [ blue, dashed,->] (0,6.5) -- (7,5);

		\node at (0,0) [green,circle,fill,inner sep=1.3pt]{};

		\end{tikzpicture}
	\end{center}
	\caption{\footnotesize Destination at time $t$: $\mathcal{A}=[u,v)$. The SFM $\{(\varphi(t),X(t)):t\geq 0\}$ starts in some phase $i\in\mathcal{S}$ and level $X(0)$ according to the initial distribution $\mu$.  }
	\label{fig:Vt}
\end{figure}

Define the operator  $V(t)=[V_{ij}(t)]_{i,j\in\mathcal{S}}$, $t>0$, such that, for any set $\mathcal{A}\in\{(u,v),(u,v],[u,v),[u,v]\}$, $v\geq u\geq 0$, the operators $V_{ij}(t)$ are given by
\begin{equation*}
\hspace*{-10pt}
\mu_i V_{ij}(t) (\mathcal{A})=\int_{x=0}^{\infty}d\mu_i(x)\mathbb{P}[\varphi(t)=j,X(t)\in\mathcal{A}\ |\ \varphi(0)=i,X(0)=x],
\end{equation*}
where $\mu_i V_{ij}(t) (\mathcal{A})$ is the total probability of the process $\{(\varphi(t),X(t)):t\geq 0\}$ being in the destination set $(j,\mathcal{A})$ at time $t$, assuming that it starts at time zero  in $i\in\mathcal{S}$ according to the measure $\mu$. Further, $\mu_i V_{ij}(t) ([u,v))=\mu_i V_{ij}(t) ([u,v])=\mu_i V_{ij}(t) (\mathcal{A}_v)-\mu_i V_{ij}(t) (\mathcal{A}_u)$. We illustrate this in Figure~\ref{fig:Vt}.

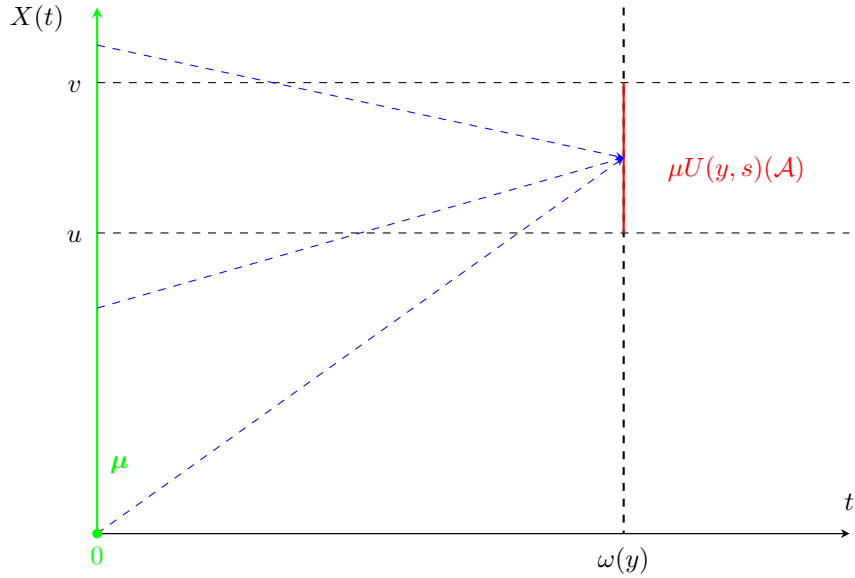
\begin{figure}
	\begin{center}
		\begin{tikzpicture}[>=stealth,redarr/.style={->}]

		\draw [thick, green, ->] (0,0) -- (0,7);
		\draw [ ->] (0,0) -- (10,0);
		\draw (-0.8,7) node[anchor=north, below=-0.17cm] {\footnotesize{\color{black} $X(t)$}};
		\draw (-0.3,4) node[anchor=north, below=-0.17cm] {\footnotesize{\color{black} $u$}};
		\draw (-0.3,6) node[anchor=north, below=-0.17cm] {\footnotesize{\color{black} $v$}};
		\draw (7,-0.2) node[anchor=north, below=-0.17cm] {\footnotesize{\color{black} $\omega(y)$}};
		\draw (0,-0.2) node[anchor=north, below=-0.17cm] {\footnotesize{\color{green} $0$}};
		
		\draw (0.3,1) node[anchor=north, below=-0.17cm] {\footnotesize{\color{green} $\bmu$}};
		
		
		\draw (8.5,5) node[anchor=north, below=-0.17cm] {\footnotesize{\color{red} $\mu U(y,s)(\mathcal{A})$}};
		
		\draw (10,0.5) node[anchor=north, below=-0.17cm] {\footnotesize{\color{black} $t$}};

		\draw [thick,dashed] (7,7) -- (7,0);
		\draw [dashed] (0,4) -- (10,4);
		\draw [dashed] (0,6) -- (10,6);
		
		\draw [red, thick] (7,6) -- (7,4);
		
		\draw [ ->] (0,0) -- (10,0);
		
		\draw [ blue, dashed,->] (0,0) -- (7,5);
		\draw [ blue, dashed,->] (0,3) -- (7,5);
		\draw [ blue, dashed,->] (0,6.5) -- (7,5);

		\node at (0,0) [green,circle,fill,inner sep=1.3pt]{};

		\end{tikzpicture}
	\end{center}
	\caption{\footnotesize Destination at time $\omega(y)$: $\mathcal{A}=[u,v)$. The SFM $\{(\varphi(t),X(t)):t\geq 0\}$ starts in some phase $i\in\mathcal{S}$ and level $X(0)$ according to the initial distribution $\mu$. The SFM $\{(\varphi(t),Y(t)):t\geq 0\}$ starts in level $Y(0)=z$ and evolves until time $\omega(y)$, see Figure~\ref{fig:omega(y)}. }
	\label{fig:Ut}
\end{figure}

\begin{lem}\label{lemma_opB_bounded}
	We have,
	\begin{equation*}\label{B}
	V(t)=e^{Bt},
	\end{equation*}
	with $B=[B_{ij}]_{i,j\in\mathcal{S}}$, where the operators $B_{ij}$ are given by the following expressions for any set $\mathcal{A}_v$, $v> 0$.
	
	For all $i,j \in {\mathcal S}$, $i\neq j$,
	\begin{eqnarray*}
	\mu_i B_{ij} (\mathcal{A}_v) &=& T_{ij}  \mu_i(\mathcal{A}_v),
	\end{eqnarray*}
	and for all $j \in {\mathcal S}$,
	\begin{eqnarray*}
			\mu_j B_{jj}(\mathcal{A}_v)&=&
		T_{jj}  \mu_j(\mathcal{A}_v)
		-
		c_j
		\nu_j(v).
	\end{eqnarray*}
\end{lem}

\begin{Remark}
The following interpretation of the operator $B$ was established in~\cite{BeaOreFF}. The operator $B$ is the infinitesimal generator with respect to time of the process $Y(\cdot)$. The term $T_{ij}  \mu_i(\mathcal{A}_v)$ represents the stochastic jumps from phase $i$ to phase $j \not= i$. This occurs with rate $T_{ij}$ and so, in order to ensure that such a transition ends up in the destination set $\mathcal{A}_v$, this is multiplied by the probability mass $\mu_i(\mathcal{A}_v)$ of starting in $\mathcal{A}_v$. The term $T_{jj}  \mu_j(\mathcal{A}_v) -c_j\nu_j(v)$ represents the stochastic jumps out of phase $j$ and the drift within the level, where $c_j\nu_j(v)$ represents the drift in and out of $\mathcal{A}_v$ when $c_j>0$ and $c_j<0$, respectively.
\end{Remark}
{\bf Proof:} The result follows directly by applying~\cite[Lemma 3]{BeaOreFF} for $\mathcal{A}\in\{[u,v),[u,v]\}$ with $u=0$, $v>u$ and letting $\mathcal{S}^0=\varnothing$, $\mathcal{S}_0=\varnothing$. Note that the second case of~\cite[Lemma 3]{BeaOreFF} does not occur here since $v>0$. The detailed proof of this result is presented in the Appendix for convenience of the reader. \hfill\rule{9pt}{9pt}

Define the operator $U(y,s)=[U_{ij}(y,s)]_{i,j\in\mathcal{S}}$, such that the operators $U_{ij}(y,s)$, for any set $\mathcal{A}\in\{(u,v),(u,v],[u,v),[u,v]\}$, $v\geq u\geq 0$, are given by
\begin{eqnarray*}
	\mu_i U_{ij}(y,s) (\mathcal{A})=\int_{x=0}^{\infty}d\mu_i(x)\mathbb{E}
	[e^{-s\omega(y)}
 I(\varphi(\omega(y))=j,X(\omega(y))\in\mathcal{A})\ |\ \varphi(0)=i,X(0)=x],
\end{eqnarray*}
where $\mu_i U_{ij}(y,s) (\mathcal{A})$ is the LST of the time taken for  the total amount of fluid that has flowed into or out of the buffer $Y$ to reach $y$ (that is, at time $\omega(y)$) and do so with the process $\{(\varphi(t),X(t)),t\geq 0\}$ in the destination set $(j,\mathcal{A})$, assuming the process starts in $i$ at time zero according to measure~$\mu$. Further, $\mu_iU_{ij}(y,s)([u,v))=\mu_iU_{ij}(y,s)([u,v])=\mu_iU_{ij}(y,s)(\mathcal{A}_v)-\mu_iU_{ij}(y,s)(\mathcal{A}_u)$. We illustrate this in Figure~\ref{fig:Ut}.

\begin{lem}\label{lemma_D} For all $y \geq 0$ and $s \in \CC$ with $\Re(s)\geq 0$,
	\begin{equation*}\label{genD}
	U(y,s)=e^{D(s)y},
	\end{equation*}
	with $D(s)=\left[D_{ij}(s)\right]_{i,j\in\mathcal{S}}$ where
	\begin{equation*}
	D_{ij}(s)=\frac{1}{|r_i|}\left[
	\left(B-sI
	\right)\right]_{ij},\\
	\end{equation*}
	for all $i,j\in\mathcal{S}$.
\end{lem}
\begin{Remark}
	The following interpretation of the operator $D(s)$ was established in~\cite{BeaOreFF}. The operator $D(s)$ is the infinitesimal generator of the LST with respect to the level in the in-out fluid $\widehat{Y}(\cdot)$ of the process~$\widetilde Y(\cdot)$. We denote $D=D(0)$.
\end{Remark}
{\bf Proof:} The result follows directly from~\cite[Lemma 4]{BeaOreFF} by letting $r_i(x)=r_i$ and $\mathcal{S}^0=\varnothing$, $\mathcal{S}_0=\varnothing$. The detailed proof of this result is presented in the Appendix. \hfill\rule{9pt}{9pt}

\section{Main result}\label{sec:mainresult}

Before we state the main result it will be convenient to introduce notations for key vectors and matrices required in our analysis.

Recall that $\mathcal{S}_+=\{j\in\mathcal{S}:c_j>0\}$ and $\mathcal{S}_-=\{j\in\mathcal{S}:c_j< 0\}$. We partition the quantities below according to $\mathcal{S}=\mathcal{S}_+\cup\mathcal{S}_-$. Also, ${\bf R}=diag(r_j)_{j\in\mathcal{S}}$, ${\bf R}_{+}=diag(r_j)_{j\in\mathcal{S}_+}$, ${\bf R}_{-}=diag(r_j)_{j\in\mathcal{S}_-}$, and so similarly, we let ${\bf C}=diag(c_j)_{j\in\mathcal{S}}$, ${\bf C}_{+}=diag(c_j)_{j\in\mathcal{S}_+}$, ${\bf C}_{-}=diag(c_j)_{j\in\mathcal{S}_-}$.

Let $\bnu(x)=[\nu_j(x)]_{j\in\mathcal{S}}$, $x>0$, be the density and ${\bf P}=[ P_j]_{j\in\mathcal{S}}$ be the point mass at $X(0)=0$ of the {\em initial probability distribution} $\bmu=[\mu_j]_{j\in\mathcal{S}}$ of $X(0)$, so that for any set $\mathcal{A}_v$, $v>0$,
\begin{equation*}
	\mu_j (\mathcal{A}_v)=\mathbb{P}(\varphi(0)=j,X(0)\in\mathcal{A}_v)=\int_{x=0}^v \nu_j(x)dx +P_j,
\end{equation*}
and denote
\begin{equation}\label{eq:initialdistr}
\bmu(\mathcal{A}_v)=[\mu_j(\mathcal{A}_v)]_{j\in\mathcal{S}}.
\end{equation}

We partition $\bnu(x)$ and ${\bf P}$ as
\begin{eqnarray*}
{\bf P}&=&
 \left[
\begin{array}{cc}
{\bf P}_+&{\bf P}_-
\end{array}
\right]
=
 \left[
 \begin{array}{cc}
 	{\bf 0}&{\bf P}_-
 \end{array}
 \right],\\
\bnu (x)&=&
\left[
\begin{array}{cc}
\bnu_+(x)&\bnu_-(x)
\end{array}
\right],
\end{eqnarray*}
where we assume ${\bf P}_+={\bf 0}$ since no mass may accumulate at $X(t)=0$ in phases $j$ with positive rates $c_j>0$ at times $t>0$. Also, for $k=0,1,2,\ldots$, we denote
\begin{equation*}
\bnu^{(k)}(0)=\lim_{x\to 0^{+}}\bnu^{(k)}(x)
=
\lim_{x\to 0^{+}}
\left[\frac{d^k}{dx^k}\nu_j(x)\right]_{j\in\mathcal{S}}
=\left[
\begin{array}{cc}
\bnu^{(k)}_+(0)&\bnu^{(k)}_-(0)
\end{array}
\right].
\end{equation*}

We define the $n$-th {\em time} derivative with respect to time evaluated at $t=0$, ${}_{B^n}{\bmu}=[{}_{B^n}\mu_j]_{j\in\mathcal{S}}$, by
\begin{equation*}
{}_{B^n}\mu_j(\mathcal{A}_v)
=\frac{d^n}{dt^n} \sum_i\mu_i [e^{Bt}]_{ij}(\mathcal{A}_v)
\Big|_{t=0}
=\sum_i\mu_i [B^n]_{ij}(\mathcal{A}_v)
,
\end{equation*}
where the operator $B$ is given in Lemma~\ref{lemma_opB_bounded}. Let ${}_{B^n}\bnu(x)=[{}_{B^n}\nu_j(x)]_{j\in\mathcal{S}}$, $x>0$, be the density and ${}_{B^n}{\bf P}=[{}_{B^n}P_j]_{i,j\in\mathcal{S}}$ be the atom of ${}_{B^n}{\bmu}=[{}_{B^n}\mu_j]_{j\in\mathcal{S}}$ such that, for any set $\mathcal{A}_v$, $v>0$,
\begin{equation*}
{}_{B^n}\mu_j(\mathcal{A}_v)
=\sum_i\mu_i [B^n]_{ij}(\mathcal{A}_v)
=\int_{x=0}^v {}_{B^n}\nu_j(x)dx +{}_{B^n}P_j,
\end{equation*}
whenever ${}_{B^n}\bnu(x)$ and ${}_{B^n}{\bf P}$ exist.

Denote, for $k,n=0,1,2\ldots$,
\begin{eqnarray*}
{}_{B^n}\nu_j^{(k)}(x)=\frac{d^k{}_{B^n}\nu_j(x)}{dx^k},\qquad&&
{}_{B^n}{\bnu}^{(k)}(x)=[{}_{B^n}\nu_j^{(k)}(x)]_{j\in\mathcal{S}},\nonumber\\
{}_{B^n}{\mu}^{(k)}_j(\mathcal{A}_v)=\int_{x=0}^v {}_{B^n}\nu_j^{(k)}(x)dx,\qquad&&
{}_{B^n}{\bmu}^{(k)}(\mathcal{A}_v)=[{}_{B^n}\mu_j^{(k)}(\mathcal{A}_v)]_{j\in\mathcal{S}}.
\end{eqnarray*}

Similarly, we define the $n$-th {\em level} derivative with respect to level evaluated at $y=0$, ${}_{D^n}{\bmu}=[{}_{D^n}\mu_j]_{j\in\mathcal{S}}$, by
	\begin{equation*}
		{}_{D^n}\mu_j(\mathcal{A}_v)
		=
		\frac{d^n}{dy^n} \sum_i\mu_i [e^{Dy}]_{ij}(\mathcal{A}_v)
		\Big|_{y=0}
		=\sum_i\mu_i [D^n]_{ij}(\mathcal{A}_v),
	\end{equation*}
where the operator $D$ is given in Lemma~\ref{lemma_D}.

Further, for $n=0,1,2,\ldots$, let ${}_{D^n}\bnu(x)=[{}_{D^n}\nu_j]_{j\in\mathcal{S}}$ be the density and ${}_{D^n}{\bf P}=[{}_{D^n}P_j]_{j\in\mathcal{S}}$ the atom of ${}_{D^n}{\bmu}=[{}_{D^n}\mu_j]_{j\in\mathcal{S}}$ in the sense that, for any set $\mathcal{A}_v$, $v>0$,
\begin{equation*}
	{}_{D^n}\mu_j(\mathcal{A}_v)
	=\sum_i\mu_i [D^n]_{ij}(\mathcal{A}_v)
	=\int_{x=0}^v {}_{D^n}\nu_j(x)dx
	+{}_{D^n}P_j,
\end{equation*}
whenever ${}_{D^n}\bnu(x)$ and ${}_{D^n}{\bf P}$ exist, and for $k,n=0,1,2\ldots$ denote
\begin{eqnarray*}
{}_{D^n}\nu_j^{(k)}(x)=\frac{d^k{}_{D^n}\nu_j(x)}{dx^k},\qquad&&
{}_{D^n}{\bnu}^{(k)}(x)=[{}_{D^n}\nu_j^{(k)}(x)]_{j\in\mathcal{S}},\nonumber\\
{}_{D^n}{\mu}^{(k)}_j(\mathcal{A}_v)=\int_{x=0}^v {}_{D^n}\nu_j^{(k)}(x)dx,\qquad&&
{}_{D^n}{\bmu}^{(k)}(\mathcal{A}_v)=[{}_{D^n}\mu_j^{(k)}(\mathcal{A}_v)]_{j\in\mathcal{S}}.
\end{eqnarray*}

Additionally, for notational convenience we introduce the following quantities which appear in the form of weights in the expressions in Lemmas~\ref{lem:mu_Dn}-\ref{cond_mu_plus_0}. Let ${\bf h}(k,n)$ be the sum of all different products in which $(-|{\bf R}|^{-1}){{\bf C}})$ appears exactly $k$ times and $(|{\bf R}|^{-1}{\bf T})$ exactly $(n-k)$ times, that is,
\begin{eqnarray*}\label{hkn}
{\bf h}(k,n)&=&
\left(-|{\bf R}|^{-1}{{\bf C}}\right)^k \left(|{\bf R}|^{-1}{\bf T}\right)^{n-k}
+
\left(-|{\bf R}|^{-1}{{\bf C}}\right)^{k-1} \left(|{\bf R}|^{-1}{\bf T}\right)^{n-k}
\left(-|{\bf R}|^{-1}{{\bf C}}\right)
\nonumber\\
&&
+\ldots
+\left(|{\bf R}|^{-1}{\bf T}\right)^{n-k}
\left(-|{\bf R}|^{-1}{{\bf C}}\right)^k,
\end{eqnarray*}
with clearly
\begin{equation*}
\sum_{k=0}^n {\bf h}(k,n)
=
(|{\bf R}|^{-1}{\bf T} - |{\bf R}|^{-1}{\bf C})^n
.
\end{equation*}

Furthermore, denote
\begin{equation*}
{\bf f}(k,n)={\bf h}(k,n)(-|{\bf R}|^{-1}{\bf C})^{-n}
\end{equation*}
and
\begin{equation*}
{\bf A}(n)=-{\bf P}_-[(|{\bf R}|^{-1}{\bf T})^n]_{- +}(-{\bf R}_+^{-1}{\bf C}_+)^{-n}
-\sum_{k=1}^{n-1} (\bnu_-)^{(k-1)}(0)[{\bf f}(k,n)]_{- +}
,
\end{equation*}
with the convention that the second term set to zero when $n=1$. For $n\geq 2$ we define
\begin{eqnarray*}
{\bf B}(n)&=& \sum_{\ell=1}^{n-1} {\bf M}(\ell,n),
\end{eqnarray*}
with
\begin{eqnarray*}
{\bf M}(\ell ,n)
&=&
(-1)^{\ell}
\sum_{1\leq k_1<n_1=k_2<n_2,\ldots,n_{\ell-1}=k_{\ell}<n_{\ell}\leq n}
{\bf A}(k_1)
\prod_{i=1}^{\ell}
[{\bf f}(k_i,n_i)]_{++}.
\end{eqnarray*}
Our first result identifies ${}_{D^n}\bmu (\mathcal{A}_v)$ which is key fact for further results and
is interesting for itself.
\begin{lem}\label{lem:mu_Dn}
Assume that the following boundary conditions
\begin{equation}\label{incond}
{}_{D^n}\bnu_+(0)={}_{D^n}{\bf P}_-
(|{\bf R}_-|^{-1}{\bf T}_{- +})
(|{\bf R}_+|^{-1}{\bf C}_+)^{-1}
\end{equation}
are met for all $n\geq 0$. Then, for any set $\mathcal{A}_v$, $v>0$, and for all $n\geq 1$, we have,
\begin{equation}\label{Dn_muA}
{}_{D^n}\bmu (\mathcal{A}_v)=
{}_{D^n}{\bmu}^{(0)}(\mathcal{A}_v)
+
{}_{D^n}{\bf P}
=
\int_{x=0}^v {}_{D^n}\bnu (x)dx
+
{}_{D^n}{\bf P}
=
\sum_{k=0}^n  {\bmu}^{(k)}(\mathcal{A}_v) {\bf h}(k,n)
+
{}_{D^n}{\bf P}
	,
\end{equation}
where
\begin{eqnarray}
{}_{D^n}\bnu (x) &=&
\left[
\begin{array}{cc}
{}_{D^n}\bnu_+(x) & {}_{D^n}\bnu_-(x)
\end{array}
\right]
=
\sum_{k=0}^n   {\bnu}^{(k)}(x) {\bf h}(k,n)
\quad\mbox{ for }x>0,
\label{Dn_nux}\\
{}_{D^n}{\bf P}
&=&
\left[
\begin{array}{cc}
{}_{D^n}{\bf P}_+ & {}_{D^n}{\bf P}_-
\end{array}
\right]
=
{\bf P}(|{\bf R}|^{-1}{\bf T})^n+\sum_{k=1}^n \bnu^{(k-1)}(0){\bf h}(k,n),
\label{Dn_P}
\\
{}_{D^n}{\bf P}_+&=&
{\bf 0}
,
\label{Dn_Pplus}
\end{eqnarray}
and
\begin{eqnarray}
{}_{D^n}\bnu_+ (0)&=&
\sum_{k=0}^n  ({\bnu}_+)^{(k)}(0)
[ {\bf h}(k,n)]_{++}
+
\sum_{k=0}^{n-1}   ({\bnu}_-)^{(k)}(0)
[{\bf h}(k,n)]_{- +}
,
\label{Dn_nuplus}
\\
{}_{D^n}\bnu_- (0)&=&
\sum_{k=0}^{n-1}   ({\bnu}_+)^{(k)}(0)
[ {\bf h}(k,n)]_{+ -}
+
\sum_{k=0}^n   ({\bnu}_-)^{(k)}(0)
[{\bf h}(k,n)]_{--}
,
\label{Dn_nubullet0}
\\
{}_{D^n}{\bf P}_-&=&
{\bf P}_-[(|{\bf R}|^{-1}{\bf T})^n]_{--}
+\sum_{k=1}^n (\bnu_{-})^{(k-1)}(0)[{\bf h}(k,n)]_{--}
+\sum_{k=1}^{n-1} (\bnu_{+})^{(k-1)}(0)[{\bf h}(k,n)]_{+-}
.
\label{Dn_Pbullet}
\end{eqnarray}
\end{lem}
\begin{Remark}
Condition~\eqref{incond}, equivalent to ${}_{D^n}{\bf P}_{-}
(|{\bf R}_{-}|^{-1}{\bf T}_{- +})={}_{D^n}\bnu_+(0)(|{\bf R}_+|^{-1}{\bf C}_+)$, ensures the cancellation of relevant terms in the recursive arguments of the proof below so that, given the existence of the density ${}_{D^n}\bnu(x)$ for all $x>0$, and that ${}_{D^n} {\bf P}_+=0$, we have
\begin{eqnarray*}
{}_{D^{n+1}} {\bf P}_+={}_{D^n}{\bf P}_{-}
(|{\bf R}_{-}|^{-1}{\bf T}_{- +})
-{}_{D^n}\bnu_+(0)(|{\bf R}_+|^{-1}{\bf C}_+)={\bf 0},
\end{eqnarray*}
which means that point mass may not accumulate in phases $j$ with positive rates $c_j>0$. Consequently, the density ${}_{D^{n+1}}\bnu(x)$ exists for all $x>0$ as well, and the recursion and Lemmas~\ref{lemma_opB_bounded}-\ref{lemma_D} can be applied again to ${}_{D^{n+1}}\bmu$. The physical interpretation of this condition is that the rate at which the point mass leaves level $0$, ${}_{D^n}{\bf P}_{-}
(|{\bf R}_{-}|^{-1}{\bf T}_{- +})$, is equal to the rate at which the density drifts away from level $0$, ${}_{D^n}\bnu_+(0)(|{\bf R}_+|^{-1}{\bf C}_+)$.
\end{Remark}
{\bf Proof:}\quad We apply mathematical induction to show that~\eqref{Dn_muA}--\eqref{Dn_Pbullet} hold for all $n\geq 1$.

{\bf Step(i).}\quad Consider $n=1$ and so ${}_{D}\bmu (\mathcal{A}_v) $ first. By Lemma~\ref{lemma_D}, for all $i,j\in\mathcal{S}$, and any set $\mathcal{A}_v$, $v>0$,
\begin{equation*}
\mu_iD_{ij}(\mathcal{A}_v)=\frac{1}{|r_i|}\mu_iB_{ij}(\mathcal{A}_v).
\end{equation*}
By Lemma~\ref{lemma_opB_bounded} we have, for $i\not= j$,
\begin{eqnarray*}
\frac{1}{|r_i|}\mu_iB_{ij}(\mathcal{A}_v)
&=&\frac{1}{|r_i|} T_{ij}\int_{x=0}^v \nu_i(x)dx + \frac{1}{|r_i|}T_{ij}P_i
\\
&=&
(\mu_i^{(0)}(\mathcal{A}_v)+P_i)[
\left(|{\bf R}|^{-1}{\bf T}\right)
]_{ij},
\end{eqnarray*}
and for $i=j$,
\begin{eqnarray*}
\frac{1}{|r_j|}\mu_jB_{jj} (\mathcal{A}_v) &=& \frac{1}{|r_j|}T_{jj}\int_{x=0}^v \nu_j(x)dx + \frac{1}{|r_j|}T_{jj}P_j\\
&&+\frac{1}{|r_j|}c_j(\nu_j(0)-\nu_j(v)) - \frac{1}{|r_j|}c_j\nu_j(0)\\
&=&
(\mu_j^{(0)}(\mathcal{A}_v)+P_j)[
\left(|{\bf R}|^{-1}{\bf T}\right)
]_{jj}
-[{\bmu}^{(1)} (\mathcal{A}_v) \left(|{\bf R}|^{-1}{{\bf C}}\right)]_j
-[{\bnu} (0) \left(|{\bf R}|^{-1}{{\bf C}}\right)]_j.
\end{eqnarray*}
Therefore, since ${}_{D^n}\mu_j(\mathcal{A}_v)
=\sum_i\mu_i D^n_{ij}(\mathcal{A}_v)$, we have
\begin{eqnarray*}\label{Dmu_A}
{}_{D}\bmu (\mathcal{A}_v) =
{\bmu}^{(0)}(\mathcal{A}_v)(|{\bf R}|^{-1}{\bf T}) + {\bf P} (|{\bf R}|^{-1}{\bf T})
-{\bmu}^{(1)} (\mathcal{A}_v) (|{\bf R}|^{-1}{\bf C})
-{\bnu} (0) (|{\bf R}|^{-1}{\bf C}),
\end{eqnarray*}
with
\begin{eqnarray*}
{}_D {\bf P}&=&
\left[
\begin{array}{cc}
{}_{D}{\bf P}_+ & {}_{D}{\bf P}_{-}
\end{array}
\right]
=
{\bf P} (|{\bf R}|^{-1}{\bf T})-{\bnu} (0) (|{\bf R}|^{-1}{\bf C})
=
{\bf P}(|{\bf R}|^{-1}{\bf T})^1+\sum_{k=1}^1 \bnu^{(k-1)}(0){\bf h}(k,1).
\end{eqnarray*}
Now, due to the boundary condition~\eqref{incond} for $n=0$, we have the following cancellation of terms,
\begin{equation*}
 {}_D {\bf P}_+={\bf P}_{-}
(|{\bf R}_{-}|^{-1}{\bf T}_{- +})
-\bnu_+(0)(|{\bf R}_+|^{-1}{\bf C}_+)={\bf 0},
\end{equation*}
and so the density of ${}_{D}\bmu$, for $x>0$, exists and is given by
\begin{equation*}
\label{D_dens}
{}_D{\bnu} (x)
= {\bnu} (x) (|{\bf R}|^{-1}{\bf T}) -  {\bnu}^{(1)} (x)
(|{\bf R}|^{-1}{\bf C})
=
\sum_{k=0}^1   {\bnu}^{(k)}(x) {\bf h}(k,n).
\end{equation*}
Moreover, the atom can be identified as
\begin{equation*}\label{D_mass}
{}_D {\bf P}_{-}
=
{\bf P}_{-} (|{\bf R}_{-}|^{-1}{\bf T}_{--})
- {\bnu}_{-} (0) (|{\bf R}_{-}|^{-1}{\bf C}_{-})
=
{\bf P}_{-}[(|{\bf R}|^{-1}{\bf T})^1]_{--}
+\sum_{k=1}^1 (\bnu_{-})^{(k-1)}(0)[{\bf h}(k,n)]_{--}
 ,
\end{equation*}
with the third term in~\eqref{Dn_Pbullet} for $n=1$, given by $\sum_{k=1}^1 (\bnu_{+})^{(k-1)}(0)[{\bf h}(k,1)]_{+-}={\bf 0}$ since $[{\bf h}(1,1)]_{+-}=[(|{\bf R}|^{-1}{\bf C})]_{+-}={\bf 0}$. That is, for any set $\mathcal{A}_v$, $v>0$,
\begin{equation*}
{}_{D}{\bmu}(\mathcal{A}_v)=\int_{x=0}^v {}_D{\bnu} (x)dx+{}_D {\bf P},
\end{equation*}
which completes the proof of \eqref{Dn_muA}--\eqref{Dn_Pbullet} when $n=1$.

{\bf Step(ii).}\quad Next, suppose~\eqref{Dn_muA}--\eqref{Dn_Pbullet} hold for some $n\geq 1$. We will show that it then follows that~\eqref{Dn_muA}--\eqref{Dn_Pbullet} hold for $(n+1)$.
We have by~\eqref{Dmu_A},
\begin{eqnarray}\label{Dnplus1_mu_A}
{}_{D^{n+1}}\bmu (\mathcal{A}_v) =
{}_{D^n}{\bmu}(\mathcal{A}_v)(|{\bf R}|^{-1}{\bf T}) + {}_{D^n}{\bf P} (|{\bf R}|^{-1}{\bf T})
-{}_{D^n}{\bmu}^{(1)} (\mathcal{A}_v) (|{\bf R}|^{-1}{\bf C})
-{}_{D^n}{\bnu} (0) (|{\bf R}|^{-1}{\bf C}),
\end{eqnarray}
with, by~\eqref{incond},
\begin{equation*}\label{Dnplus1_Pplus}
{}_{D^{n+1}} {\bf P}_+={}_{D^n}{\bf P}_{-}
(|{\bf R}_{-}|^{-1}{\bf T}_{- +})
-{}_{D^n}\bnu_+(0)(|{\bf R}_+|^{-1}{\bf C}_+)={\bf 0}.
\end{equation*}
Thus~\eqref{Dn_Pplus} holds for $(n+1)$, and the density of ${}_{D^{n+1}}\bmu$, for $x>0$, exists and is given by
\begin{eqnarray}\label{D2_dens}
{}_{D^{n+1}}{\bnu} (x) &=& {}_{D^n}{\bnu} (x) (|{\bf R}|^{-1}{\bf T})
-  {}_{D^n}{\bnu}^{(1)} (x)
(|{\bf R}|^{-1}{\bf C}).
\end{eqnarray}
Moreover, we have
\begin{eqnarray}\label{DnplusB}
{}_{D^{n+1}} {\bf P}_{-} &=& {}_{D^n}{\bf P}_{-} (|{\bf R}_{-}|^{-1}{\bf T}_{--})
- {}_{D^n}{\bnu}_{-} (0) (|{\bf R}_{-}|^{-1}{\bf C}_{-}).
\end{eqnarray}
Thus, for any set $\mathcal{A}_v$, $v>0$,
\begin{equation*}
{}_{D^{n+1}}{\bmu}(\mathcal{A}_v)=\int_{x=0}^v {}_{D^{n+1}}{\bnu} (x)dx+{}_{D^{n+1}}{\bf P}.
\end{equation*}
Next, assuming that~\eqref{Dn_nux} holds for some $n\geq 1$, we have
\begin{eqnarray*}
{}_{D^{n+1}}{\bnu} (x) &=& {}_{D^n}{\bnu} (x) (|{\bf R}|^{-1}{\bf T})
-  {}_{D^n}{\bnu}^{(1)} (x)
(|{\bf R}|^{-1}{\bf C})
\nonumber\\
&=&
\sum_{k=0}^n   {\bnu}^{(k)}(x) {\bf h}(k,n)(|{\bf R}|^{-1}{\bf T})
+
\sum_{k=0}^n  {\bnu}^{(k+1)}(x) {\bf h}(k,n)
(-|{\bf R}|^{-1}{\bf C})
\nonumber\\
&=&
\sum_{k=0}^n   {\bnu}^{(k)}(x) {\bf h}(k,n)(|{\bf R}|^{-1}{\bf T})
+
\sum_{k=1}^{n+1}  {\bnu}^{(k)}(x) {\bf h}(k-1,n)
(-|{\bf R}|^{-1}{\bf C})
\nonumber\\
&=&
 {\bnu}^{(0)}(x) {\bf h}(0,n+1)
+
\sum_{k=1}^n   {\bnu}^{(k)}(x)
[{\bf h}(k,n)(|{\bf R}|^{-1}{\bf T})+{\bf h}(k-1,n)
(-|{\bf R}|^{-1}{\bf C})]
\nonumber\\
&&
+
   {\bnu}^{(n+1)}(x) {\bf h}(n+1,n+1)
\nonumber\\
&=&
\sum_{k=0}^{n+1}   {\bnu}^{(k)}(x) {\bf h}(k,n+1),
\end{eqnarray*}
and hence~\eqref{Dn_muA}--\eqref{Dn_nux} holds for all $(n+1)$. Also,~\eqref{Dn_nuplus}--\eqref{Dn_nubullet0} for $(n+1)$ follow immediately by noting $[ {\bf h}(n,n)]_{- +}={\bf 0}$ and $[ {\bf h}(n,n)]_{+-}={\bf 0}$.

Furthermore, assuming that~\eqref{Dn_P} holds for some $n\geq 1$, then by~\eqref{Dn_nux} and~\eqref{Dnplus1_mu_A},
\begin{eqnarray*}
{}_{D^{n+1}}{\bf P} &=&
{}_{D^n}{\bf P} (|{\bf R}|^{-1}{\bf T})
-{}_{D^n}{\bnu} (0) (|{\bf R}|^{-1}{\bf C})
\nonumber\\
&=&
{\bf P}(|{\bf R}|^{-1}{\bf T})^n(|{\bf R}|^{-1}{\bf T})+
\sum_{k=1}^n \bnu^{(k-1)}(0){\bf h}(k,n)(|{\bf R}|^{-1}{\bf T})
+\sum_{k=0}^n  {\bnu}^{(k)}(0) {\bf h}(k,n)(-|{\bf R}|^{-1}{\bf C})
\nonumber\\
&=&
{\bf P}(|{\bf R}|^{-1}{\bf T})^n(|{\bf R}|^{-1}{\bf T})+
\sum_{k=1}^n \bnu^{(k-1)}(0){\bf h}(k,n)(|{\bf R}|^{-1}{\bf T})
+\sum_{k=1}^{n+1}  {\bnu}^{(k-1)}(0) {\bf h}(k-1,n)(-|{\bf R}|^{-1}{\bf C})
\nonumber\\
&=&
{\bf P}(|{\bf R}|^{-1}{\bf T})^{n+1}
+\sum_{k=1}^n
\bnu^{(k-1)}(0)[{\bf h}(k,n)(|{\bf R}|^{-1}{\bf T})
+{\bf h}(k-1,n)(-|{\bf R}|^{-1}{\bf C})]
\nonumber\\
&&
+\bnu^{(n)}(0){\bf h}(n,n)(-|{\bf R}|^{-1}{\bf C})
\nonumber\\
&=&
{\bf P}(|{\bf R}|^{-1}{\bf T})^{n+1}
+
\sum_{k=1}^{n+1} \bnu^{(k-1)}(0){\bf h}(k,n+1),
\label{Dnplus1P}
\end{eqnarray*}
and so~\eqref{Dn_P} holds for $(n+1)$. Further,~\eqref{Dn_Pbullet} follows immediately by noting $[ {\bf h}(n,n)]_{- +}={\bf 0}$ and $[ {\bf h}(n,n)]_{+-}={\bf 0}$.
This completes the proof of this lemma. \hfill\rule{9pt}{9pt}

Note we state a condition that needs to be met for all $\bnu_+^{(n)}(0)$, $n\geq 1$, in order for~\eqref{incond} to be satisfied.
\begin{lem}\label{cond_mu_plus_0}
The boundary condition~\eqref{incond} is equivalent to the following recursive expression for $(\bnu_+)^{(n)}(0)$, $n\geq 1$, in terms of $(\bnu_+)^{(k)}(0)$, $0\leq k\leq n-1$,
		\begin{eqnarray}
		(\bnu_+)^{(n)}(0)
		&=&
			{\bf A}(n+1)
			-
			\sum_{k=1}^n
			(\bnu_+)^{(k-1)}(0)
			[{\bf f}(k,n+1)]_{++}
		.
		\label{Dn_nuplus0_rec}
		\end{eqnarray}
		
		Further,~\eqref{Dn_nuplus0_rec} is equivalent to
	\begin{eqnarray}
	(\bnu_+)^{(n)}(0)&=&
	{\bf A}(n+1)+{\bf B}(n+1).
	\label{Dn_nuplus0}
	\end{eqnarray}

\end{lem}
{\bf Proof:}\quad First, we prove~\eqref{Dn_nuplus0_rec}. By Lemma~\ref{lem:mu_Dn}, Equation~\eqref{incond} holds true for $n$ if and only if
\begin{eqnarray*}
{\bf 0}
&=&
{}_{D^n}{\bf P}_{-}
(|{\bf R}_{-}|^{-1}{\bf T}_{- +})
-{}_{D^n}\bnu_+(0)(|{\bf R}_+|^{-1}{\bf C}_+)
\nonumber\\
&=&	{}_{D^{n+1}} {\bf P}_+
\nonumber\\
&=&
\left[
{\bf P}(|{\bf R}|^{-1}{\bf T})^{n+1}
+\sum_{k=1}^{n+1} \bnu^{(k-1)}(0){\bf h}(k,n+1)
\right]_+
\nonumber\\
&=&
{\bf P}_{-}
\left[
(|{\bf R}|^{-1}{\bf T})^{n+1}
\right]_{- +}
+\sum_{k=1}^{n+1} (\bnu_+)^{(k-1)}(0)
\left[
{\bf h}(k,n+1)
\right]_{++}
+\sum_{k=1}^n (\bnu_{-})^{(k-1)}(0)
\left[
{\bf h}(k,n+1)
\right]_{- +}
\nonumber\\
&=&
{\bf P}_{-}
\left[
(|{\bf R}|^{-1}{\bf T})^{n+1}
\right]_{- +}
+\sum_{k=1}^n (\bnu_+)^{(k-1)}(0)
\left[
{\bf h}(k,n+1)
\right]_{++}
+ (\bnu_+)^{(n)}(0)
\left[
{\bf h}(n+1,n+1)
\right]_{++}
\nonumber\\
&&
+\sum_{k=1}^n (\bnu_{-})^{(k-1)}(0)
\left[
{\bf h}(k,n+1)
\right]_{- +}.
\end{eqnarray*}
Thus, by noting that $\left[{\bf h}(n+1,n+1)\right]_{++}=\left( - |{\bf R}_+|^{-1}{\bf C}_+\right)^{n+1}$ and rearranging the above equation,
we have
		\begin{eqnarray*}
		(\bnu_+)^{(n)}(0)
		&=&
		-{\bf P}_{-}
		\left[
		(|{\bf R}|^{-1}{\bf T})^{n+1}
		\right]_{- +}
		\left(-|{\bf R}_+|^{-1}{\bf C}_+ \right)^{-(n+1)}
			\nonumber\\
			&&
		-\sum_{k=1}^n (\bnu_{-})^{(k-1)}(0)
		\left[
		{\bf h}(k,n+1)
		\right]_{- +}\left(-|{\bf R}_+|^{-1}{\bf C}_+ \right)^{-(n+1)}
		\nonumber\\
		&&
		-\sum_{k=1}^{n} (\bnu_+)^{(k-1)}(0)
		\left[
		{\bf h}(k,n+1)
		\right]_{+ +}\left(-|{\bf R}_+|^{-1}{\bf C}_+ \right)^{-(n+1)}
		\nonumber\\
		&=&
		{\bf A}(n+1)
		-
		\sum_{k=1}^n
		(\bnu_+)^{(k-1)}(0)
		[{\bf f}(k,n+1)]_{++}
		,
		\end{eqnarray*}
which proves~\eqref{Dn_nuplus0_rec}.

Next, we prove~\eqref{Dn_nuplus0} using mathematical induction.

{\bf Step(i).}\quad First, consider $n=1$. By~\eqref{Dn_nuplus0_rec} for $n=1$, we have
\begin{eqnarray*}
{\bnu}^{(1)}_+ (0)
&=&
{\bf A}(2)
-
\sum_{k=1}^1
(\bnu_+)^{(k-1)}(0)
[{\bf h}(k,2)]_{++}
(-|{\bf R}_+|^{-1}{\bf C}_+)^{-2}
\nonumber\\
&=&
{\bf A}(2)
-
{\bf A}(1)[{\bf f}(1,2)]_{++}
\nonumber\\
&=& {\bf A}(2)+{\bf B}(2),
\label{nuplus1_0}
\end{eqnarray*}
which completes the proof of~\eqref{Dn_nuplus0} for $n=1$.

{\bf Step(ii).}\quad Suppose now that~\eqref{Dn_nuplus0} holds for all $k=1,\ldots ,n-1$ for some $n\geq 2$. Then, by~\eqref{Dn_nuplus0_rec},
\begin{eqnarray*}
\bnu_+^{(n)}(0)
&=&
{\bf A}(n+1)
-
\sum_{k=1}^{n} (\bnu_+)^{(k-1)}(0)[{\bf f}(k,n+1)]_{+ +}
,
\end{eqnarray*}
and, by the inductive assumption,
\begin{eqnarray*}\label{Bnplus1}
\lefteqn{
	-
	\sum_{k=1}^{n}
	\bnu_+^{(k-1)}(0)
	[{\bf f}(k,n+1)]_{++}
	=
	-
	\sum_{k=1}^{n}
	\left(
	{\bf A}(k)+{\bf B}(k)
	\right)
	[{\bf f}(k,n+1)]_{++}
}
\nonumber\\
&=&
	-
	\sum_{k=1}^{n}
	\left(
	{\bf A}(k)+\sum_{r=1}^{k-1}{\bf M}(r,k)
	\right)
	[{\bf f}(k,n+1)]_{++}
	\nonumber\\
&=&
-
\sum_{k=1}^{n}
{\bf A}(k)
[{\bf f}(k,n+1)]_{++}
-
\sum_{k=1}^{n}
\sum_{r=1}^{k-1}{\bf M}(r,k)
[{\bf f}(k,n+1)]_{++}
\nonumber\\
&=&
{\bf M}(1,n+1)
-\sum_{r=1}^{n-1}\sum_{k=r+1}^n
{\bf M}(r,k)
[{\bf f}(k,n+1)]_{++}
\nonumber\\
&=&
{\bf M}(1,n+1)
+\sum_{s=2}^{(n+1)-1}{\bf M}(s,n+1)
\nonumber\\
&=&\sum_{s=1}^{(n+1)-1}{\bf M}(s,n+1)
\nonumber\\
&=&
{\bf B}(n+1).
\end{eqnarray*}
This completes the proof.
\hfill\rule{9pt}{9pt}

We are now ready to prove our first main result identifying the generator $D$ defined in Lemma~\ref{lemma_D}.
\begin{thm}\label{exp_distr}
Suppose the initial distribution $\bmu$ of $X(0)$ has a density given by
\begin{equation}\label{muexpo}
\nu_i(x)=p_i\lambda e^{-\lambda x},
\end{equation} for some $\lambda>0$, $0\leq p_i\leq 1$, such that the boundary conditions~\eqref{incond} are met.
Then, for any $y>0$,
	\begin{eqnarray*}
	\bmu e^{Dy}(\mathcal{A}_v)&=&
	-\bmu (\bar{\mathcal{A}_v})
	e^{(|{\bf R}|^{-1}{\bf T}+\lambda |{\bf R}|^{-1}{\bf C})y}
	+\bmu ([ 0,\infty))
	e^{(|{\bf R}|^{-1}{\bf T})y},
	\end{eqnarray*}
	where
	\begin{eqnarray*}
	\bmu (\bar{\mathcal{A}_v})&=&
	\bmu^{(0)}([0,\infty))
	-\bmu^{(0)} (\mathcal{A}_v)
	=e^{-\lambda v}\frac{\bnu(0)}{\lambda}
	\end{eqnarray*}
is the initial ditribution of starting outside set $\mathcal{A}_v$, and so in the set $\bar{\mathcal{A}_v}=(v,+\infty)$.	
\end{thm}

\begin{Remark}
	\label{rem:meanings}
	We have the following physical interpretations of interest in applications. Note that
	\begin{eqnarray}
	\lim_{v\to\infty}\bmu e^{Dy}(\mathcal{A}_v)&=&
	\bmu e^{Dy}([ 0,\infty))=
	\bmu ([ 0,\infty))e^{(|{\bf R}|^{-1}{\bf T})y}
	,
	\end{eqnarray}
	as expected, since $\bmu ([ 0,\infty))$ is the initial distribution of the process $\varphi(\cdot)$ and $\bmu ([ 0,\infty))e^{(|{\bf R}|^{-1}{\bf T})y}$ is the distribution of the process $\varphi(\cdot)$ at time $\omega(y)$, by~\cite{BOT2,BOT,BOT5}. Further,
	\begin{eqnarray*}
	\lim_{v\to 0}\bmu e^{Dy}(\mathcal{A}_v)&=&
	\left[
	\sum_{i\in\mathcal{S}^{\ell},\ell\in\{+,-\}}
	\int_{x=0}^{\infty}d\mu_i^{\ell}(x)\mathbb{P}[\varphi(\omega(y))=j,X(\omega(y))=0\ |\ \varphi(0)=i,X(0)=x]\right]_{j\in\mathcal{S}^+\cup\mathcal{S}^-}
	\nonumber\\
	&=&
	-\bmu^{(0)} ([ 0,\infty))
	e^{(|{\bf R}|^{-1}{\bf T}+\lambda |{\bf R}|^{-1}{\bf C})y}
	+\bmu ([ 0,\infty))
	e^{(|{\bf R}|^{-1}{\bf T})y}
	\nonumber\\
	&=&
	-\frac{\bnu (0)}{\lambda}
	e^{(|{\bf R}|^{-1}{\bf T}+\lambda |{\bf R}|^{-1}{\bf C})y}
	+\bmu ([ 0,\infty))
	e^{(|{\bf R}|^{-1}{\bf T})y}
	\end{eqnarray*}
	is the distribution of $\varphi(\cdot)$ when $X=0$ at time $\omega(y)$;
	\begin{eqnarray*}
	\lefteqn{
	\lim_{v\to\infty}\bmu e^{Dy}(\mathcal{A}_v)
	-\lim_{v\to 0}\bmu e^{Dy}(\mathcal{A}_v)
}
	\nonumber\\
	&=&
	\left[
	\sum_{i\in\mathcal{S}^{\ell},\ell\in\{+,-\}}
	\int_{x=0}^{\infty}d\mu_i^{\ell}(x)\mathbb{P}[\varphi(\omega(y))=j,X(\omega(y))>0\ |\ \varphi(0)=i,X(0)=x]\right]_{j\in\mathcal{S}^+\cup\mathcal{S}^-}
	\nonumber\\
	&=&
	\frac{\bnu (0)}{\lambda}
	e^{(|{\bf R}|^{-1}{\bf T}+\lambda |{\bf R}|^{-1}{\bf C})y}
	\end{eqnarray*}
	is the distribution of $\varphi(\cdot)$ when $X>0$ at time $\omega(y)$; and
\begin{eqnarray*}
\lefteqn{
\lim_{y\to\infty}\bmu e^{Dy}(\mathcal{A}_v)
}
\nonumber\\
&=&
\left[
\sum_{i\in\mathcal{S}^{\ell},\ell\in\{+,-\}}
\lim_{y\to\infty}\int_{x=0}^{\infty}d\mu_i^{\ell}(x)\mathbb{P}
[\varphi(\omega(y))=j,X(\omega(y))\leq v\ |\ \varphi(0)=i,X(0)=x]\right]_{j\in\mathcal{S}^+\cup\mathcal{S}^-}
\nonumber\\
&=&
\left[
\sum_{i\in\mathcal{S}^{\ell},\ell\in\{+,-\}}
\lim_{t\to\infty}\int_{x=0}^{\infty}d\mu_i^{\ell}(x)\mathbb{P}
[\varphi(t)=i,X(t)\leq v\ |\ \varphi(0)=j,X(0)=x]\right]_{j\in\mathcal{S}^+\cup\mathcal{S}^-}
\end{eqnarray*}
is the stationary distribution of the process $\{(\varphi(t),X(t)):t\geq 0\}$.
\end{Remark}
{\bf Proof:}\quad Due to the assumption~\eqref{muexpo}, we have $\bnu^{(k)}(x)=(-\lambda)^k\bnu(x)$ and $\bmu^{(k)} (\mathcal{A}_v)=(-\lambda)^k\bmu^{(0)}(\mathcal{A}_v)$. This key property of the initial distribution $\bmu=[\mu_j]_{j\in\mathcal{S}}$ of $X(0)$ results in the the following analysis. We have,
\begin{eqnarray*}
	{}_{D^n}\bmu (\mathcal{A}_v)
	&=&
	\sum_{k=0}^n (-1)^k {\bmu}^{(k)}(\mathcal{A}_v) {\bf h}(k,n)
	+
	{\bf P}(|{\bf R}|^{-1}{\bf T})^n+\sum_{k=1}^n (-1)^k\bnu^{(k-1)}(0){\bf h}(k,n)
	\nonumber\\
	&=&
	\sum_{k=0}^n (-1)^k(-\lambda)^k \bmu^{(0)}(\mathcal{A}_v){\bf h}(k,n)
	+
	{\bf P}(|{\bf R}|^{-1}{\bf T})^n+\sum_{k=1}^n (-1)^k(-\lambda)^{k-1}\bnu(0){\bf h}(k,n)
	\nonumber\\
	&=&
	\bmu^{(0)}(\mathcal{A}_v)\sum_{k=0}^n \lambda^k {\bf h}(k,n)
	+
	{\bf P}(|{\bf R}|^{-1}{\bf T})^n
	-\frac{\bnu(0)}{\lambda}\sum_{k=1}^n \lambda^k{\bf h}(k,n)
	\nonumber\\
	&=&
	\bmu^{(0)}(\mathcal{A}_v)\sum_{k=0}^n \lambda^k {\bf h}(k,n)
	+
	{\bf P}(|{\bf R}|^{-1}{\bf T})^n
	-\frac{\bnu(0)}{\lambda}\sum_{k=0}^n \lambda^k{\bf h}(k,n)
	+\frac{\bnu(0)}{\lambda}{\bf h}(0,n)
	\nonumber\\			
	&=&
	\bmu^{(0)}(\mathcal{A}_v)(|{\bf R}|^{-1}{\bf T}+\lambda |{\bf R}|^{-1}{\bf C})^n
	+{\bf P}(|{\bf R}|^{-1}{\bf T})^n
	-\frac{\bnu(0)}{\lambda}(|{\bf R}|^{-1}{\bf T}+\lambda |{\bf R}|^{-1}{\bf C})^n
	\nonumber\\
	&&
	+\frac{\bnu(0)}{\lambda}(|{\bf R}|^{-1}{\bf T})^n
	\nonumber\\
	&=&	
	\left(\bmu^{(0)}(\mathcal{A}_v)-\frac{\bnu(0)}{\lambda}\right)
	(|{\bf R}|^{-1}{\bf T}+\lambda |{\bf R}|^{-1}{\bf C})^n
	+\left({\bf P}+\frac{\bnu(0)}{\lambda} \right)(|{\bf R}|^{-1}{\bf T})^n,
\end{eqnarray*}
and so
\begin{eqnarray}
\bmu e^{D y}(\mathcal{A}_v)&=&
\sum_{n=0}^{\infty}\frac{y^n}{n!}{}_{D^n}\bmu (\mathcal{A}_v)
\nonumber\\
	&=&
	\left(\bmu^{(0)}(\mathcal{A}_v)-\frac{\bnu(0)}{\lambda}\right)
	e^{(|{\bf R}|^{-1}{\bf T}+\lambda |{\bf R}|^{-1}{\bf C})y}
	+\left({\bf P}+\frac{\bnu(0)}{\lambda} \right)
	e^{(|{\bf R}|^{-1}{\bf T})y}.
\end{eqnarray}

Then the result follows, since
\begin{eqnarray*}
	\bmu (\bar{\mathcal{A}_v})
	&=&
	\int_{x=0}^{\infty}
	\bnu(0)
	e^{-\lambda x}dx
	-
	\int_{x=0}^{v}
	\bnu(0)
	e^{-\lambda x}dx
	=
	\frac{\bnu(0)}{\lambda}
	-\bmu^{(0)}(\mathcal{A}_v)
	=
	e^{-\lambda v}\frac{\bnu(0)}{\lambda}
\end{eqnarray*}
and
\begin{eqnarray*}
	\bmu ([ 0,\infty))&=&
	\int_{x=0}^{\infty}
	\bnu(0)
	e^{-\lambda x}dx
	+{\bf P}
	=
	\left({\bf P}+\frac{\bnu(0)}{\lambda} \right).
\end{eqnarray*}
\hfill\rule{9pt}{9pt}

Recall that the main result stated in Theorem \ref{exp_distr} holds under assumed boundary conditions~\eqref{incond}. In the next theorem below,  we give sufficient conditions for these boundary conditions to hold.

\begin{thm}\label{lem_ex}
Suppose that $\gamma |r_i|=|c_i|$ for all $i$, for some $\gamma>0$. Let $(|{\bf R}_+|^{-1}{\bf T}_{+-})$ be a nonnegative matrix such that $(|{\bf R}_+|^{-1}{\bf T}_{+-}){\bf 1}=(b+\beta){\bf 1}$ for some $b,\beta >0$, and $(|{\bf R}_{-}|^{-1}{\bf T}_{- +})$ be a nonnegative matrix such that $(|{\bf R}_{-}|^{-1}{\bf T}_{- +}){\bf 1}=b{\bf 1}$. Assume that
\begin{eqnarray*}
|{\bf R}|^{-1}{\bf T}&=&
\left[
\begin{array}{cc}
(|{\bf R}_+|^{-1}{\bf T}_{++})&(|{\bf R}_+|^{-1}{\bf T}_{+-})\\
(|{\bf R}_{-}|^{-1}{\bf T}_{- +})&(|{\bf R}_{-}|^{-1}{\bf T}_{- -})
\end{array}
\right]
=
\left[
\begin{array}{cc}
-(b+\beta){\bf I}&(|{\bf R}_+|^{-1}{\bf T}_{+-})\\
(|{\bf R}_{-}|^{-1}{\bf T}_{- +})&-b{\bf I}
\end{array}
\right]
,
\end{eqnarray*}
and that the initial distribution is given by
\begin{eqnarray*}
{\bf P}_{-}&=&\mbox{ is such that }
{\bf P}_{-}\geq {\bf 0},{\bf P}_{-}{\bf 1}\leq \frac{\lambda\gamma}{b+\lambda\gamma},
\\
\bnu_+(0)&=&{\bf P}_{-}(|{\bf R}_{-}|^{-1}{\bf T}_{- +})
(|{\bf R}_+|^{-1}{\bf C}_+)^{-1}
,
\\
\bnu_+(x)&=&e^{-\lambda x}\bnu_+(0),\ x>0,
\\
\bnu_{-}(0)&=&\mbox{ is such that }
\bnu_{-}(0)\geq {\bf 0},\bnu_{-}(0){\bf 1}
=\lambda -\frac{b+\lambda\gamma}{\gamma}{\bf P}_{-}{\bf 1}
,\\
\bnu_{-}(x)&=&e^{-\lambda x}\bnu_{-}(0),\ x>0,
\end{eqnarray*}
where $\lambda=\beta/\gamma$. Then the boundary conditions~\eqref{incond} are satisfied for all $n\geq 0$.	
\end{thm}

\begin{Remark}
	Under the assumptions of Theorem~\ref{lem_ex}, the process $\{(\varphi(t),X(t)):t\geq 0\}$ is stable due to the structure of its fluid generator $|{\bf C}|^{-1}{\bf T}=|{\bf R}|^{-1}{\bf T}/\gamma$, which by~\cite{BOT2,BOT5} implies that stability condition given by $\sum_{i\in\mathcal{S}_+\cup\mathcal{S}_-} c_i\pi_i<0$, $\pi_i=\lim_{t\to\infty}\mathbb{P}(\varphi(t)=i)$, is met. The stability condition for $\{(\varphi(t),Y(t)):t\geq 0\}$ given by $\sum_{i\in\mathcal{S}^+\cup\mathcal{S}^-} r_i\pi_i<0$, may be met or not, depending on the parameters $(|{\bf R}_+|^{-1}{\bf T}_{+-})$ and $(|{\bf R}_{-}|^{-1}{\bf T}_{- +})$ in its fluid generator $|{\bf R}|^{-1}{\bf T}$.
\end{Remark}
{\bf Proof:}\quad We have $\bnu_{-}(0){\bf 1}\geq 0$ since ${\bf P}_{-}{\bf 1}\leq \frac{\lambda\gamma}{b+\lambda\gamma}$, and
\begin{eqnarray*}
\lefteqn{\int_0^{\infty}\bnu_+(x)dx{\bf 1}+\int_0^{\infty}\bnu_{-}(x)dx{\bf 1}
+{\bf P}_{-}{\bf 1}
=\frac{1}{\lambda}\bnu_+(0){\bf 1}
+\frac{1}{\lambda}\bnu_{-}(0){\bf 1}+{\bf P}_{-}{\bf 1}
}
\nonumber\\
&=&\frac{1}{\lambda}{\bf P}_{-}(|{\bf R}_{-}|^{-1}{\bf T}_{- +})
(|{\bf R}_+|^{-1}{\bf C}_+)^{-1}{\bf 1}
+\frac{1}{\lambda}
\left(
\lambda -\frac{b+\lambda\gamma}{\gamma}{\bf P}_{-}{\bf 1}
\right)
+{\bf P}_{-}{\bf 1}
\nonumber\\
&=&
\frac{b}{\lambda\gamma}{\bf P}_{-}{\bf 1}
+
\left(1-\frac{b+\lambda\gamma}{\lambda\gamma}{\bf P}_{-}{\bf 1}\right)
+\frac{\lambda\gamma}{\lambda\gamma}{\bf P}_{-}{\bf 1}
\nonumber\\
&=&1,
\end{eqnarray*}
and so the initial distribution is well defined.

The boundary condition~\eqref{incond} is clearly met for $n= 0$. We verify that the boundary conditions~\eqref{incond} are met for all $n\geq 0$, by mathematical induction. Suppose the boundary condition~\eqref{incond} is met for some $n\geq 0$ (inductive assumption). Then by~\eqref{D2_dens},
\begin{eqnarray*}
{}_{D^{n+1}}{\bnu} (x) &=& {}_{D^n}{\bnu} (x) (|{\bf R}|^{-1}{\bf T})
-  {}_{D^n}{\bnu}^{(1)} (x)
(|{\bf R}|^{-1}{\bf C})
\nonumber\\
&=&
{}_{D^n}{\bnu} (x) \left( |{\bf R}|^{-1}{\bf T}+\lambda |{\bf R}|^{-1}{\bf C} \right)
\nonumber\\
&=&
\left[
\begin{array}{cc}
{}_{D^n}{\bnu}_+ (x)&{}_{D^n}{\bnu}_{-} (x)
\end{array}
\right]
\left[
\begin{array}{cc}
-b{\bf I}&(|{\bf R}_+|^{-1}{\bf T}_{+-})\\
(|{\bf R}_{-}|^{-1}{\bf T}_{- +})&-(b+\beta){\bf I}
\end{array}
\right]
\end{eqnarray*}
and so
\begin{eqnarray*}
{}_{D^{n+1}}{\bnu}_+ (0)&=&
{}_{D^{n}}{\bnu}_+ (0)(-b)+{}_{D^{n}}{\bnu}_{-} (0)(|{\bf R}_{-}|^{-1}{\bf T}_{- +}).
\end{eqnarray*}
Thus, by inductive assumption,
\begin{eqnarray*}
{}_{D^{n+1}}{\bnu}_+ (0)&=&{}_{D^{n}}{\bnu}_+ (0)(-b)+
{}_{D^{n}}{\bnu}_{-} (0)(|{\bf R}_{-}|^{-1}{\bf T}_{- +})
\nonumber\\
&=&
{}_{D^n}{\bf P}_{-}(|{\bf R}_{-}|^{-1}{\bf T}_{- +})(|{\bf R}_+|^{-1}{\bf C}_+)^{-1}(-b)+
{}_{D^{n}}{\bnu}_{-} (0)(|{\bf R}_{-}|^{-1}{\bf T}_{- +})
\nonumber\\
&=&
{}_{D^n}{\bf P}_{-}(|{\bf R}_{-}|^{-1}{\bf T}_{- +})\gamma^{-1}(-b)
+{}_{D^{n}}{\bnu}_{-} (0)(|{\bf R}_{-}|^{-1}{\bf T}_{- +}),
\end{eqnarray*}
and, by~\eqref{DnplusB}, we have
\begin{eqnarray*}
{}_{D^{n+1}} {\bf P}_{-}(|{\bf R}_{-}|^{-1}{\bf T}_{- +})(|{\bf R}_+|^{-1}{\bf C}_+)^{-1} &=&
\left(
{}_{D^n}{\bf P}_{-} (|{\bf R}_{-}|^{-1}{\bf T}_{--})
- {}_{D^n}{\bnu}_{-} (0) (|{\bf R}_{-}|^{-1}{\bf C}_{-})
\right)
\nonumber\\
&&\times
(|{\bf R}_{-}|^{-1}{\bf T}_{- +})(|{\bf R}_+|^{-1}{\bf C}_+)^{-1}
\nonumber\\
&=&
\left(
{}_{D^n}{\bf P}_{-}(-b{\bf I})
-{}_{D^n}{\bnu}_{-} (0)(-\gamma{\bf I})
\right)(|{\bf R}_{-}|^{-1}{\bf T}_{- +})\gamma^{-1}
\nonumber\\
&=&{}_{D^{n+1}}{\bnu}_+ (0),
\end{eqnarray*}
which implies that the boundary condition~\eqref{incond} is met for $(n+1)$. \hfill\rule{9pt}{9pt}

\subsection{Numerical examples}\label{sec:numeex1}

Below, we construct simple examples to illustrate the application of the main result established above to the SFFMs with a range of different behaviours. In Examples~\ref{ex2}-\ref{ex3} we consider SFFMs in which the sign of the fluid rates is the same or the opposite in both fluids, respectively. That is, there are two possible directions of the movement in the quadrant depicted in Figure~\ref{fig:SFFMquadrant}, in each of these SFFMs. Next, in Examples~\ref{ex4}-\ref{ex5} we consider SFFMs with four possible directions of the movement in the quadrant.

Also, different behaviours as far as the stability of the fluids are presented in these examples. In Example~\ref{ex2} both fluids are stable, while in Example~\ref{ex3} only $Y(\cdot)$ is stable. In Example~\ref{ex4} process $X(\cdot)$ is stable, while $Y(\cdot)$ is null recurrent. In Example~\ref{ex4} both fluids are stable.

\begin{exa}\label{ex2}
\rm Assume $\mathcal{S}=\{1,2\}$, $\mathcal{S}_+=\mathcal{S}^+=\{1\}$, $\mathcal{S}_-=\mathcal{S}^-=\{2\}$, and so the signs of the fluid rates are the same in both fluids.
\begin{figure}
	\begin{center}
		\begin{tikzpicture}[>=stealth,redarr/.style={->}]
		\node at (11,6) [black,circle,fill,inner sep=1.3pt]{};
		\draw [black, dashed] (11,6) -- (11.5,6);
		\draw [black, dashed] (11,6) -- (10.5,6);
		\draw [black, dashed] (11,6) -- (11,6.5);
		\draw [black, dashed] (11,6) -- (11,5.5);
		\draw [red,thick, ->] (11,6) -- (11.5,6.5);
		\draw [red,thick, ->] (11,6) -- (10.5,5.5);
		\draw [black, dashed] (11,6) -- (11.5,5.5);
		\draw [black, dashed] (11,6) -- (10.5,6.5);
		
		\draw (11.6,6.6) node[anchor=north, below=-0.17cm] {\tiny{\color{red} $1$}};
		\draw (10.4,5.4) node[anchor=north, below=-0.17cm] {\tiny{\color{red} $2$}};
		\end{tikzpicture}
		\caption{\footnotesize Evolution of the SFFM $\{(\varphi(t),X(t),Y(t)):t\geq 0\}$ in Example~\ref{ex2}. When $\varphi(t)=1$ both levels are increasing. When $\varphi(t)=2$ both levels are decreasing . No other directions are possible.}
	\end{center}
\end{figure}
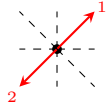

Assume $|r_i|=|c_i|$ for $i=1,2$. Let $\beta>0$ and
\begin{eqnarray*}
|{\bf R}|^{-1}{\bf T}&=&
\left[
\begin{array}{cc}
-(b+\beta)&b+\beta\\b&-b
\end{array}
\right]
=
|{\bf C}|^{-1}{\bf T}.
\end{eqnarray*}
Assume initial distribution given by,
\begin{eqnarray*}
{\bf P}_{-}&=&p \mbox{ such that }0<p<\frac{\beta}{b+\beta},
\nonumber\\
\bnu_+(x)&=&\bnu_+(0)e^{-\beta x}
,\ x>0,
\ \bnu_+(0)=pb,
\nonumber\\
\bnu_{-}(x)&=&\bnu_{-}(0)e^{-\beta x}
,\ x>0,
\ \bnu_{-}(0)=\beta(1-p)-pb
>0 \mbox{ since }p<\frac{\beta}{b+\beta}
.
\end{eqnarray*}
Then both $\{(\varphi(t),X(t)):t\geq 0\}$ and $\{(\varphi(t),Y(t)):t\geq 0\}$ are stable due to the structure of the fluid generator $|{\bf R}|^{-1}{\bf T}=
|{\bf C}|^{-1}{\bf T}$. Moreover, by Theorem~\ref{lem_ex}, the initial distribution is well defined and the boundary conditions~\eqref{incond} are met for all $n\geq 0$.

Further, by some calculations involving standard algebra, we have,
\begin{eqnarray*}
|{\bf R}|^{-1}{\bf T}
&=&
{\bf B}{\bf D}{\bf B}^{-1}
\nonumber\\
&=&
\left[
\begin{array}{cc}
1 & 1\\
1 & \frac{-b}{b+\beta}
\end{array}
\right]
\left(
\left(
-(2b+\beta)
\right)
\left[
\begin{array}{cc}
0 & 0\\
0 & 1
\end{array}
\right]
\right)
\left[
\begin{array}{cc}
1 & 1\\
1 & \frac{-b}{b+\beta}
\end{array}
\right]^{-1}
\nonumber\\
&=&
\left[
\begin{array}{cc}
1 & 1\\
1 & \frac{-b}{b+\beta}
\end{array}
\right]
\left(
\left(
-(2b+\beta)
\right)
\left[
\begin{array}{cc}
0 & 0\\
0 & 1
\end{array}
\right]
\right)
\left(
\frac{b+\beta}{2b+\beta}
\left[
\begin{array}{cc}
1 & 1\\
1 & \frac{-b}{b+\beta}
\end{array}
\right]
\right),
\end{eqnarray*}
and thus
\begin{eqnarray*}
(|{\bf R}|^{-1}{\bf T})^{n+1}
&=&
{\bf B} {\bf D}^{n+1} {\bf B}^{-1}
\nonumber\\
&=&
\left[
\begin{array}{cc}
-(b+\beta) & b+\beta\\
b & -b
\end{array}
\right]
(-(2b+\beta))^n.
\end{eqnarray*}
Similarly, with ${\bf\Lambda}=\beta{\bf I}$, we have
\begin{eqnarray*}
|{\bf R}|^{-1}{\bf T}+{\bf \Lambda}|{\bf R}|^{-1} {\bf C}
&=&
\left[
\begin{array}{cc}
-b&b+\beta\\b&-(b+\beta)
\end{array}
\right]
\nonumber\\
&=&
\widehat{\bf B}{\bf D}{\widehat{\bf B}}^{-1}
\nonumber\\
&=&
\left[
\begin{array}{cc}
1 & -1\\
\frac{b}{b+\beta} & 1
\end{array}
\right]
\left(
\left(
-(2b+\beta)
\right)
\left[
\begin{array}{cc}
0 & 0\\
0 & 1
\end{array}
\right]
\right)
\left[
\begin{array}{cc}
1 & -1\\
\frac{b}{b+\beta} & 1
\end{array}
\right]^{-1}
\nonumber\\
&=&
\left[
\begin{array}{cc}
1 & -1\\
\frac{b}{b+\beta} & 1
\end{array}
\right]
\left(
\left(
-(2b+\beta)
\right)
\left[
\begin{array}{cc}
0 & 0\\
0 & 1
\end{array}
\right]
\right)\left(
\frac{b+\beta}{2b+\beta}
\left[
\begin{array}{cc}
1 & 1\\
\frac{-b}{b+\beta} & 1
\end{array}
\right]
\right),
\end{eqnarray*}
and so
\begin{eqnarray*}
(|{\bf R}|^{-1}{\bf T}+{\bf \Lambda}|{\bf R}|^{-1}{\bf C})^{n+1}
&=&
\widehat{\bf B} {\bf D}^{n+1} {\widehat{\bf B}}^{-1}
\nonumber\\
&=&
\left[
\begin{array}{cc}
-b & b+\beta\\
b & -(b+\beta)
\end{array}
\right]
(-(2b+\beta))^n,
\end{eqnarray*}
where ${\bf D}$ records eigenvalues and ${\bf B}$, $\widehat{\bf B}$ corresponding eigenvectors of $|{\bf R}|^{-1}{\bf T}$
and $|{\bf R}|^{-1}{\bf T}+{\bf \Lambda}|{\bf R}|^{-1} {\bf C}$, with
\begin{eqnarray*}
{\bf D}=diag(0,-(2b+\beta)),\
{\bf B}=\left[
\begin{array}{cc}
1 & 1\\
1 & \frac{-b}{b+\beta}
\end{array}
\right],\
\widehat{\bf B}=\left[
\begin{array}{cc}
1 & -1\\
\frac{b}{b+\beta} & 1
\end{array}
\right].
\end{eqnarray*}
Then by Theorem~\ref{exp_distr},
we have,
\begin{eqnarray*}
	\bmu e^{Dy}(\mathcal{A}_v)&=&
-e^{-\beta v}
\frac{\bnu(0)}{\beta}
e^{(|{\bf R}|^{-1}{\bf T}+\beta |{\bf R}|^{-1}{\bf C})y}
+\left({\bf P}+\frac{\bnu(0)}{\beta} \right)
e^{(|{\bf R}|^{-1}{\bf T})y}
,\nonumber
\end{eqnarray*}
where
\begin{eqnarray*}
-e^{-\beta v}\frac{\bnu(0)}{\beta}
&=&
-e^{-\beta v}\left[
\begin{array}{cc}
	p\frac{b}{\beta} & \left((1-p)-p\frac{b}{\beta}\right)
\end{array}
\right]
,
\\
\left({\bf P}+\frac{\bnu(0)}{\beta} \right)
&=&
\left[
\begin{array}{cc}
	p\frac{b}{\beta} & \left(1-p\frac{b}{\beta}\right)
\end{array}
\right]
,
\\
e^{(|{\bf R}|^{-1}{\bf T}+\beta |{\bf R}|^{-1}{\bf C})y}
&=&
\left[
\begin{array}{cc}
-b & b+\beta\\
b & -(b+\beta)
\end{array}
\right]
(1+e^{-(2b+\beta)y})
,
\\
e^{|{\bf R}|^{-1}{\bf T}y}
&=&
\left[
\begin{array}{cc}
-(b+\beta) & b+\beta\\
b & -b
\end{array}
\right]
(1+e^{-(2b+\beta)y})
,\nonumber
\end{eqnarray*}
and so
\begin{eqnarray*}
\bmu e^{Dy}(\mathcal{A}_v)&=&
-e^{-\beta v}
\left[
\begin{array}{cc}
	\left(
	(1-p)b-p\frac{2b^2}{\beta}
	\right)	
	&
	\left(
	-(1-p)(b+\beta)+p\frac{2b^2}{\beta}
	\right)	
\end{array}
\right]
(1+e^{-(2b+\beta)y})
\\
&&
\quad
+
\left[
\begin{array}{cc}
	\left(
	1-pb-p\frac{2b^2}{\beta}
	\right)
	&
	\left(
	1+pb+p\frac{2b^2}{\beta}
	\right)
\end{array}
\right]
(1+e^{-(2b+\beta)y}),
\end{eqnarray*}
where $[\bmu e^{Dy}(\mathcal{A}_v)]_j$ is the probability that at time $\omega(y)$ the observed phase is $j$ and the level $X(\omega(y))$ is within the set $\mathcal{A}_v$, also see Figure~\ref{fig:Ut}.

\end{exa}

\begin{exa}\label{ex3}
\rm We modify Example~\ref{ex2} and assume that $\mathcal{S}=\{1,2\}$, $\mathcal{S}_+=\mathcal{S}^-=\{1\}$, $\mathcal{S}_-=\mathcal{S}^+=\{2\}$, and so the signs of the fluid rates are now opposite in both fluids.
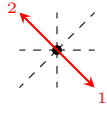
\begin{figure}
	\begin{center}
		\begin{tikzpicture}[>=stealth,redarr/.style={->}]
		\node at (11,6) [black,circle,fill,inner sep=1.3pt]{};
		\draw [black, dashed] (11,6) -- (11.5,6);
		\draw [black, dashed] (11,6) -- (10.5,6);
		\draw [black, dashed] (11,6) -- (11,6.5);
		\draw [black, dashed] (11,6) -- (11,5.5);
		\draw [black, dashed] (11,6) -- (11.5,6.5);
		\draw [black, dashed] (11,6) -- (10.5,5.5);
		\draw [red,thick, ->] (11,6) -- (11.5,5.5);
		\draw [red,thick, ->] (11,6) -- (10.5,6.5);
		
		\draw (11.6,5.4) node[anchor=north, below=-0.17cm] {\tiny{\color{red} $1$}};
		\draw (10.4,6.6) node[anchor=north, below=-0.17cm] {\tiny{\color{red} $2$}};
		\end{tikzpicture}
		\caption{\footnotesize Evolution of the SFFM $\{(\varphi(t),X(t),Y(t)):t\geq 0\}$ in Example~\ref{ex3}. When $\varphi(t)=1$ the level in $X$ is increasing while the level in $Y$ is decreasing. When $\varphi(t)=2$ the level in $X$ is decreasing while the level in $Y$ is increasing. No other directions are possible.}
	\end{center}
\end{figure}

Note that the analysis in Example~\ref{ex2} still holds. Indeed, the signs of the rates $r_i$ do not change the results as we take the absolute values of these rates in all expressions. Due to the structure of the fluid generator ${\bf RT}=|{\bf C}|^{-1}{\bf T}$, $\{(\varphi(t),Y(t)):t\geq 0\}$ is stable while $\{(\varphi(t),X(t)):t\geq 0\}$ is unstable.
\end{exa}

\begin{exa}\label{ex4}
\rm Assume that $\mathcal{S}=\{1,2,3,4\}$, $\mathcal{S}_+=\{1,2\}$, $\mathcal{S}_-=\{3,4\}$, $\mathcal{S}^+=\{1,4\}$, $\mathcal{S}^-=\{2,3\}$ so that signs are the same for phases $1$ and $3$, and opposite for $2$ and $4$.
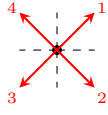
\begin{figure}
	\begin{center}
		\begin{tikzpicture}[>=stealth,redarr/.style={->}]
		\node at (11,6) [black,circle,fill,inner sep=1.3pt]{};
		\draw [black, dashed] (11,6) -- (11.5,6);
		\draw [black, dashed] (11,6) -- (10.5,6);
		\draw [black, dashed] (11,6) -- (11,6.5);
		\draw [black, dashed] (11,6) -- (11,5.5);
		\draw [red,thick, ->] (11,6) -- (11.5,6.5);
		\draw [red,thick, ->] (11,6) -- (10.5,5.5);
		\draw [red,thick, ->] (11,6) -- (11.5,5.5);
		\draw [red,thick, ->] (11,6) -- (10.5,6.5);
		
		\draw (11.6,6.6) node[anchor=north, below=-0.17cm] {\tiny{\color{red} $1$}};
		\draw (10.4,5.4) node[anchor=north, below=-0.17cm] {\tiny{\color{red} $3$}};
		\draw (11.6,5.4) node[anchor=north, below=-0.17cm] {\tiny{\color{red} $2$}};
		\draw (10.4,6.6) node[anchor=north, below=-0.17cm] {\tiny{\color{red} $4$}};
		\end{tikzpicture}
	\end{center}
	\caption{\footnotesize
		Evolution of the SFFM $\{(\varphi(t),X(t),Y(t)):t\geq 0\}$ in Examples~\ref{ex4}--\ref{ex5}. When $\varphi(t)=1,3$ the levels in $X$ and $Y$ move in the same directions (increasing or decreasing, respectively). When $\varphi(t)=2,4$ the levels in $X$ and $Y$ move in the opposite directions (increasing-decreasing or decreasing-increasing, respectively). No other directions are possible.}
\end{figure}

Further, let $|r_i|=|c_i|=1$ for all $i$, and, with ${\bf E}$ denoting a matrix of ones of appropriate size, let
	\begin{eqnarray*}
	|{\bf R}|^{-1}{\bf T}&=&
	\left[
	\begin{array}{cc}
	-(b+\beta){\bf I}&((b+\beta)/2){\bf E}\\
	(b/2){\bf E}&-b{\bf I}
	\end{array}
	\right]
	={\bf T}
	=|{\bf C}|^{-1}{\bf T}
	,\\
	{\bf P}_{-}&=&(p/2)[1\ \ 1],\\
	\bnu_+(x)&=&(pb/2)e^{-\beta x}[1\ \ 1],
	\bnu_+(0)=(pb/2)[1\ \ 1],
	\\
	\bnu_{-}(x)&=&(\beta(1-p)/2 -pb/2)e^{-\beta x}[1\ \ 1],
	\bnu_{-}(0)=(\beta(1-p)/2 -pb/2)[1\ \ 1].
	\end{eqnarray*}
By Theorem~\ref{lem_ex}, the initial distribution is well defined and the boundary conditions~\eqref{incond} are met for all $n\geq 0$. Further, by some calculations involving standard algebra,
	\begin{eqnarray*}
	|{\bf R}|^{-1}{\bf T}&=&{\bf B}{\bf D}{\bf B}^{-1},
	(|{\bf R}|^{-1}{\bf T})^n={\bf B}{\bf D}^n{\bf B}^{-1}\\
	|{\bf R}|^{-1}{\bf T}+\beta|{\bf R}|^{-1}{\bf C}&=&\widehat{\bf B}{\bf D}\widehat{\bf B}^{-1},
	(|{\bf R}|^{-1}{\bf T}+\beta|{\bf R}|^{-1}{\bf C})^n=\widehat{\bf B}{\bf D}^n\widehat{\bf B}^{-1}
	\end{eqnarray*}	
	where ${\bf D}$ records eigenvalues and ${\bf B}$, $\widehat{\bf B}$ corresponding eigenvectors of $|{\bf R}|^{-1}{\bf T}$, $|{\bf R}|^{-1}{\bf T}+\beta|{\bf R}|^{-1}{\bf C}$, with
	\begin{eqnarray*}
	{\bf D}&=&diag(-b,-(b+\beta),0,-(2b+\beta)),\\
	{\bf B}&=&
	\left[
	\begin{array}{cccc}
	0&1&1&1\\
	0&-1&1&1\\
	1&0&1&-\frac{b}{b+\beta}\\
	-1&0&1&-\frac{b}{b+\beta}
	\end{array}	
	\right],\ \
	{\bf B}^{-1}=
	\left[
	\begin{array}{cccc}
	0&0&1/2&-1/2\\
	1/2&-1/2&0&0\\
	\frac{b}{2(2b+\beta)}&\frac{b}{2(2b+\beta)}
	&\frac{b+\beta}{2(2b+\beta)}&\frac{b+\beta}{2(2b+\beta)}\\
	\frac{b+\beta}{2(2b+\beta)}&\frac{b+\beta}{2(2b+\beta)}&
	-\frac{b+\beta}{2(2b+\beta)}&-\frac{b+\beta}{2(2b+\beta)}
	\end{array}	
	\right],
	\\
	\widehat{\bf B}&=&
	\left[
	\begin{array}{cccc}
	1&0&1&1\\
	-1&0&1&1\\
	0&1&1&-1\\
	0&-1&1&-1
	\end{array}	
	\right],\ \
	\widehat{\bf B}^{-1}=
	\left[
	\begin{array}{cccc}
	1/2&-1/2&0&0\\
	0&0&1/2&-1/2\\
	1/4&1/4&1/4&1/4\\
	1/4&1/4&-1/4&-1/4
	\end{array}	
	\right].
	\end{eqnarray*}	
Then, by Theorem~\ref{exp_distr}, $\bmu e^{Dy}(\mathcal{A}_v)$ recording the probabilities $[\bmu e^{Dy}(\mathcal{A}_v)]_j$ that at time $\omega(y)$ the observed phase is $j$ and the level $X(\omega(y))$ is within the set $\mathcal{A}_v$, is given by,
\begin{eqnarray*}
\bmu e^{Dy}(\mathcal{A}_v)&=&
-e^{-\beta v}
\frac{\bnu(0)}{\beta}
{\bf B}e^{{\bf D}y}{\bf B}^{-1}
+\left({\bf P}+\frac{\bnu(0)}{\beta} \right)
\widehat{\bf B}e^{{\bf D}y}\widehat{\bf B}^{-1}
,
\end{eqnarray*}
where,
\begin{eqnarray*}
	e^{{\bf D}y}&=&diag(e^{-by},e^{-(b+\beta)y},0,e^{-(2b+\beta)y}),
	\\
	-e^{-\beta v}\frac{\bnu(0)}{\beta}
	&=&
	-e^{-\beta v}
	\frac{1}{\beta}
	\left[
	\begin{array}{cccc}
		(pb/2)& (pb/2)& (\beta(1-p)/2 -pb/2)& (\beta(1-p)/2 -pb/2)
	\end{array}	
	\right]
	,
	\\
	\left({\bf P}+\frac{\bnu(0)}{\beta} \right)
	&=&
	\left[
	\begin{array}{cc}
		p\frac{b}{\beta} & \left(1-p\frac{b}{\beta}\right)
	\end{array}
	\right]
	,	
\end{eqnarray*}	
and so
\begin{eqnarray*}
-e^{-\beta v}
\frac{\bnu(0)}{\beta}
{\bf B}e^{{\bf D}y}{\bf B}^{-1}&=&
			-e^{-\beta v}\frac{1}{\beta}
			\left(
			pb -\frac{b}{b+\beta}
			\left(
			\beta(1-p)-pb
			\right)
			\right)
			e^{-(2b+\beta)y}
			\\
			&&
			\quad\times
			\left[
			\begin{array}{cccc}
				\frac{b+\beta}{2(2b+\beta)}
				&\frac{b+\beta}{2(2b+\beta)}
				&-
				\frac{b+\beta}{2(2b+\beta)}
				&-\frac{b+\beta}{2(2b+\beta)}
			\end{array}	
			\right],
\end{eqnarray*}
and
\begin{eqnarray*}
\left({\bf P}+\frac{\bnu(0)}{\beta} \right)
\widehat{\bf B}e^{{\bf D}y}\widehat{\bf B}^{-1}
&=&
	\frac{1}{\beta}
	\left[
	\begin{array}{c}
		e^{-by}	(pb/4)
		+
		e^{-(2b+\beta)y}(1/4)(\beta/2 -pb/2)
		\\[1ex]
		-e^{-by}	(pb/4)
		+
		e^{-(2b+\beta)y}(1/4)(\beta/2 -pb/2)
		\\[1ex]
		e^{-(b+\beta)y}(pb/4)
		-e^{-(b+\beta)y}(1/4)(\beta/2 -pb/2)
		\\[1ex]
		e^{-(b+\beta)y}(pb/4)
		-e^{-(b+\beta)y}(1/4)(\beta/2 -pb/2)
	\end{array}	
	\right]^T
.
\end{eqnarray*}

We note that here $\{(\varphi(t),X(t)):t\geq 0\}$ is stable with $\pi_1=\pi_2=\frac{b}{2(2b+\beta)}$, $\pi_3=\pi_4=\frac{b+\beta}{2(2b+\beta)}$, however $\{(\varphi(t),Y(t)):t\geq 0\}$ is null recurrent since $\pi_1+\pi_4=\pi_2+\pi_3$.

\end{exa}

\begin{exa}\label{ex5}
\rm We modify ${\bf T}=[q_{ij}]$ in Example~\ref{ex4} so that $\{(\varphi(t),Y(t)):t\geq 0\}$ is stable as well, in the following manner. For some $r>0.5$, we
\begin{enumerate}[(a)]
\item change rates to $q_{31}=q_{41}=(1-r)b$, $q_{32}=q_{42}=rb$, so that $\pi_3=\pi_4$ still but $\pi_1<\pi_2$; or
\item change rates to $q_{14}=q_{24}=r(b+\beta)$, $q_{13}=q_{23}=(1-r)(b+\beta)$, so that $\pi_1=\pi_2$ still but $\pi_4<\pi_3$.
\end{enumerate}
Consider (a). Let $\widetilde{\bf E}={\bf E}\times diag(1-r,r)$ and
\begin{eqnarray*}
|{\bf R}|^{-1}{\bf T}&=&
\left[
\begin{array}{cc}
-(b+\beta){\bf I}&((b+\beta)/2){\bf E}\\
b\widetilde{\bf E}&-b{\bf I}
\end{array}
\right]
={\bf T}
=|{\bf C}|^{-1}{\bf T}
,\\
{\bf P}_{-}&=&(p/2)[1\ \ 1],\\
\bnu_+(x)&=&(pb)e^{-\beta x}[1-r\ \ r],
\bnu_+(0)=(pb)[1-r\ \ r],
\\
\bnu_{-}(x)&=&(\beta(1-p)/2 -pb/2)e^{-\beta x}[1\ \ 1],
\bnu_{-}(0)=(\beta(1-p)/2 -pb/2)[1\ \ 1].
\end{eqnarray*}
By Theorem~\ref{lem_ex}, the initial distribution is well defined and the boundary conditions~\eqref{incond} are met for all $n\geq 0$. The stationary distribution of ${\bf T}$ is
\begin{eqnarray*}
\bpi
&=&
\left[
\begin{array}{cccc}
(1-r)\frac{b}{2b+\beta}&r\frac{b}{2b+\beta}
&\frac{1}{2}\frac{b}{2b+\beta}&\frac{1}{2}\frac{b}{2b+\beta}
\end{array}
\right].
\end{eqnarray*}
Since $\beta>0$ and $r>0.5$, we have
\begin{eqnarray*}
\pi_1+\pi_2=\frac{b}{2b+\beta}&<&
\pi_3+\pi_4=\frac{b+\beta}{2b+\beta},\\
\pi_1+\pi_4=\frac{(1-r)b+((b+\beta)/2)}{2b+\beta}&<&
\pi_2+\pi_3=\frac{rb+((b+\beta)/2)}{2b+\beta}.
\end{eqnarray*}	
Thus $\{(\varphi(t),X(t)):t\geq 0\}$ and $\{(\varphi(t),Y(t)):t\geq 0\}$ are stable. Stability for both processes in (b) follows by symmetry.
 	\end{exa}

\section{First return to level zero: measure $\bmu{\Phi}(\mathcal{A}_v)$}\label{sec:returnzero}

The distribution of the busy period is a key quantity analysed in the classic literature of the stochastic fluid models~\cite{1994Abusy,asmu95,2001BST,BOT,rama99} and in the queueing theory in general, and it is the key quantity required in the computation of other transient and stationary measures of interest. Therefore, below we consider the distribution of the process $X(\cdot)$ at the time the busy period ends in $Y(\cdot)$. In the context of the application example presented in the Introduction, we are interested in the distribution of level in buffer $X$ at the time when buffer $Y$ becomes empty.

	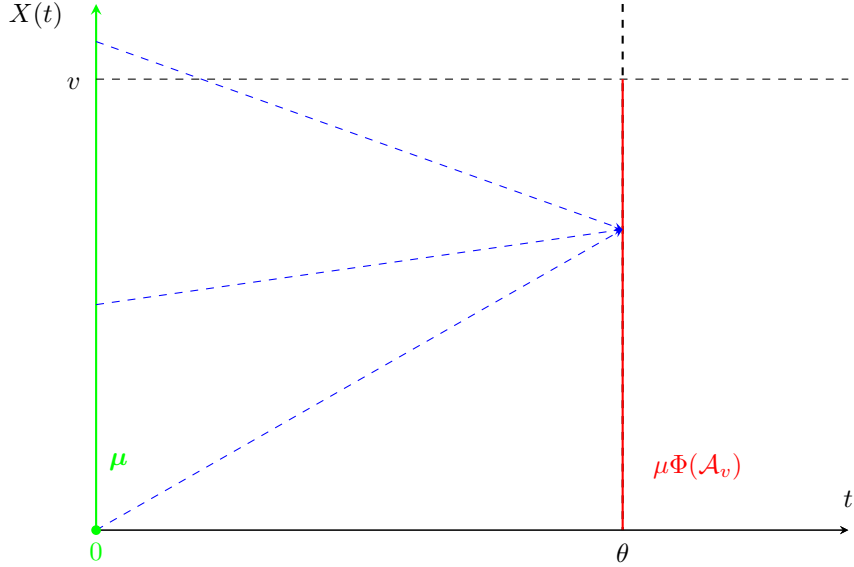
\begin{figure}
		\begin{center}
			\begin{tikzpicture}[>=stealth,redarr/.style={->}]

			\draw [thick, green, ->] (0,0) -- (0,7);
			\draw [ ->] (0,0) -- (10,0);
			\draw (-0.8,7) node[anchor=north, below=-0.17cm] {\footnotesize{\color{black} $X(t)$}};
			\draw (-0.3,6) node[anchor=north, below=-0.17cm] {\footnotesize{\color{black} $v$}};
			\draw (7,-0.2) node[anchor=north, below=-0.17cm] {\footnotesize{\color{black} $\theta$}};
			\draw (0,-0.2) node[anchor=north, below=-0.17cm] {\footnotesize{\color{green} $0$}};
			
			\draw (0.3,1) node[anchor=north, below=-0.17cm] {\footnotesize{\color{green} $\bmu$}};
			\draw (8,1) node[anchor=north, below=-0.17cm] {\footnotesize{\color{red} $\mu{\Phi}(\mathcal{A}_v)$}};
			
			\draw (10,0.5) node[anchor=north, below=-0.17cm] {\footnotesize{\color{black} $t$}};

			\draw [thick,dashed] (7,7) -- (7,0);
			
			\draw [dashed] (0,6) -- (10,6);
			
			\draw [red, thick] (7,6) -- (7,0);
			
			\draw [ ->] (0,0) -- (10,0);
			
			\draw [ blue, dashed,->] (0,0) -- (7,4);
			\draw [ blue, dashed,->] (0,3) -- (7,4);
			\draw [ blue, dashed,->] (0,6.5) -- (7,4);

			\node at (0,0) [green,circle,fill,inner sep=1.3pt]{};

			\end{tikzpicture}
		\end{center}
		\caption{\footnotesize Destination at time $\theta$: $\mathcal{A}_v=[0,v]$. The SFM $\{(\varphi(t),X(t)):t\geq 0\}$ starts in some phase $i\in\mathcal{S}$ and level $X(0)$ according to the initial distribution $\mu$. The SFM $\{(\varphi(t),\widetilde Y(t)):t\geq 0\}$ starts from level $0$ in phase $i$ and first returns to level $0$ at time $ \theta$ and does so in some phase $j\in \mathcal{S}$.}
		\label{fig:Phi}
	\end{figure}

Consider process $\{(\varphi(t),X(t),\widetilde Y(t)):t\geq 0\}$ with unbounded level $\widetilde Y(t)\in(-\infty,+\infty)$ and rates $d\widetilde Y(t)/dt=r_{\varphi(t)}$, defined in the Introduction. We introduce operators with physical interpretations similar to matrices ${\bf \Psi}$ and ${\bf \Xi}$ for the SFMs discussed in Section~\ref{sec:stop_SFM} (see Figure~\ref{fig:theta}).

Define operator $\Phi=[\Phi_{ij}]_{i\in\mathcal{S}^{+},j\in\mathcal{S}^{-}}$, partitioned according to $\mathcal{S}=\mathcal{S}^+\cup\mathcal{S}^-$,
\begin{eqnarray*}
\Phi&=&
\left[
\begin{array}{cc}
0& \Psi\\
\Xi& 0
\end{array}
\right],
\end{eqnarray*}
such that, with
\begin{equation}\label{eq:theta}
\theta=\inf\{t>0:\widetilde Y(t)=0\},
\end{equation}
we have
\begin{eqnarray*}
\mu_i{\Phi}_{ij}(\mathcal{A}_v)&=&
\int_{x=0}^{\infty}d\mu_i(x)
\mathbb{P}(
\varphi( \theta)=j, X(\theta)\in \mathcal{A}_v
\ |\ \varphi(0)=i,X(0)=x,\widetilde Y(0)=0).
\end{eqnarray*}
The quantity $\mu_i{\Phi}_{ij}(\mathcal{A}_v)$ records the probability that the fluid level $\widetilde Y(\cdot)$ first returns to level $0$ and does so in phase $j$ and with $X(\cdot)\in \mathcal{A}_v$, given the process starts at time zero in phase $i$ and in level $X(0)$ distributed according to the initial distribution $\mu$. We illustrate this in Figure~\ref{fig:Phi}.

\begin{Remark}
Note that the quantity $\Psi$ is the same in both $\{(\varphi(t),X(t),\widetilde Y(t)):t\geq 0\}$ and $\{(\varphi(t),X(t), Y(t)):t\geq 0\}$, since the behaviour of $\widetilde Y(0)$ and $Y(0)$ above level $0$ is analogous. That is, for all $i\in\mathcal{S}^+$, $j\in\mathcal{S}^-$ we have
\begin{eqnarray*}
\mu_i{\Phi}_{ij}(\mathcal{A}_v)&=&
\int_{x=0}^{\infty}d\mu_i(x)
\mathbb{P}(
\varphi( \theta)=j, X(\theta)\in \mathcal{A}_v
\ |\ \varphi(0)=i,X(0)=x,Y(0)=0).
\end{eqnarray*}

\end{Remark}

Denote $\bmu{\Phi}(\mathcal{A}_v)$ such that $[\bmu{\Phi}(\mathcal{A}_v)]_j=\sum_{i\in \mathcal{S}^+\cup\mathcal{S}^-}\mu_i{\Phi}_{ij}(\mathcal{A}_v)$. To evaluate $\bmu{\Phi}(\mathcal{A}_v)$ we build on the methods developed in the previous section and
theory of the SFMs described in Section~\ref{Prelim}. We partition $\bmu$ according to $\mathcal{S}^+\cup\mathcal{S}^-$ so that $\bmu =\left[\begin{array}{cc}\bmu^{+}&\bmu^{-} \end{array} \right],$ and assume analogous partitioning for all other quantities in the expressions below.

Let $\widetilde {\bf f}_y(x)$ be the density matrix that is the inverse of the LST matrix $e^{(|{\bf R}|^{-1}({\bf T}-s \widecheck{\bf R}))y}$
described in Section~\ref{Prelim}, and let $\widehat {\bf f}_{\lambda;y}(y/2)dy$ be the inverse of the LST matrix $e^{(|{\bf R}|^{-1}({\bf T}-s \widecheck{\bf R})+\lambda |{\bf R}|^{-1}{\bf C})y}$. Further, let ${\bf M}=\int_{y=0}^{\infty}\widetilde {\bf f}_y(y/2)dy$ as defined in Section~\ref{Prelim}, and let $\widehat{\bf M}_{\lambda}=\int_{y=0}^{\infty}\widehat {\bf f}_{\lambda;y}(y/2)dy$ and ${\bf\Phi}_{\lambda}=
{\bf I}-
({\bf I}+\widehat{\bf M}_{\lambda})^{-1}
$.

Our final main result identifies $\bmu{\Phi}(\mathcal{A}_v)$.
\begin{thm}\label{th:main}
Assume that $\{(\varphi(t),Y(t)):t\geq 0\}$ and $\{(\varphi(t),X(t)):t\geq 0\}$ are not null recurrent. Suppose the original distribution $\bmu$ of $X(0)$ has a density $\nu_i(x)=p_i\lambda e^{-\lambda x}$, for some $\lambda>0$, $0\leq p_i\leq 1$, such that the boundary condition~\eqref{incond} is met. Then, for any set $\mathcal{A}_v=[0,v]$, $v>0$,
\begin{eqnarray}
\bmu{\Phi}(\mathcal{A}_v)
&=&
\left[
\begin{array}{cc}
\bmu^{+}\Psi(\mathcal{A}_v)
&
\bmu^{-}\Xi(\mathcal{A}_v)
\end{array}
\right]
\label{eqfirst}
\\
&=&
-\bmu (\bar{\mathcal{A}_v})
{\bf\Phi}_{\lambda}
+\bmu([0,\infty))
{\bf\Phi}
=
-e^{-\lambda v}\frac{\bnu(0)}{\lambda}
{\bf\Phi}_{\lambda}
+\left({\bf P}+\frac{\bnu(0)}{\lambda} \right)
{\bf\Phi}
.
\label{eqsecond}
\end{eqnarray}	
\end{thm}
\begin{Remark}
	 Note that $\lim_{v\to\infty}\bmu{\Phi}(\mathcal{A}_v)=\bmu([0,\infty))
	 {\bf\Phi}$, as expected, since $\bmu([0,\infty))
	 {\bf\Phi}$ is the probability (vector) that the process returns to level $Y(\theta)=0$ in some phase $\varphi(\theta)\in\mathcal{S}^{+}\cup\mathcal{S}^{-}$, assuming start from $Y(0)=0$ in some phase $\varphi(0)\in\mathcal{S}^{+}\cup\mathcal{S}^{-}$ according to the initial distribution $\bmu([0,\infty))$.

	 The term $\bmu (\bar{\mathcal{A}_v}){\bf\Phi}_{\lambda}=\left(
	\bmu^{(0)}([0,\infty))
	-\bmu^{(0)} (\mathcal{A}_v)
	\right)
	{\bf\Phi}_{\lambda}$ is the probability (vector) that the process returns to level $Y(\theta)=0$ in some level $X(\theta)>v$ and phase $\varphi(\theta)\in\mathcal{S}^{+}\cup\mathcal{S}^{-}$, assuming start from $Y(0)=0$ in some level $X(0)\geq 0$ and phase $\varphi(0)\in\mathcal{S}^{+}\cup\mathcal{S}^{-}$ according to the initial distribution $\bmu$.
	
	Further,
	\begin{eqnarray*}
	\lim_{v\to 0}\bmu{\Phi}(\mathcal{A}_v)&=&
	\left[ \sum_{i\in \mathcal{S}^+\cup\mathcal{S}^-}	
	\int_{x=0}^{\infty}d\mu_i(x)
	\mathbb{P}(
	\varphi( \theta)=j, X( \theta)=0
	\ |\ \varphi(0)=i,X(0)=x,\widetilde Y(0)=0)
	\right]_{j\in \mathcal{S}^+\cup\mathcal{S}^-}
	\nonumber\\
	&=&
	-\frac{\bnu(0)}{\lambda}
	{\bf\Phi}_{\lambda}
	+\left({\bf P}+\frac{\bnu(0)}{\lambda} \right)
	{\bf\Phi}
	,
	\end{eqnarray*}
	where $[\lim_{v\to 0}\bmu{\Phi}(\mathcal{A}_v)]_j=0$ for $j\in\mathcal{S}_+$ since no mass may exist at level $X(\cdot)=0$ in phases $j\in\mathcal{S}_+$.

	Matrix ${\bf M}$ can be obtained using Equation~\eqref{Meq} in Section~\ref{Prelim} and algorithms for ${\bf \Phi}$ (see e.g.~\cite{BOT5} and references there). We compute $\widehat{\bf M}_{\lambda}$ using analogous methods, with $({\bf T}+\lambda{\bf C})$ in place of ${\bf T}$, since $|{\bf R}|^{-1}({\bf T}-s \widecheck{\bf R})+\lambda |{\bf R}|^{-1}{\bf C}
	=	|{\bf R}|^{-1}(({\bf T}+\lambda{\bf C})-s \widecheck{\bf R})$. Also note that, with ${\bf R}^+=diag(r_j)_{j\in\mathcal{S}^+}$,
${\bf C}^+=diag(c_j)_{j\in\mathcal{S}^+}$, ${\bf C}^-=diag(c_j)_{j\in\mathcal{S}^-}$, we have
	\begin{eqnarray*}
	\left(|{\bf R}|^{-1}({\bf T}-s \widecheck{\bf R})+\lambda |{\bf R}|^{-1}{\bf C}\right){\bf 1}
	&=&
	|{\bf R}|^{-1}
	\left(
	{\bf T}
	-s
	\left[\begin{array}{cc} {\bf R}^{+}&{\bf 0}\\{\bf 0}& {\bf 0}\end{array} \right]
	+\lambda \left[\begin{array}{cc} {\bf C}^{+}&{\bf 0}\\{\bf 0}& {\bf C}^{-}\end{array} \right]
	\right){\bf 1}
	\\
	&=&
	-s
	\left[\begin{array}{c} {\bf 1}^{+}\\{\bf 0}^-\end{array} \right]
	+\lambda \left[\begin{array}{cc} {\bf C}^{+}{\bf 1}&\\ {\bf C}^{-}{\bf 1}\end{array} \right].
	\end{eqnarray*}
Thus $|{\bf R}|^{-1}({\bf T}-s \widecheck{\bf R})+\lambda |{\bf R}|^{-1}{\bf C}$ has strictly positive sums for all rows $i\in\mathcal{S}^-\cap \mathcal{S}_+$, and for all rows $i\in\mathcal{S}^+\cap\mathcal{S}_+$ with $ \lambda c_i > s$.
\end{Remark}
{\bf Proof:}\quad First, note that $$[\bmu {\Phi}(\mathcal{A}_v)]_j=\sum_{i\in \mathcal{S}^+\cup\mathcal{S}^-}\mu_i{\Phi}_{ij}(\mathcal{A}_v)=(\bmu^{+}\Psi(\mathcal{A}_v))I(j\in\mathcal{S}^{+})
+(\bmu^{-}\Xi(\mathcal{A}_v))I(j\in\mathcal{S}^{-})$$
and thus~\eqref{eqfirst} follows.
Next, for any $s>0$, $y>0$, by argument similar to the proof of Theorem~\ref{exp_distr},
\begin{eqnarray*}
\left[\bmu e^{Dy}(\mathcal{A}_v,s)\right]_j&=&
\left[-\bmu (\bar{\mathcal{A}_v})
e^{(|{\bf R}|^{-1}({\bf T}-s \widecheck{\bf R})+\lambda |{\bf R}|^{-1}{\bf C})y}
+\bmu([0,\infty))
e^{(|{\bf R}|^{-1}({\bf T}-s \widecheck{\bf R}))y}\right]_j
\end{eqnarray*}
is the Laplace-Stieltjes transform of the distribution of the total upward shift in $Y(\cdot)$ accumulated by the time the in-out fluid of the process $Y(\cdot)$ first reaches level $y$ and does so in phase $j$ and with $X(\cdot)\in \mathcal{A}_v$, given that the process starts at time zero according to distribution $\bmu$. That is,
\begin{equation*}
\left[\bmu e^{D y}(\mathcal{A}_v,s)\right]_j=\mathbb{E}
(e^{-sh_+(\omega(y))}  I(\varphi(\omega(y))=j, X(\omega(y)\in\mathcal{A}_v))).
\end{equation*}
For $0\leq x\leq y$, we define the density $\bmu f_y(x)$ as the inverse of the LST $\bmu e^{Dy}(\mathcal{A}_v,s)$, so that
\begin{equation*}
\bmu e^{Dy}(\mathcal{A}_v,s)
=\int_{x=0}^{\infty}e^{-sx}\bmu f_y(\mathcal{A}_v,x)dx,
\end{equation*}
where $[\bmu f_y(x)]_j$ is the probability density that the total upward shift in $Y(\cdot)$ at time $\omega(y)$ is $x$ and the phase is $\varphi(\omega(y))=j$  and with $X(\cdot)\in \mathcal{A}_v$, given that the process starts at time zero according to distribution $\bmu$, where $f_y(x)$ is the corresponding operator. We have
\begin{eqnarray*}
\bmu f_y(\mathcal{A}_v,x)&=&
-\bmu (\bar{\mathcal{A}_v})
\widehat {\bf f}_{\lambda;y}(x)
+\bmu([0,\infty))
\widetilde {\bf f}_y(x).
\end{eqnarray*}
Consider
\begin{eqnarray*}
\bmu M (\mathcal{A}_v)
&=&
-\bmu (\bar{\mathcal{A}_v})
\widehat{\bf M}_{\lambda}
+\bmu([0,\infty))
{\bf M}
=
-\left(
\bmu^{(0)}([0,\infty))
-\bmu^{(0)} (\mathcal{A}_v)
\right)
\widehat{\bf M}_{\lambda}
+\bmu([0,\infty))
{\bf M},
\end{eqnarray*}
and note that $[\bmu M (\mathcal{A}_v)]_j$ is the expected number of visits to level $0$ and doing so in phase $j$ and with $X(\cdot)\in \mathcal{A}_v$, given that the process starts at time zero according to distribution $\bmu$, and $M$ is the corresponding operator. We have
\begin{eqnarray*}
\bmu M (\mathcal{A}_v)&=&
\int_{x=0}^v {}_{M}\bnu (x)dx
+
{}_{M}{\bf P},
\end{eqnarray*}
where,
\begin{eqnarray*}
{}_{M}\bnu (x)&=&\bnu (x) \widehat{\bf M}_{\lambda},\\
{}_{M}{\bf P}&=&
-\bmu^{(0)}([0,\infty))\widehat{\bf M}_{\lambda}
+\bmu([0,\infty))
{\bf M}
=
\frac{\bnu(0)}{\lambda}
\left(
{\bf M}
-
\widehat{\bf M}_{\lambda}
\right)
+ {\bf P}{\bf M},
\end{eqnarray*}
and by induction, $\bmu M^n $ has density ${}_{M^n}\bnu (x)$ and mass ${}_{M^n}{\bf P}$ given by,
\begin{eqnarray}
{}_{M^n}\bnu (x)&=& \bnu (x)
\left(
\widehat{\bf M}_{\lambda}
\right)^n
\label{densM}
,\\
{}_{M^n}{\bf P}&=&
\frac{\bnu (0)}{\lambda}
\sum_{k=0}^{n-1}
\left( \widehat{\bf M}_{\lambda} \right)^k
\left(
{\bf M}
-\widehat{\bf M}_{\lambda}
\right)
\left(
{\bf M}
\right)^{n-1-k}
+ {\bf P}
\left(
{\bf M}
\right)^n.
\label{massM}
\end{eqnarray}
Indeed, by above, the expressions~\eqref{densM}-\eqref{massM} for $\bmu M^n$ are clearly true for $n=1$. Suppose that they hold for some $n\geq 1$. Then,
\begin{eqnarray*}
{}_{M^{n+1}}\bnu (x)&=& {}_{M^n}\bnu (x)
\widehat{\bf M}_{\lambda}
=
\bnu (x)
\left(
\widehat{\bf M}_{\lambda}
\right)^{n+1},
\end{eqnarray*}
and
\begin{eqnarray*}
{}_{M^{n+1}}{\bf P}&=&
\frac{{}_{M^n}\bnu(0)}{\lambda}
\left(
{\bf M}
-
\widehat{\bf M}_{\lambda}
\right)
+ {}_{M^n}{\bf P}{\bf M}
\nonumber\\
&=&
\frac{\bnu(0)}{\lambda}
\left(
\widehat{\bf M}_{\lambda}
\right)^n
\left(
{\bf M}
-
\widehat{\bf M}_{\lambda}
\right)
\nonumber\\
&&
+
\left(
\frac{\bnu (0)}{\lambda}
\sum_{k=0}^{n-1}
\left( \widehat{\bf M}_{\lambda} \right)^k
\left(
{\bf M}
-\widehat{\bf M}_{\lambda}
\right)
\left(
{\bf M}
\right)^{n-1-k}
+ {\bf P}
\left(
{\bf M}
\right)^n
\right){\bf M}
\nonumber\\
&=&
\frac{\bnu (0)}{\lambda}
\sum_{k=0}^{(n+1)-1}
\left( \widehat{\bf M}_{\lambda} \right)^k
\left(
{\bf M}
-\widehat{\bf M}_{\lambda}
\right)
\left(
{\bf M}
\right)^{(n+1)-1-k}
+ {\bf P}
\left(
{\bf M}
\right)^{n+1},
\end{eqnarray*}
which completes the proof of~\eqref{densM}-\eqref{massM} for all $n\geq 1$.

Now, by the physical interpretations similar to those used for a one-dimensional stochastic fluid model summarised in Section~\ref{Prelim}, we consider a Markov process observed at the times the SFFM $\{(\varphi(t),X(t),\widetilde Y(t)):t\geq 0\}$ visits level $\widetilde Y(\cdot)=0$, with one-step transition operator $\Phi$, and note that by standard theory of Markov chains we have $M=\sum_{n=1}^{\infty}(\Phi)^n$, and so
\begin{eqnarray*}
M&=&\Phi +M\Phi=\Phi(I-\Phi)^{-1}=\Phi\sum_{n=0}^{\infty}(\Phi)^n,\\
\Phi&=&(I+M)^{-1}M=\sum_{n=0}^{\infty}(-M)^nM =-\sum_{n=1}^{\infty}(-1)^n M^n ,
\end{eqnarray*}
and thus
\begin{eqnarray*}
\bmu{\Phi}(\mathcal{A}_v)&=&
-\sum_{n=1}^{\infty} (-1)^n \bmu M^n (\mathcal{A}_v)
=
-\sum_{n=1}^{\infty} (-1)^n
\left(
{}_{M^n}\bmu^{(0)} (\mathcal{A}_v)
+{}_{M^n}{\bf P}
\right).
\end{eqnarray*}
Consequently an explicit expression of $\bmu{\Phi}(\mathcal{A}_v)$ follows by~\eqref{densM}-\eqref{massM}. Indeed, we have
\begin{eqnarray*}
\bmu{\Phi}(\mathcal{A}_v)
&=&
-\bmu^{(0)} (\mathcal{A}_v)\sum_{n=1}^{\infty}
\left(
-
\widehat{\bf M}_{\lambda}
\right)^n
-{\bf P}\sum_{n=1}^{\infty}
\left(-
{\bf M}
\right)^n
\nonumber\\
&&
+\frac{\bnu (0)}{\lambda}\sum_{n=1}^{\infty}
\sum_{k=0}^{n-1}
\left( - \widehat{\bf M}_{\lambda} \right)^k
\left(
{\bf M}
-\widehat{\bf M}_{\lambda}
\right)
\left(
-
{\bf M}
\right)^{n-1-k}
\nonumber\\
&=&
\bmu^{(0)} (\mathcal{A}_v)
\left(
{\bf I}-
({\bf I}+\widehat{\bf M}_{\lambda})^{-1}
\right)
+{\bf P}
\left(
{\bf I}-
({\bf I}+{\bf M})^{-1}
\right)
\nonumber\\
&&
+\frac{\bnu (0)}{\lambda}
\left(
{\bf I}
+
\widehat{\bf M}_{\lambda}
\right)^{-1}
\left(
-{\bf I}
+{\bf I}
+
{\bf M}
-\widehat{\bf M}_{\lambda}
\right)
\left(
{\bf I}
+
{\bf M}
\right)^{-1}
\nonumber\\
&=&
\bmu^{(0)} (\mathcal{A}_v)
{\bf\Phi}_{\lambda}
+{\bf P}
{\bf\Phi}
+\bmu^{(0)}([0,\infty))
\left(
{\bf I}
+
\widehat{\bf M}_{\lambda}
\right)^{-1}
\left(
-{\bf I}
-\widehat{\bf M}_{\lambda}
\right)
\left(
{\bf I}
+
{\bf M}
\right)^{-1}
\nonumber\\
&&
+\bmu^{(0)}([0,\infty))
\left(
{\bf I}
+
\widehat{\bf M}_{\lambda}
\right)^{-1}
\left(
{\bf I}
+
{\bf M}
\right)
\left(
{\bf I}
+
{\bf M}
\right)^{-1}
\nonumber\\
&=&
\bmu^{(0)} (\mathcal{A}_v)
{\bf\Phi}_{\lambda}
+{\bf P}
{\bf\Phi}
-\bmu^{(0)}([0,\infty))
\left(
{\bf I}
+
{\bf M}
\right)^{-1}
+\bmu^{(0)}([0,\infty))
\left(
{\bf I}
+
\widehat{\bf M}_{\lambda}
\right)^{-1}
\nonumber\\
&=&
\bmu^{(0)} (\mathcal{A}_v)
{\bf\Phi}_{\lambda}
+{\bf P}
{\bf\Phi}
+
\bmu^{(0)}([0,\infty))
\left(
{\bf I}-
\left(
{\bf I}
+
{\bf M}
\right)^{-1}
\right)
-\bmu^{(0)}([0,\infty))
\nonumber\\
&&-\bmu^{(0)}([0,\infty))
\left(
{\bf I}-
\left(
{\bf I}
+
\widehat{\bf M}_{\lambda}
\right)^{-1}
\right)
+\bmu^{(0)}([0,\infty))
\nonumber\\
&=&
\bmu^{(0)} (\mathcal{A}_v)
{\bf\Phi}_{\lambda}
+{\bf P}
{\bf\Phi}
+
\bmu^{(0)}([0,\infty)){\bf\Phi}
-\bmu^{(0)}([0,\infty)){\bf\Phi}_{\lambda}
,
\end{eqnarray*}
which gives~\eqref{eqsecond}.
\hfill\rule{9pt}{9pt}

\subsection{Numerical examples}\label{sec:numeex2}

The physical intepretations of the key quantities in the examples below are:
\begin{itemize}
	\item $[{\bf\Psi}]_{ij}$ is the probability that given start from level $Y(0)=0$ and phase $\varphi(0)=i$, the fluid $Y(\cdot)$ will return to level $0$ at time $\theta$ and does so on phase $\varphi(\theta)=j$;
	\item $[{\bf\Xi}]_{ij}$ is the probability that given start from level $\widetilde Y(0)=0$ and phase $\varphi(0)=i$, the fluid $\widetilde Y(\cdot)$ will return to level $0$ at time $\theta$ and does so on phase $\varphi(\theta)=j$;
	\item $[\bmu([0,\infty))
	]_j$ is the probability that $\varphi(0)=j$;
	\item $[\bmu([0,\infty))
	{\bf\Phi}]_j$ is the probability that $\varphi(\theta)=j$;
	\item $[e^{-\lambda v}\frac{\bnu(0)}{\lambda}
	{\bf\Phi}_{\lambda}]_j$ is the probability that $\varphi(\theta)=j$ and $X(\theta)\notin \mathcal{A}_v$;
	\item $[\bmu{\Phi}(\mathcal{A}_v)]_j$ is the probability that $\varphi(\theta)=j$ and $X(\theta)\in \mathcal{A}_v$;
	\item $[\lim_{y\to\infty}\bmu e^{Dy}(\mathcal{A}_v)]_j$ is the stationary probability that in the long run $\varphi=j$ and $X\in \mathcal{A}_v$;
	\item $[\bmu e^{Dy}(\mathcal{A}_v)]_j$ is the probability that $\varphi(\omega(y))=j$ and $X(\omega(y))\in \mathcal{A}_v$;
	\item $\lim_{v\to 0}[\bmu e^{Dy}(\mathcal{A}_v)]_j$ is the probability that $\varphi(\omega(y))=j$ and $X(\omega(y))=0$;
	\item $[\lim_{v\to\infty}\bmu e^{Dy}(\mathcal{A}_v)
	-\lim_{v\to 0}\bmu e^{Dy}(\mathcal{A}_v)]_j$ is the probability that $\varphi(\omega(y))=j$ and $X(\omega(y))>0$;
\end{itemize}
given the initial distribution $\bmu$, with $\omega(y)$ defined in~\eqref{eq:omegay}, $\theta$ defined in~\eqref{eq:theta}, and $\theta\equiv \theta_0$ defined in~\eqref{eq:thetaz} whenever $\varphi(0)\in\mathcal{S}^+$.

\begin{exa}\label{sec3_exa}
\rm Consider Example~\ref{ex2}, where $\mathcal{S}=\{1,2\}$, $\mathcal{S}_+=\mathcal{S}^+=\{1\}$, $\mathcal{S}_-=\mathcal{S}^-=\{2\}$, $|r_i|=|c_i|=1$ for all $i$. Let $p=0.2$, $b=1$, $d=1$. We have
\begin{eqnarray*}
&&|{\bf R}|^{-1}{\bf T}=|{\bf C}|^{-1}{\bf T}
=
\left[
\begin{array}{cc}
-2&2\\1&-1
\end{array}
\right]
,\
|{\bf R}|^{-1}{\bf T}+\lambda |{\bf R}|^{-1}{\bf C}
=
\left[
\begin{array}{cc}
-1&2\\
1&-2
\end{array}
\right]
,
\\
&&{\bf P}=
\left[
\begin{array}{cc}
0&0.2
\end{array}
\right]
,\
\bnu(0)=
\left[
\begin{array}{cc}
0.2&0.6
\end{array}
\right],
\end{eqnarray*}
which gives, by Theorem~\ref{th:main},
\begin{eqnarray*}
&&{\bf\Phi}=
\left[
\begin{array}{cc}
0&{\bf\Psi}\\
{\bf\Xi}&0
\end{array}
\right]
=
\left[
\begin{array}{cc}
0&1\\
0.5&0
\end{array}
\right],
{\bf\Phi}_{\lambda}=
\left[
\begin{array}{cc}
0&{\bf\Psi}_{\lambda}\\
{\bf\Xi}_{\lambda}&0
\end{array}
\right]
=
\left[
\begin{array}{cc}
0&1\\
0.5&0
\end{array}
\right]
,
\\
&&\bmu([0,\infty))
=
\left[
\begin{array}{cc}
0.2&0.8
\end{array}
\right]
,
\
\bmu([0,\infty))
{\bf\Phi}
=
\left[
\begin{array}{cc}
0.4&0.2
\end{array}
\right],
\\
&&
\frac{\bnu(0)}{\lambda}=
\left[
\begin{array}{cc}
0.2&0.6
\end{array}
\right]
,
\
\frac{\bnu(0)}{\lambda}
{\bf\Phi}_{\lambda}
=\left[
\begin{array}{cc}
0.3&0.2
\end{array}
\right],
\\
&&\bmu{\Phi}(\mathcal{A}_v)
=
\left[
\begin{array}{cc}
0.4&0.2
\end{array}
\right]
-
e^{-\lambda v}\left[
\begin{array}{cc}
0.3&0.2
\end{array}
\right].
\end{eqnarray*}
Moreover, by Theorem~\ref{exp_distr},
\begin{eqnarray*}
\lim_{y\to\infty}\bmu e^{Dy}(\mathcal{A}_v)&=&
\lim_{y\to\infty}
\left(
-e^{-\lambda v}\frac{\bnu(0)}{\lambda}
e^{(|{\bf R}|^{-1}{\bf T}+\lambda |{\bf R}|^{-1}{\bf C})y}
+\bmu ([ 0,\infty))
e^{(|{\bf R}|^{-1}{\bf T})y}
\right)
,
\nonumber\\
&=&
-e^{- v}
\left[
\begin{array}{cc}
0.3333  &  0.3333
\end{array}
\right]
+
\left[
\begin{array}{cc}
0.3333  &  0.6667
\end{array}
\right].
\end{eqnarray*}
We plot the values of $\bmu e^{Dy}(\mathcal{A}_v)$, $\lim_{v\to 0}\bmu e^{Dy}(\mathcal{A}_v)$, $\lim_{v\to\infty}\bmu e^{Dy}(\mathcal{A}_v)
-\lim_{v\to 0}\bmu e^{Dy}(\mathcal{A}_v)$ and $\bmu{\Phi}(\mathcal{A}_v)$ in Figure~\ref{fig:muDyAlimits}.
\begin{figure}
	\centering
	\begin{tabular}{@{}cc@{}}
		\includegraphics[width=.4\textwidth]{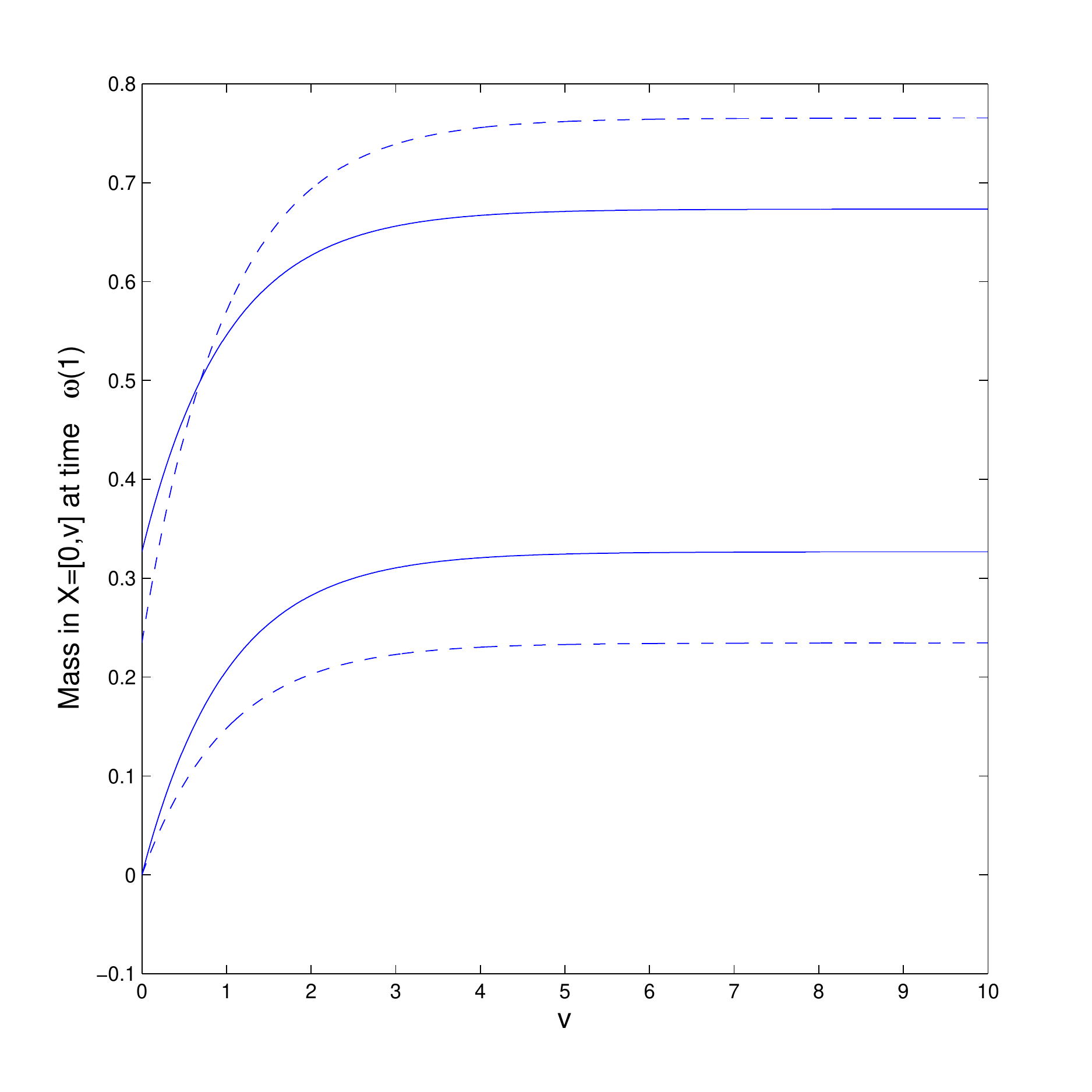} &
		\includegraphics[width=.4\textwidth]{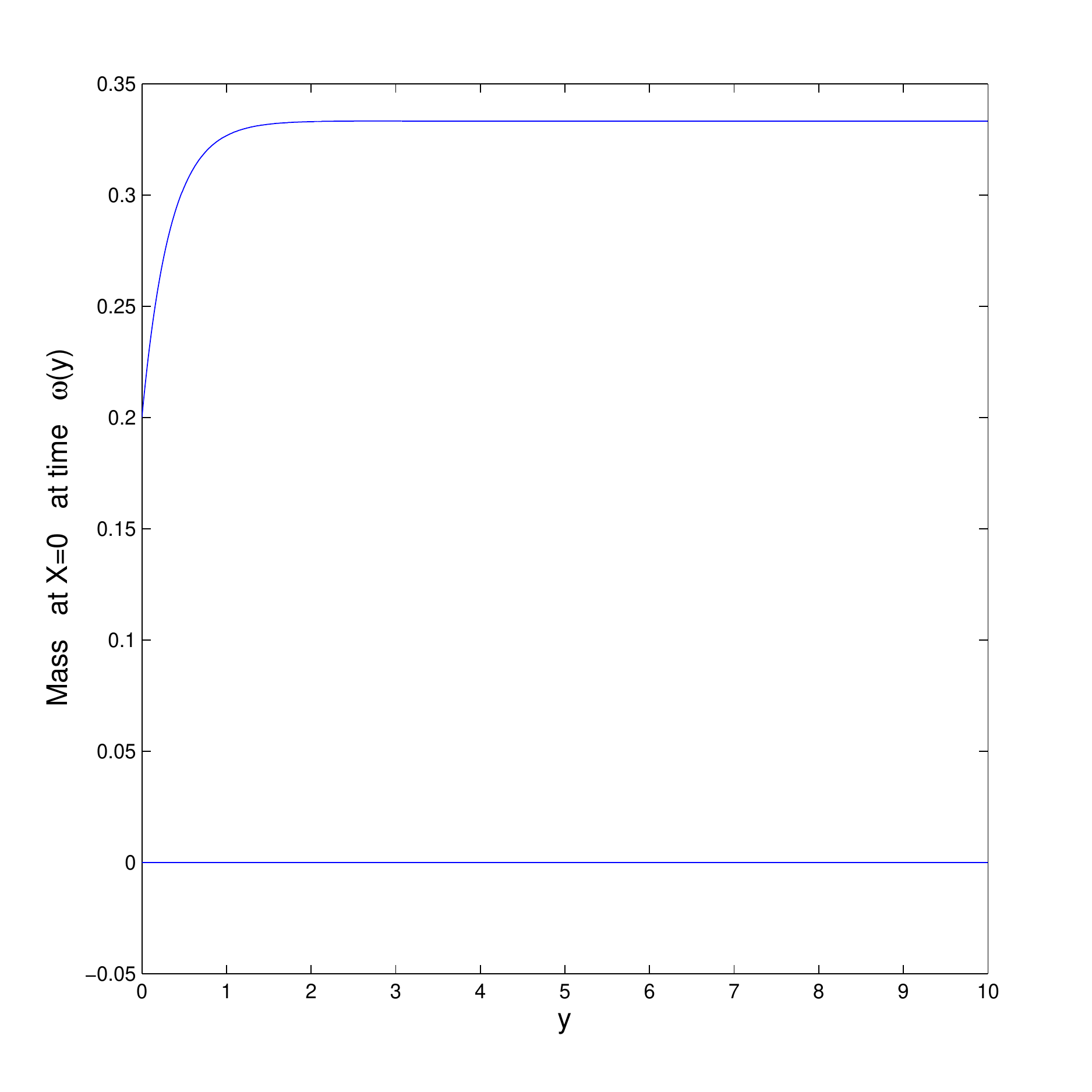} \\
		\footnotesize (a) $[\bmu e^{Dy}(\mathcal{A}_v)]_j$ &
		\footnotesize (b) $\lim_{v\to 0}[\bmu e^{Dy}(\mathcal{A}_v)]_j$\\
	\end{tabular}
	\begin{tabular}{@{}cc@{}}
		\includegraphics[width=.4\textwidth]{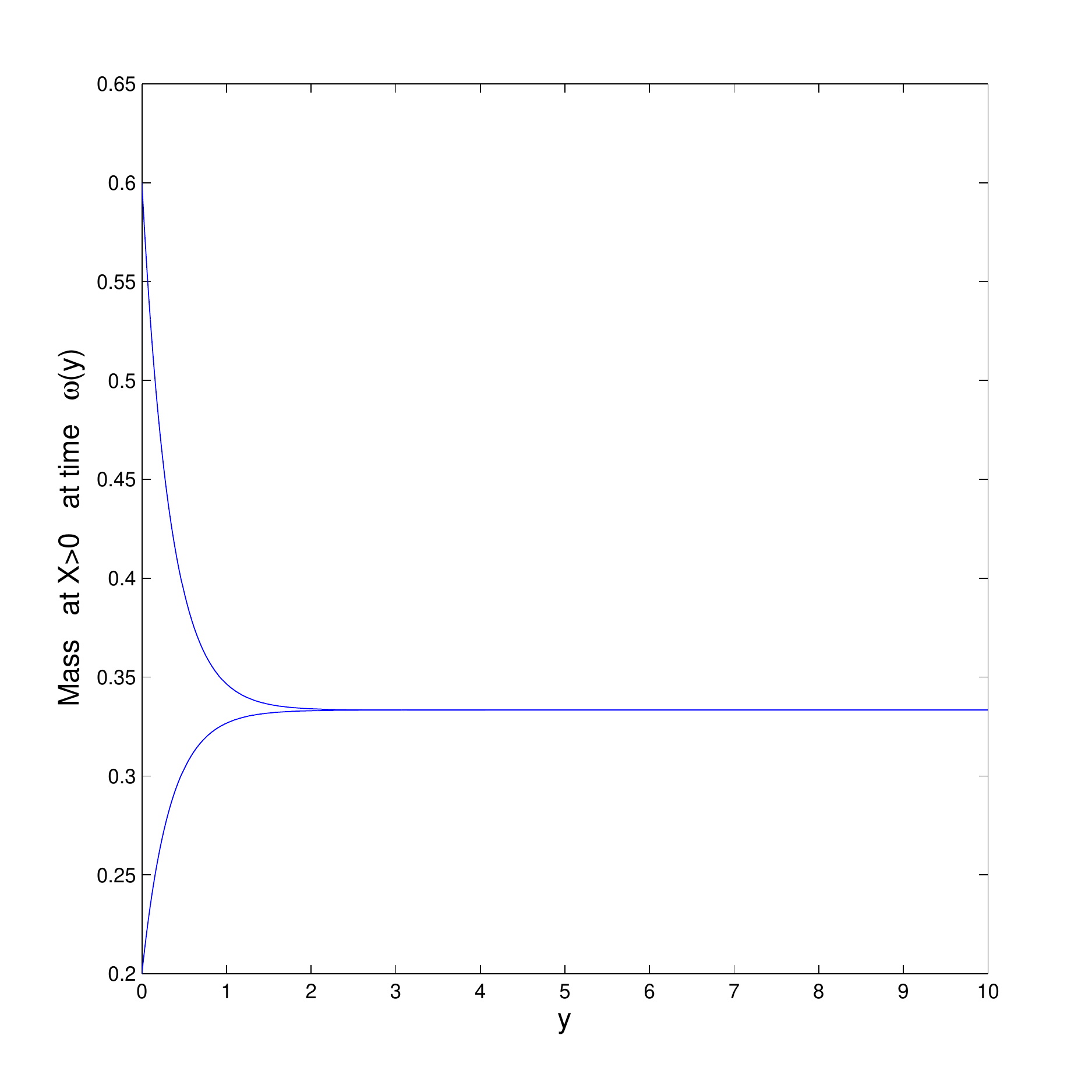} &
		\includegraphics[width=.4\textwidth]{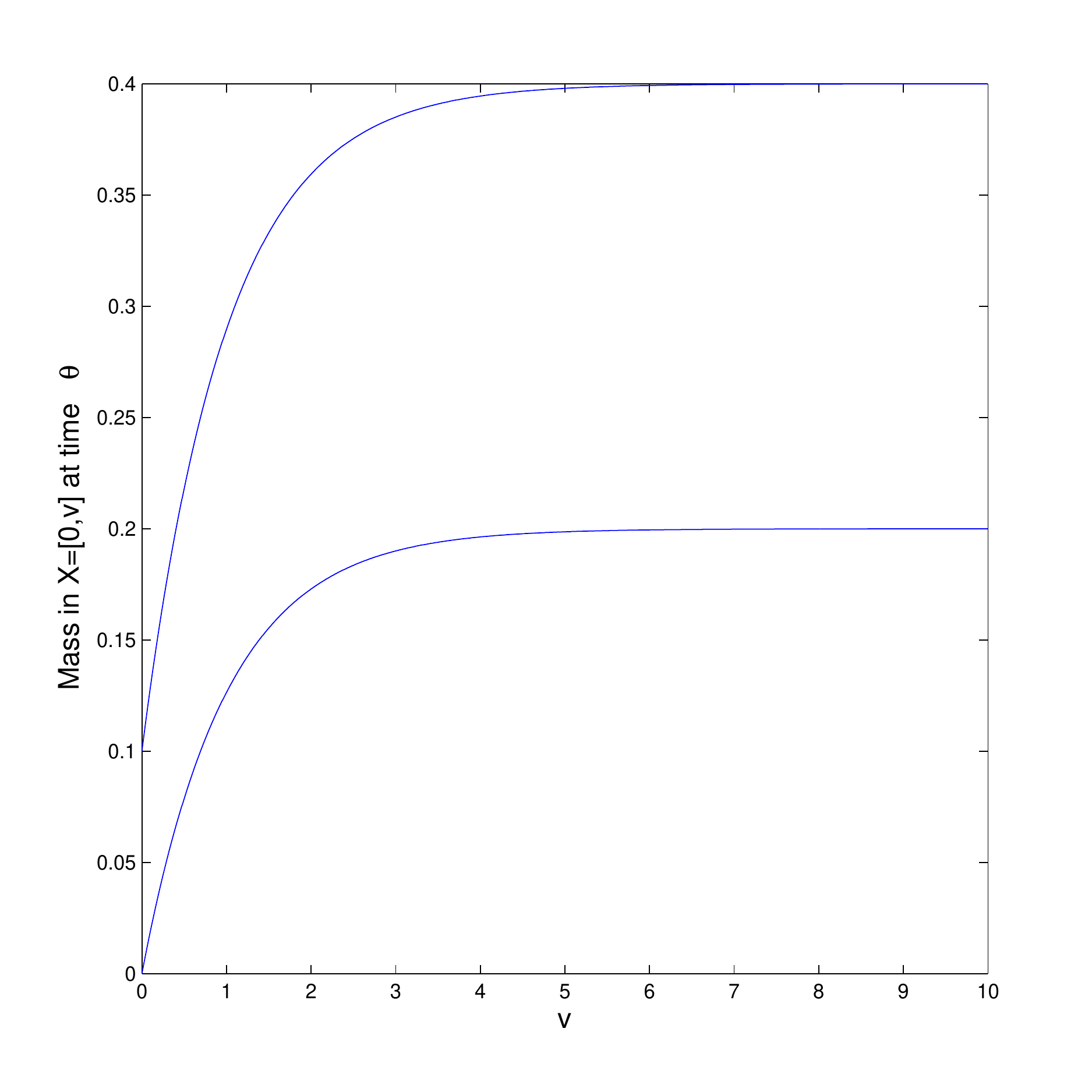} \\
			\footnotesize (c) $[\lim_{v\to\infty}\bmu e^{Dy}(\mathcal{A}_v)
			-\lim_{v\to 0}\bmu e^{Dy}(\mathcal{A}_v)]_j$ &
			\footnotesize (d) $[\bmu{\Phi}(\mathcal{A}_v)]_j$\\
	\end{tabular}
	\caption{\footnotesize Example~\ref{sec3_exa}: (a)\quad $[\bmu e^{Dy}(\mathcal{A}_v)]_j$ for $y=0.1\mbox{ (dashed line), and }y=1$;\quad (b)\quad $\lim_{v\to 0}[\bmu e^{Dy}(\mathcal{A}_v)]_j$ (mass at $X=0$ at time $\omega(y)$ is zero for $j=1\in\mathcal{S}_+$); \quad (c)\quad  $[\lim_{v\to\infty}\bmu e^{Dy}(\mathcal{A}_v)
		-\lim_{v\to 0}\bmu e^{Dy}(\mathcal{A}_v)]_j$ (mass at $X>0$ at time $\omega(y)$);\quad (d)\quad $[\bmu{\Phi}(\mathcal{A}_v)]_j=[\bmu{\Psi}(\mathcal{A}_v)]_j$ for $j\in\mathcal{S}^-$ and $[\bmu{\Phi}(\mathcal{A}_v)]_j=[\bmu{\Xi}(\mathcal{A}_v)]_j$ for $j\in\mathcal{S}^+$.}
	\label{fig:muDyAlimits}
\end{figure}

\end{exa}

\begin{exa}\label{sec3_exb} \rm Consider now Example~\ref{ex5}(a), where $\mathcal{S}=\{1,2,3,4\}$, $\mathcal{S}_+=\{1,2\}$, $\mathcal{S}_-=\{3,4\}$, $\mathcal{S}^+=\{1,4\}$, $\mathcal{S}^-=\{2,3\}$, $|r_i|=|c_i|=1$ for all $i$. Let $b=\beta=1$. Then
	\begin{eqnarray*}
	&&|{\bf R}|^{-1}
	=
	\left[
	\begin{array}{cc}
	-2{\bf I}&{\bf E}\\
	\widetilde{\bf E}&-{\bf I}
	\end{array}
	\right]
=
\left[
\begin{array}{cccc}
-2&0&1&1\\
0&-2&1&1\\
1-r&r&-1&0\\
1-r&r&0&-1
\end{array}
\right]
	={\bf T}
	=|{\bf C}|^{-1}{\bf T}
	,\nonumber\\
	&&{\bf P}_{-}=[p/2\ \ p/2],\
	\bnu_+(0)=e^{- x}p[1-r\ \ r],
	\
	\bnu_{-}(x)=e^{- x}(1/2 -p)[1\ \ 1].
	\end{eqnarray*}
We partition the quantities below according to $\mathcal{S}^+\cup\mathcal{S}^-=\{1,4\}\cup\{2,3\}$ so that they have a suitable form for application of Theorem~\ref{th:main}. We have,
\begin{eqnarray*}
&&|{\bf R}|^{-1}{\bf T}
	=
	\left[
	\begin{array}{cccc}
	-2&1&0&1\\
	1-r&-1&r&0\\
	0&1&-2&1\\
	1-r&0&r&-1
	\end{array}
	\right]
	,\
|{\bf R}|^{-1}{\bf T}+\beta|{\bf R}|^{-1}{\bf C}=	
	\left[
	\begin{array}{cccc}
	-1&1&0&1\\
	1-r&-2&r&0\\
	0&1&-1&1\\
	1-r&0&r&-2
	\end{array}
	\right]
	,
\nonumber\\
&&
{\bf P}=
[0\ \ p/2\ \ 0\ \ p/2],\
\bnu(0)=[p(1-r)\ \ (1/2-p)\ \ pr\ \ (1/2-p)].
\end{eqnarray*}
Let $r=0.6$, $p=0.2$. By Theorem~\ref{th:main},
\begin{eqnarray*}
&&{\bf\Phi}=
\left[
\begin{array}{cc}
0&{\bf\Psi}\\
{\bf\Xi}&0
\end{array}
\right]
=
\left[
\begin{array}{cccc}
         0&         0&    0.2662&    0.7338\\
         0&         0&    0.4314&    0.5686\\
         0.1774&    0.7190&         0&         0\\
         0.2935&    0.5686&         0&         0
\end{array}
\right],
\nonumber\\
&&{\bf\Phi}_{\lambda}=
\left[
\begin{array}{cccc}
        0&         0&    0.6354&    0.7292\\
        0&         0&    0.3962&    0.2077\\
        0.4236&    0.6603&         0&         0\\
        0.2917&    0.2077&         0&         0
\end{array}
\right]
,
\nonumber\\
&&\bmu([0,\infty))
=
\left[
\begin{array}{cccc}
0.08&    0.40&   0.12&    0.40
\end{array}
\right]
,
\
\bmu([0,\infty))
{\bf\Phi}
=
\left[
\begin{array}{cccc}
0.1387&    0.3137&    0.1939&    0.2861
\end{array}
\right],
\nonumber\\
&&
\frac{\bnu(0)}{\lambda}=
\left[
\begin{array}{cccc}
0.08&    0.30&    0.12&    0.30
\end{array}
\right]
,
\
\frac{\bnu(0)}{\lambda}
{\bf\Phi}_{\lambda}
=\left[
\begin{array}{cccc}
0.1383&    0.1415&    0.1697&    0.1206
\end{array}
\right],
\nonumber\\
&&\bmu{\Phi}(\mathcal{A}_v)
=
\left[
\begin{array}{cccc}
0.1387&    0.3137&    0.1939&    0.2861
\end{array}
\right]
-
e^{-\lambda v}
\left[
\begin{array}{cccc}
0.1383&    0.1415&    0.1697&    0.1206
\end{array}
\right].
\nonumber\\
\end{eqnarray*}
and by Theorem~\ref{exp_distr},
\begin{eqnarray*}
\lim_{y\to\infty}\bmu e^{Dy}(\mathcal{A}_v)&=&
\lim_{y\to\infty}
\left(
-e^{-\lambda v}\frac{\bnu(0)}{\lambda}
e^{({\bf R}{\bf T}+\lambda {\bf R}{\bf C})y}
+\bmu ([ 0,\infty))
e^{({\bf R}{\bf T})y}
\right)
,
\nonumber\\
&=&
-e^{- v}
\left[
\begin{array}{cccc}
0.1333&    0.1667&    0.2000&    0.1667
\end{array}
\right]
+
\left[
\begin{array}{cccc}
0.1333&    0.3333&    0.2000&    0.3333
\end{array}
\right].
\nonumber\\
\end{eqnarray*}
We plot the values of $\bmu e^{Dy}(\mathcal{A}_v)$, $\lim_{v\to 0}\bmu e^{Dy}(\mathcal{A}_v)$, $\lim_{v\to\infty}\bmu e^{Dy}(\mathcal{A}_v)
-\lim_{v\to 0}\bmu e^{Dy}(\mathcal{A}_v)$ and $\bmu{\Phi}(\mathcal{A}_v)$ in Figure~\ref{fig:muDyAlimits2}.
\begin{figure}
	\centering
	\begin{tabular}{@{}cc@{}}
		\includegraphics[width=.4\textwidth]{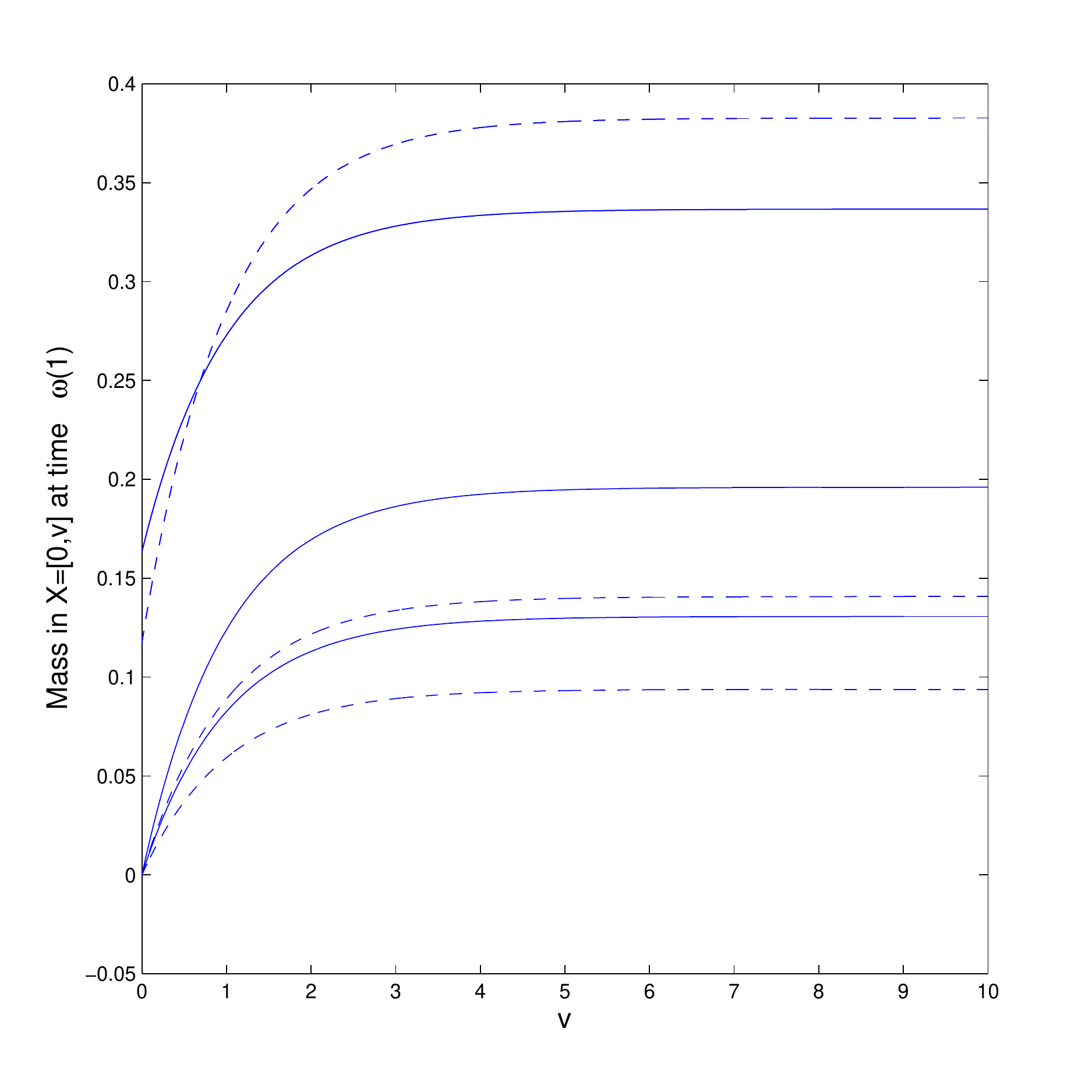} &
		\includegraphics[width=.4\textwidth]{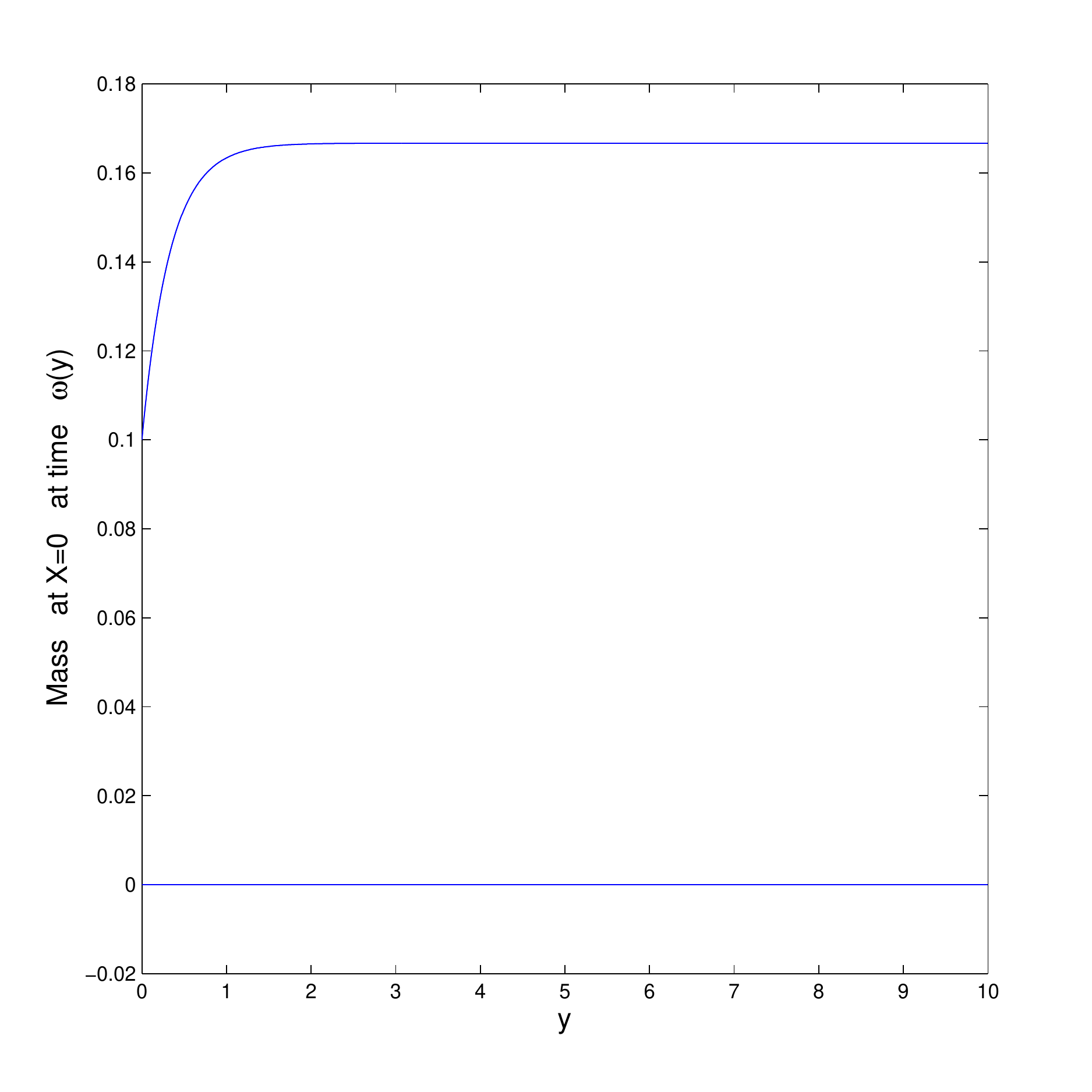} \\
		\footnotesize (a) $[\bmu e^{Dy}(\mathcal{A}_v)]_j$ &
		\footnotesize (b) $\lim_{v\to 0}[\bmu e^{Dy}(\mathcal{A}_v)]_j$\\
	\end{tabular}
	\begin{tabular}{@{}cc@{}}
		\includegraphics[width=.4\textwidth]{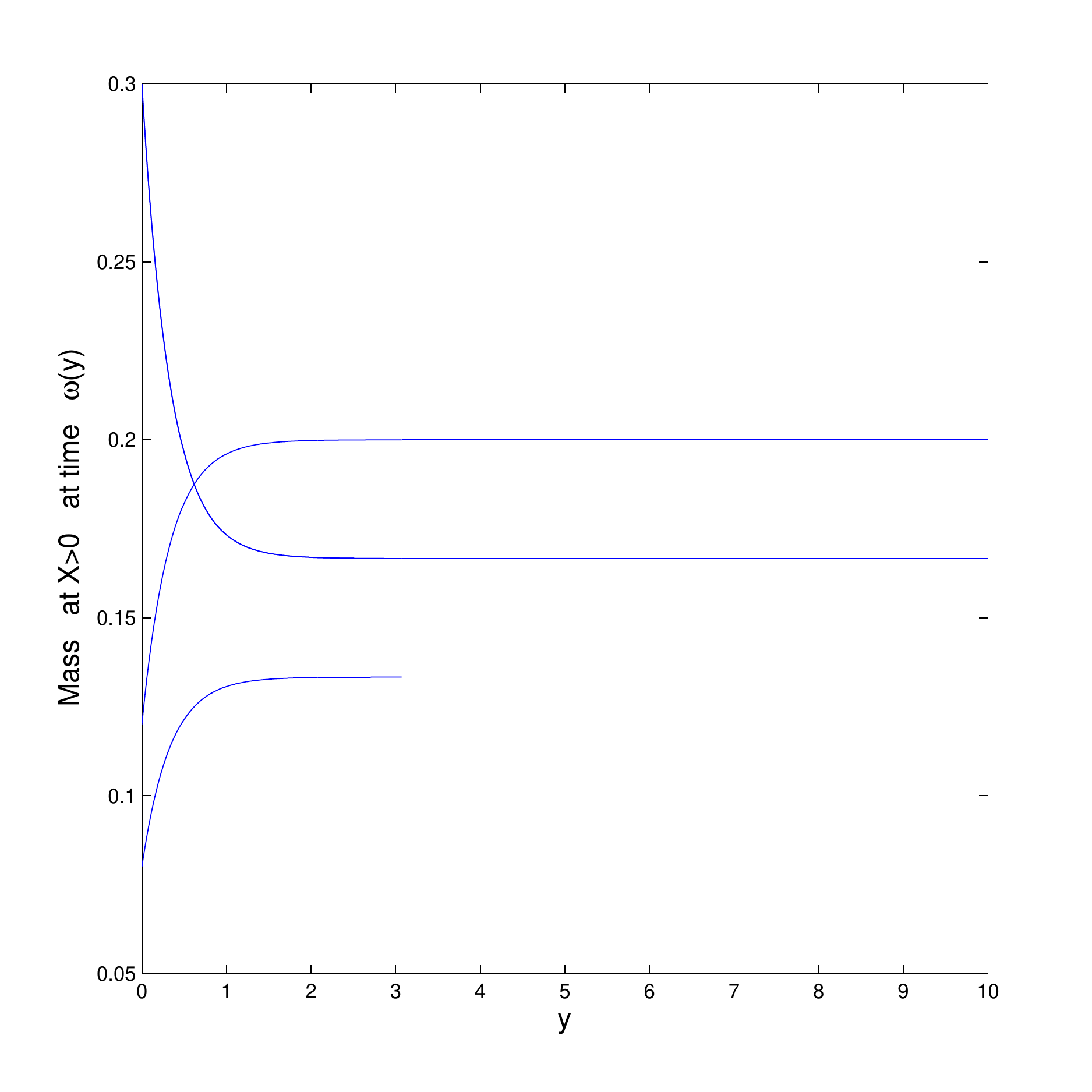} &
		\includegraphics[width=.4\textwidth]{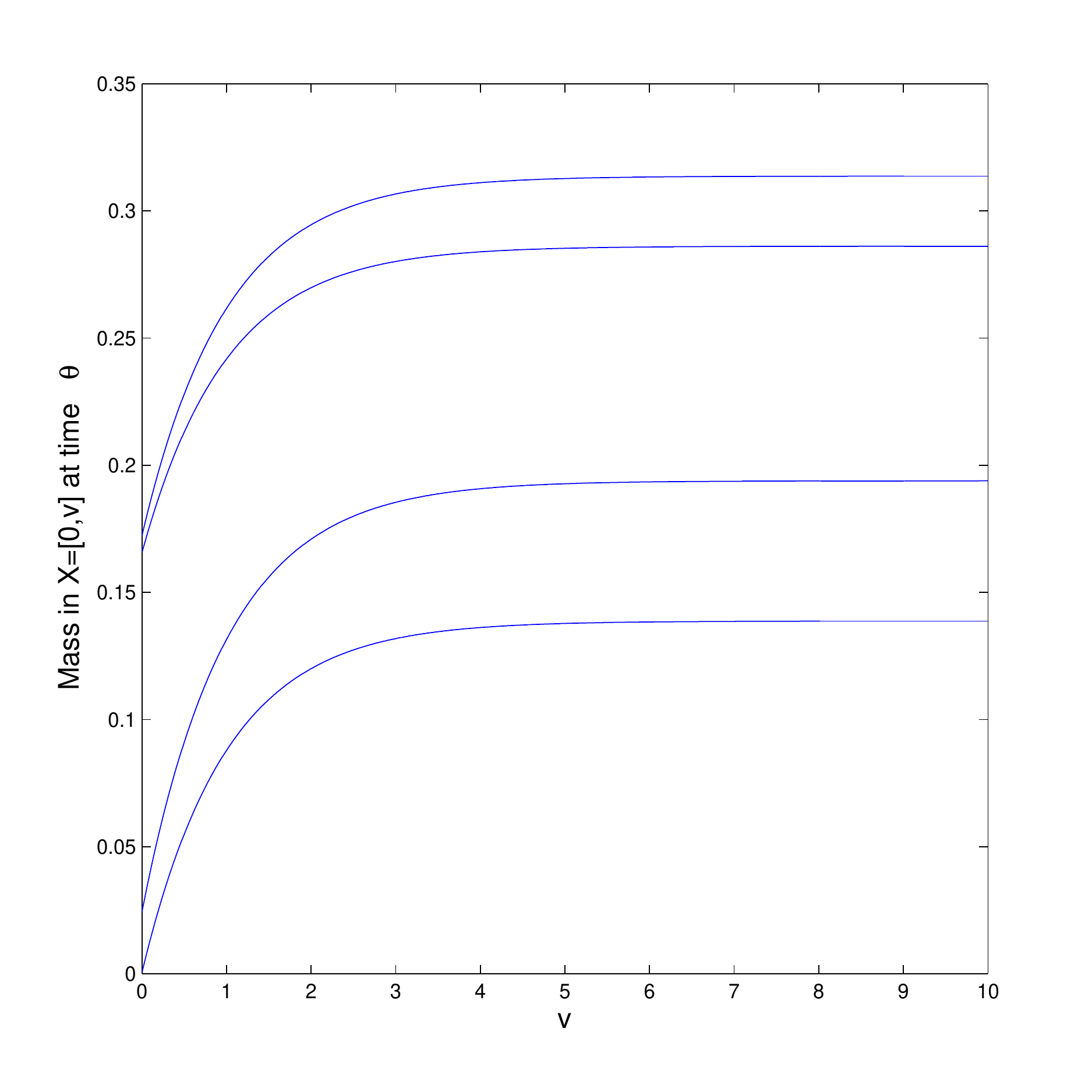} \\
		\footnotesize (c) $[\lim_{v\to\infty}\bmu e^{Dy}(\mathcal{A}_v)
		-\lim_{v\to 0}\bmu e^{Dy}(\mathcal{A}_v)]_j$ &
		\footnotesize (d) $[\bmu{\Phi}(\mathcal{A}_v)]_j$\\
	\end{tabular}
	\caption{\footnotesize Example~\ref{sec3_exb}: (a)\quad $[\bmu e^{Dy}(\mathcal{A}_v)]_j$ for $y=0.1\mbox{ (dashed line), and }y=1$;\quad (b)\quad $\lim_{v\to 0}[\bmu e^{Dy}(\mathcal{A}_v)]_j$ (mass at $X=0$ at time $\omega(y)$ is zero for $j=1\in\mathcal{S}_+$); \quad (c)\quad  $[\lim_{v\to\infty}\bmu e^{Dy}(\mathcal{A}_v)
		-\lim_{v\to 0}\bmu e^{Dy}(\mathcal{A}_v)]_j$ (mass at $X>0$ at time $\omega(y)$);\quad (d)\quad $[\bmu{\Phi}(\mathcal{A}_v)]_j=[\bmu{\Psi}(\mathcal{A}_v)]_j$ for $j\in\mathcal{S}^-$ and $[\bmu{\Phi}(\mathcal{A}_v)]_j=[\bmu{\Xi}(\mathcal{A}_v)]_j$ for $j\in\mathcal{S}^+$.}
	\label{fig:muDyAlimits2}
\end{figure}

\end{exa}

\section{Conclusions}\label{sec:conc}

In this paper we have provided the first theoretical step towards practical applications of the stochastic fluid-fluid models (SFFMs), noting that such applications rely on the ability to perform matrix computations for the quantities of the interest.

We have demonstrated how to apply matrix-analytic methods in modelling of the SFFMs when the intensities of the growth $r_i$, $i\in \mathcal{S}$, of the second level process $Y(t)$, do not depend on the level of first level process $X(t)$. In particular, we have developed results for the key quantity $\Psi$, which in the theory of stochastic fluid models is required for the computation of many other transient and stationary quantities of interest~\cite{BeaOreFF,BOT}, and we derived expressions for other useful quantities as well.

One can treat our main results as the first building block towards developing this method for more general SFFMs. As a simple example, we can consider piece-wise constant rates $r_i(x)$ with $r_i(x)=r_i^1$ if $x<L$ and $r_i(x)=r_i^2$ otherwise, where $r_i^1, r_i^2$ and $L$ are fixed rates and threshold, respectively. It would be very interesting to handle other than exponential initial distributions $\nu_i(x)$ as well, see the assumption made in~\eqref{muexpo}. One can think for example of a mixture of exponential distributions~\cite{2021YTT} for each state $i\in \mathcal{S}$. These topics are planed for further investigations of the transient and stationary analysis of the SFFMS.

\appendix

\section{Appendix}

{\bf Proof of Lemma~\ref{lemma_opB_bounded}.} Since we have $V(u+w)=V(u)V(w)$ for all $u,w>0$, and $V(0) = I$, it follows that $V(t)=e^{Bt}$, where $B=\lim_{t\to 0^{+}}\frac{ V(t)-I}{t}$, and $B=[B_{ij}]_{i,j\in\mathcal{S}}$, $B_{ij}=\lim_{t\to 0^{+}}\frac{ V_{ij}(t)-I_{ij}}{t}$, with
\begin{eqnarray*}
		\mu_i B_{ij}(\mathcal{A}_v)
	&=&
	\lim_{t\to 0^{+}}
	\Bigg\{
	\Bigg(
	\int_{x=0}^{\infty}d\mu_i(x)\mathbb{P}
	\left(
	\varphi(t)=j,X(t)\in\mathcal{A}_v|\varphi(0)=i,X(0)=x
	\right)
	\\
	&&\quad -\mu_i(\mathcal{A}_v)I(i=j)
	\Bigg)\Big/t \Bigg\}.
\end{eqnarray*}

Then for $i\not= j$, we consider a sample path contributing to $V_{ij}(t)$ in which the process starts from $\varphi(0)=i$, $X(0)=x$, then remains in phase $i$ until a transition from $i$ to $j$ at time $h$, for some $0<h<t$, and then remains in $j$ until time $t$. We note that other events contributing to $V_{ij}(t)$ occur with probability $o(t)$.

Therefore, for all $i\neq j$, by applying $\int_{h=0}^t g(h)dh= t\times g(t) + o(t)$ for small $t$, we obtain
\begin{eqnarray*}
	\mu_i B_{ij}(\mathcal{A}_v)	
		&=&
		\lim_{t\to 0^{+}}
		\Bigg\{
		\Bigg(
		\int_{x=0}^{\infty}d\mu_i(x)\mathbb{P}
		\left(
		\varphi(t)=j,X(t)\in\mathcal{A}_v\ |\ \varphi(0)=i,X(0)=x
		\right)
		\Bigg)\Big/t \Bigg\}
		\\
	&=&
	\lim_{t\to 0^{+}}
	\Bigg\{
	\Bigg(
	\int_{x=0}^{\infty} \nu_i(x)
	\mathbb{P}
	\left(
	\varphi(t)=j,X(t)\in\mathcal{A}_v\ |\ \varphi(0)=i,X(0)=x
	\right)dx
	\Bigg)\Big/t \Bigg\}
\\	
	&&	+\lim_{t\to 0^{+}}
	\Bigg\{
	\Bigg(
	p_i(0)
	\mathbb{P}
	\left(
	\varphi(t)\!=\!j,X(t)\!\in\!\mathcal{A}_v\ |\ \varphi(0)\!=\!i,X(0)\!=\!0
	\right)
	\Bigg)\Big/t \Bigg\}		
	\\
	&=&\lim_{t\to 0^{+}}
	\Bigg\{\Bigg(
	\int_{x=0}^{\infty} \!\!\! \nu_i(x) \!\!
	\int_{h=0}^t \!\!\!\!
	e^{T_{ii}h}T_{ij}e^{T_{jj}(t-h)}
	I(x\!+\!c_ih\!+\!c_j(t\!-\!h)\in\mathcal{A}_v)dhdx
	+ o(t)
	\Bigg)\Big/t
	\Bigg\}
	\\	
&&+
	p_i(0)
	\lim_{t\to 0^{+}}
	\Bigg\{\Bigg(
	\int_{h=0}^t
	e^{T_{ii}h}T_{ij}e^{T_{jj}(t-h)}
\Big[ I(c_j>0)I(c_j(t-h)\in\mathcal{A}_v)+I(c_j<0)
\Big]
dh + o(t)\Bigg)\Big/t
\Bigg\}
\\
	&=&\lim_{t\to 0^{+}}
	\Bigg\{\Bigg(
	\int_{x=0}^{\infty} \!\!\! \nu_i(x) \!\!
	\times
t\times
	e^{T_{ii}t}T_{ij}e^{T_{jj}(t-t)}
	I(x\!+\!c_it\!+\!c_j(t\!-\!t)\in\mathcal{A}_v)dx\Bigg)\Big/t
	\Bigg\}
	\\	
	&&+
	p_i(0)
	\lim_{t\to 0^{+}}
	\Bigg(
	t\times
	e^{T_{ii}t}T_{ij}e^{T_{jj}(t-t)}
	\Big[ I(c_j>0)I(c_j(t-t)\in\mathcal{A}_v)+I(c_j<0)
	\Big]
	\Big/t
	\Bigg),
\end{eqnarray*}		
and so,	
	\begin{eqnarray*}
		\mu_i B_{ij}(\mathcal{A}_v)	
		&=&\lim_{t\to 0^{+}}
	\int_{x=0}^{\infty}\nu_i(x)
	e^{T_{ii}t}T_{ij}I(x+c_it\in \mathcal{A}_v)dx
	+p_i(0)\lim_{t\to 0^{+}}e^{T_{ii}t}T_{ij}\\
	\\
	&=&T_{ij}
	\int_{x=0 }^v\nu_i(x)dx + T_{ij} p_i(0)
	\\
	&=&T_{ij}  \mu_i(\mathcal{A}_v).
\end{eqnarray*}

Further, for $j \in {\mathcal S}$, we consider a sample path contributing to $V_{jj}(t)$ in which the process starts from $\varphi(0)=j$, $X(0)=x$, and remains in phase $j$ until time $t$. We note that other events contributing to $V_{jj}(t)$ occur with probability $o(t)$.

Therefore, for all $j \in {\mathcal S}$,
\begin{eqnarray*}
	\mu_j B_{jj}(\mathcal{A})	
	&=&
	\lim_{t\to 0^{+}}
	\Bigg\{
	\Bigg(
	\int_{x=0}^{\infty}d\mu_i(x)\mathbb{P}
	\left(
	\varphi(t)=j,X(t)\in\mathcal{A}_v\ |\ \varphi(0)=j,X(0)=x
	\right)
	\\
	&&\quad -\mu_j(\mathcal{A}_v)
	\Bigg)\Big/t \Bigg\}	
	\\
	&=&
	\lim_{t\to 0^{+}}
	\Bigg\{
	\Bigg(
	\int_{x=0}^{\infty}
	\nu_j(x)
	\mathbb{P}[
	\varphi(t)=j,X(t)\in\mathcal{A}_v\ |\ \varphi(0)=j,X(0)=x]dx\\
	&&
	\hspace*{71pt}
	-\int_{x=0}^v\nu_j(x)dx
	\Bigg)\Big/t \Bigg\}
\\
	&&+
	\lim_{t\to 0^{+}}
	\Bigg\{
	\Bigg(
	p_j(0)
	\mathbb{P}
	[\varphi(t)=j,X(t)\in\mathcal{A}_v\ |\ \varphi(0)=j,X(0)=x]
	\\
	&&
	\hspace*{71pt}
	-
	p_j^m(0)
	\Bigg)\Big/t \Bigg\},
\end{eqnarray*}
where, by L'Hospital's rule, the first term in the above is equal to
\begin{eqnarray*}
	\lefteqn{
		\lim_{t\to 0^{+}}
		\Bigg\{
		\Bigg(
		\int_{x=0}^{\infty}
		\nu_j(x)
		e^{T_{jj}t}
		\Big[
		I(x+c_jt>0)I(x+c_jt\in\mathcal{A}_v) } \\
	&&
	+I(x+c_jt\leq 0)
	\Big]
	dx
	-\int_{x=0}^v\nu_j(x)dx
	+ o(t)
	\Bigg)
	\Big/t
	\Bigg\}
	\\
	&=&
	\lim_{t\to 0^{+}}
	\Bigg\{
	T_{jj}
	\int_{x=0}^{\infty}
	\nu_j(x)
	e^{T_{jj}t}
	\Big[
	I(x+c_jt>0)I(x+c_jt\in\mathcal{A}_v)
	+I(x+c_jt\leq 0)
	\Big]
	\Bigg\}
	\\
	&=&
	T_{jj}
	\int_{x=0}^{v}\nu_j(x)dx
	,
\end{eqnarray*}
and the second term is equal to
\begin{eqnarray*}
	p_j(0)
	\lim_{t\to 0^{+}}
	\Bigg\{
	\Bigg(
	e^{T_{jj}t}
	I\left(c_j t \in \mathcal{A}_v \right)
	-
	1
	+ o(t)
	\Bigg)
	\Big/t
	\Bigg\}
	&=&
	p_j(0)
	\lim_{t\to 0^{+}}
	\left\{
	T_{jj}
	e^{T_{jj}t}
 I\left(c_j t \in \mathcal{A}_v\right)
		\right\}
	\\
	&=&
	T_{jj}p_j(0),
\end{eqnarray*}

Consequently,
\begin{eqnarray*}
	\mu_j B_{jj}(\mathcal{A})&=&
	T_{jj}
	\int_{x=0}^{v}\nu_j(x)dx
	+T_{jj}p_j(0)
	=T_{jj}  \mu_j(\mathcal{A}_v),
\end{eqnarray*}	
and so the result follows. \hfill\rule{9pt}{9pt}
\bigskip

\noindent{\bf Proof of Lemma~\ref{lemma_D}.} Since we have $U(u+w,s)=U(u,s)U(w,s)$ for all $u,w>0$, and $U(0,s) = I$, it follows that $U(y,s)=e^{D(s)y}$, where $D(s)=\lim_{y\to 0^{+}}\frac{ U(y,s)-I}{t}$, and $D(s)=[D(s)_{ij}]_{i,j\in\mathcal{S}}$, $D(s)_{ij}=\lim_{y\to 0^{+}}\frac{ U_{ij}(y,s)-I_{ij}}{y}$, with
\begin{eqnarray*}
	\mu_i D_{ij}(s)(\mathcal{A}_v)
	&=&
	\lim_{y\to 0^{+}}
	\Bigg\{
	\Bigg(
	\int_{x=0}^{\infty}d\mu_i(x)\mathbb{P}
	\left(
	\varphi(\omega(y))=j,X(\omega(y))\in\mathcal{A}_v
	\ |\
	\varphi(0)=i,X(0)=x
	\right)
	\\
	&&\quad -\mu_i(\mathcal{A}_v)I(i=j)
	\Bigg)\Big/y \Bigg\}.
\end{eqnarray*}

Denote
\begin{eqnarray*}
	t^*(i,x,y)&=&\left(\
	\omega(y)
	\ |\
	\varphi(0)=i,X(0)=x,Y(0)=0,\varphi(u)=i,0<u\leq \omega(y)
	\
	\right),\\
	x^*(i,x,y)&=&\left(\
	X(\omega(y))
	\ |\
	\varphi(0)=i,X(0)=x,Y(0)=0,\varphi(u)=i,0<u\leq \omega(y)
	\
	\right),
\end{eqnarray*}
which are interpreted as random variables $\omega(y)$ and $X(\omega(y))$, respectively, given  that the process starts in $\varphi(0)=i$, $X(0)=x$, $Y(0)=0$, and there is no change in phase at least until time $\omega(y)$.

Then for $i\not= j$, we consider a sample path contributing to $U_{ij}(y,s)$ in which the process starts from $\varphi(0)=i$, $X(0)=x$, $Y(0)=0$, then remains in phase $i$ until a transition from $i$ to $j$ at time $\omega(u)$, for some $0<u<y$, and then remains in $j$ until time $\omega(y)$. We note that other events contributing to $U_{ij}(y,s)$ occur with probability $o(y)$.

Therefore, for $i\not= j$, by L'Hospital's rule, and since $t^{*}(i,x,0)=0$, $x^*(i,x,0)=x$, we have,
\begin{eqnarray*}
	\lefteqn{
		\mu_iD_{ij}(s)(\mathcal{A})=
		\lim_{y\to 0^{+}}
		\mu_i\frac{U_{ij}(y,s)-I_{ij}} {y}
		(\mathcal{A}_v)
	}
	\\
	&=&
	\lim_{y\to 0^{+}}
	\Bigg\{
	\Bigg(
	\int_{u=0}^y
	\int_{x=0}^{\infty}d\mu_i(x)
	e^{-st^{*}(i,x,u)}e^{B_{ii}t^{*}(i,x,u)} \\
	&\times&
	\lim_{w\to 0^{+}}\Big(e^{-st^{*}(i,x^*(i,x,u),w)}V_{ij}(t^{*}(i,x^*(i,x,u),w)) \Big/w\Big)
	\\
	&\times& e^{-st^{*}(j,x^*(i,x,u),y-u)}
	e^{B_{jj}t^{*}(j,x^*(i,x,u),y-u)}
	(x,\mathcal{A}_v)dtdu
	+o(y)
	\Bigg)
	\Big/y
	\Bigg\}
	\\
%
&=&
		\lim_{y\to 0^{+}}
		\Bigg\{
		\int_{x=0}^{\infty} d\mu_i(x)
		\lim_{w\to 0^{+}}\Big(e^{-st^{*}(i,x,w)}V_{ij}(t^{*}(i,x,w))\Big/w\Big)
	\\
	&& \times
	e^{-st^{*}(j,x,y)}e^{B_{jj}t^{*}(j,x,y)}
	(x,\mathcal{A}_v)
	\Bigg\}
	\\
	&=&
	\int_{x=0}^{\infty}d\mu_i(x)
	\lim_{w\to 0^{+}} \Bigg(\frac{\partial t^{*}(i,x,w)}{\partial w}
	\Bigg[-s e^{-st^{*}(i,x,w)} V_{ij}(t^{*}(i,x,w))
	\\
	&& + e^{-st^{*}(i,x,w)} \frac{\partial V_{ij}(t^{*}(i,x,w))}{\partial w}\Bigg] \Bigg)
	(x,\mathcal{A}_v)
	\\
	&=&
	\int_{x=0}^{\infty}d\mu_i(x) \frac{1}{|r_i|}
	\Big[0 +  B_{ij} \Big]
	(x,\mathcal{A}_v)
	\\
	&=&
		\frac{1}{|r_i|}
		\mu_i [B]_{ij}(\mathcal{A}_v).
\end{eqnarray*}

Further, for $j \in {\mathcal S}$, we consider a sample path contributing to $U_{jj}(y,s)$ in which the process starts from $\varphi(0)=j$, $X(0)=x$, $Y(0)=0$, and remains in phase $j$ until time $\omega(u)$. We note that other events contributing to $U_{jj}(y,s)$ occur with probability $o(y)$.

Therefore, for all $j \in {\mathcal S}$, by L'Hopital's rule,
\begin{eqnarray*}
	\lefteqn{
		\mu_jD_{jj}(s)(\mathcal{A}_v)=
		\lim_{y\to 0^{+}}
		\mu_j^+\frac{U_{jj}(y,s)-I_{jj}} {y}
		(\mathcal{A}_v)
	}
	\\
	&=&
	\lim_{y\to 0^{+}}
	\Bigg\{
	\Bigg(
	\int_{x=0}^{\infty}d\mu_j(x)
	e^{-st^*(j,x,y)}e^{B_{jj}t^*(j,x,y)}(x,\mathcal{A}_v)-\mu_j(\mathcal{A}_v)
	+o(y)\Bigg)
	\Big/y
	\Bigg\}
\\
&=&
		\lim_{y\to 0^{+}}
		\Bigg\{
		\int_{x=0}^{\infty}d\mu_j(x)
		\frac{\partial t^*(j,x,y)}{\partial y}
		(B_{jj}-sI_{jj})
		e^{(B_{jj}-sI_{jj})t^*(j,x,y)}
		(x,\mathcal{A}_v)
		\Bigg\}
	\\
	&=&
	\int_{x=0}^{\infty}d\mu_j(x)
	\frac{1}{|r_j|}
	(B_{jj}-sI_{jj})
	(x,\mathcal{A}_v)
	\\
	&=&
		\frac{1}{|r_j|}\mu_j\left[\left(B-sI\right)\right]_{jj}(\mathcal{A}_v),
\end{eqnarray*}
and so the result follows. \hfill\rule{9pt}{9pt}

\bibliographystyle{abbrv}
\bibliography{refs}


\end{document}